\documentclass{amsart}
\usepackage{graphicx,amsthm}
\usepackage{verbatim,amsmath,amscd,
amssymb,color}
\usepackage{enumerate}
\usepackage[all]{xy}
\xyoption{all}
\numberwithin{equation}{section}

\newcommand{\cB}{{\mathcal B}}

\newcommand{\cR}{{\mathcal R}}

\newcommand{\cV}{{\mathcal V}}

\newcommand{\frg}{\mathfrak g}

\newcommand{\frt}{\mathfrak t}

\newcommand{\bbC}{\mathbb C}

\newcommand{\bbR}{\mathbb R}

\newcommand{\bbX}{\mathbb X}
\newcommand{\bbZ}{\mathbb Z}

\newcommand{\wt}{\mathrm{wt}}

\newcommand{\lbr}{\begin{bmatrix}}
\newcommand{\rbr}{\end{bmatrix}}
\newcommand{\forb}{\bigcirc\kern-2.8ex \because}
\newcommand{\forbb}{\bigcirc\kern-3.0ex \because}
\newcommand{\forbbb}{\bigcirc\kern-3.1ex \because}
\newcommand{\ld}{\ldots}

\newcommand\C{\mathbb C}

\def\ge{\frg}

\def\al{\alpha}

\def\beneme{\begin{enumerate}}
\def\beq{\begin{equation}}
\def\beqn{\begin{eqnarray}}
\def\beqnn{\begin{eqnarray*}}

\def\bbra#1,#2,#3{\left\{\begin{array}{c}\hspace{-5pt}
#1;#2\\ \hspace{-5pt}#3\end{array}\hspace{-5pt}\right\}}
\def\cd{\cdots}
\def\CC{\mathbb{C}}

\def\del{\delta}
\def\Del{\Delta}

\def\Delre{\Delta^{\rm re}}

\def\eit{\tilde{e}_i}
\def\eneme{\end{enumerate}}

\def\eeq{\end{equation}}
\def\eeqn{\end{eqnarray}}
\def\eeqnn{\end{eqnarray*}}
\def\fit{\tilde{f}_i}

\def\gau#1,#2{\left[\begin{array}{c}\hspace{-5pt}#1\\
\hspace{-5pt}#2\end{array}\hspace{-5pt}\right]}

\def\ify{\infty}
\def\io{\iota}

\def\lan{\langle}
\def\lar{\longrightarrow}

\def\lm{\lambda}
\def\Lm{\Lambda}
\def\mapright#1{\smash{\mathop{\longrightarrow}\limits^{#1}}}
\def\mapleft#1{\smash{\mathop{\longleftarrow}\limits^{#1}}}

\def\nd{\noindent}
\def\nn{\nonumber}

\def\ot{\otimes}

\def\osigma{\ovl\sigma}
\def\ovl{\overline}

\def\qq{\qquad}
\def\q{\quad}
\def\qed{\hfill\framebox[2mm]{}}
\newcommand{\QED}{\end{proof}}
\newcommand{\Proof}{\begin{proof}}

\def\ra{\rightarrow}
\def\ran{\rangle}

\def\tropical{tropical }
\def\Tropical{Tropical }

\def\til{\tilde}
\def\tm{\times}
\def\tt{\frt}
\def\TY(#1,#2,#3){#1^{(#2)}_{#3}}
\def\TTY(#1,#2,#3){#1_{#3}}
\def\STY(#1,#2,#3){#1^{#2}_{#3}}

\def\uq{U_q(\ge)}
\def\uqp{U'_q(\ge)}

\def\uqm{U^-_q(\ge)}

\def\vep{\varepsilon}
\def\vp{\varphi}

\def\W1{W(\varpi_1)}

\def\wt{{\rm wt}}
\def\wtil{\widetilde}
\def\what{\widehat}

\def\ZZ{\mathbb Z}

\def\m@th{\mathsurround=0pt}

\def\fsquare(#1,#2){v_{#2}}

\newtheorem{thm}{Theorem}[section]
\newtheorem{pro}[thm]{Proposition}
\newtheorem{lem}[thm]{Lemma}
\newtheorem{ex}[thm]{Example}
\newtheorem{cor}[thm]{Corollary}

\theoremstyle{definition}
\newtheorem{df}[thm]{Definition}

\newcommand{\cmt}{\marginpar}
\newcommand{\seteq}{\mathbin{:=}}
\newcommand{\cl}{\colon}
\newcommand{\be}{\begin{enumerate}}
\newcommand{\ee}{\end{enumerate}}
\newcommand{\bnum}{\be[{\rm (i)}]}
\newcommand{\enum}{\ee}

\newcommand{\set}[2]{\left\{#1\,\vert\,#2\right\}}
\newcommand{\sbigoplus}{{\mbox{\small{$\bigoplus$}}}}
\newcommand{\ba}{\begin{array}}
\newcommand{\ea}{\end{array}}
\newcommand{\on}{\operatorname}
\newcommand{\eq}{\begin{eqnarray}}
\newcommand{\eneq}{\end{eqnarray}}
\newcommand{\hs}{\hspace*}

\begin{document}
\font\germ=eufm10
\def\bl{\bullet}
%\makeatletter
\def\aaa{@}
%\vskip3cm

%\title{\bf Tropical R 
%and Affine Geometric Crystals}
\title{Tropical R maps
and Affine Geometric Crystals}

\author{Masaki K\textsc{ashiwara}}
\address{Research Institute for 
Mathematical Sciences, 
Kyoto University, Kitashiwakawa, Sakyo-ku,
Kyoto 606, Japan}
\email{masaki{\aaa}kurims.kyoto-u.ac.jp}
\author{Toshiki N\textsc{akashima}}
\address{Department of Mathematics, 
Sophia University, Kioicho 7-1, Chiyoda-ku, Tokyo 102-8554,
Japan}
\email{toshiki@mm.sophia.ac.jp}
\thanks{supported in part by JSPS Grants 
in Aid for Scientific Research \#18340007(M.K.), 
\#19540050(T.N.), \#20540016(M.O.)}
\author{Masato O\textsc{kado}}
\address{Department of Mathematical 
Science,
Graduate School of Engineering Science, 
Osaka University, Toyonaka, 
Osaka 560-8531, Japan}
\email{okado{\aaa}sigmath.es.osaka-u.ac.jp}
%    General info
\date{}


\begin{abstract}
By modifying the method in \cite{KNO}, 
certain affine geometric crystals are 
realized in affinization
of the fundamental representation 
$W(\varpi_1)_l$ 
and the \tropical R maps for the 
affine geometric crystals are
described explicitly. 
We also define prehomogeneous geometric crystals and 
show that for a positive geometric crystal, 
the connectedness of the corresponding ultra-discretized crystal
is the sufficient condition for prehomogeneity of 
the positive geometric crystal.
Moreover, the uniqueness of 
 \tropical R maps is discussed.
\end{abstract}

\maketitle

\tableofcontents

%%%%%%%%% section 1 %%%%%%%%%%%

\section{\bf Introduction}

An R-matrix appears as a solution to 
the Yang-Baxter equation, which is a key
to solve integrable lattice models in 
statistical mechanics. To understand
R-matrices
representation theoretically, 
Drinfeld \cite{D} and Jimbo \cite{J} introduced the
quantized universal enveloping algebra. 
Once its representation theory 
is established, an R-matrix
is interpreted as an intertwiner 
between the tensor product of
finite-dimensional
modules and the one with the order of 
the tensor product being reversed. In
\cite{KMN1}, \cite{KMN2}
the notions of finite-dimensional modules 
and R-matrices acquired a
combinatorial
version by using the theory of crystal bases. 
In the papers, the objects,
perfect crystals and combinatorial R-matrices, 
are introduced and play an important
role to show that some physical quantities 
for particular vertex-type solvable
models are equal to affine Lie algebra characters. 
Combinatorial R-matrices also play
an essential role in some 
integrable cellular automata \cite{HKOTY}.

The notion of geometric crystals for 
semi-simple algebraic groups has been 
initiated by Berenstein and Kazhdan 
 \cite{BK}, and
it is extended for Kac-Moody groups in 
\cite{N}. Geometric crystals
are constructed on some geometric 
objects, such as algebraic varieties, and 
has an analogous structure to crystals.
The relation between crystals and 
geometric crystals is not only a
simple analogy, but also a more direct 
functorial connection, called 
ultra-discretization/tropicalization.
Indeed, in \cite{N}, the 
geometric crystals on 
Schubert varieties of
the corresponding Kac-Moody group
are ultra-discretized 
to certain crystals (see \ref{schubert}).
Applying this result, 
we constructed for an affine Lie algebra $\ge$
the affine geometric crystal 
$\cV=\cV(\ge)$ in
the fundamental $\ge$-representation $W(\varpi_1)$
in \cite{KNO} .
Ultra-discretizing the geometric crystal 
$\cV(\ge)$,
we obtained the limit of the 
coherent family of perfect crystals
$B_\ify(\ge^L)$ where $\ge^L$ is the 
Langlands dual of $\ge$ (see \cite{KKM}).

A tropical R map is an analogous 
object to the set-theoretic R 
(\cite{D2})
 and defined as follows
(see \S\ref{tro-r}):
For a 
family of geometric crystals 
$X\seteq\{X_\lm\}_{\lm\in\Lm}$
(parametrized by $\lm\in\Lm$) with a
product structure, a birational map
$R_{\lm\mu}:X_\lm\times X_\mu
\to X_\mu\times X_\lm$ $(\lm,\mu\in\Lm)$
is said to be a 
\tropical R map, if it satisfies the 
Yang-Baxter equation and preserves 
the geometric crystal structure.
In \cite{Y} and \cite{KOTY}, 
an explicit form of \tropical R map is given 
for the geometric crystal $\cB_l$
$(l\in \bbC^\times)$
of type $\TY(A,1,n)$ and $\TY(D,1,n)$.
They obtained the \tropical R map as a 
unique solution $(x',y')$
of the equation $M_l(x,z)M_m(y,z)
=M_m(x',z)M_l(y',z)$ where 
$x\in\cB_l,y\in\cB_m$, $z$ is an 
indeterminate and $M_l(x,z)$ is 
a square matrix called $M$-{\it matrix}
(see \S\ref{Mmat-sec}).

In this paper, we construct the affine 
geometric crystal $\cV_l$ in 
$W(\varpi_1)_l$, the affinization of the fundamental 
representation $W(\varpi_1)$, for $\ge=
\TY(A,1,n),\TY(B,1,n),\TY(D,1,n),$
$\TY(D,2,n+1),\TY(A,2,2n-1)$ and 
$\TY(A,2,2n)$. 
The geometric crystals $\cB_l$ constructed in 
\cite{KOTY} are isomorphic to our geometric crystals $\cV_l$. 
Thus, in these two cases 
we have the \tropical R map for $\cV_l$ 
through the isomorphism.
By virtue of the method 'folding', we obtain 
the \tropical R maps for the other cases; 
$\TY(B,1,n), \TY(D,2,n+1),\TY(A,2,2n-1)$ and 
$\TY(A,2,2n)$ from the \tropical R map for 
$\TY(D,1,N)$ with a suitable integer $N$.
(The case $\TY(C,1,n)$ can
not be obtained from the folding of 
$\TY(D,1,N)$. This case will be
discussed elsewhere.)
Here note that the geometric crystal $\cV(\ge)$ 
is a special case of $\cV(\ge)_l$ for $l=1$.

Let us explain the construction of the 
\tropical R map more precisely.
First, we take a Dynkin diagram automorphism
$\sigma$ for $\TY(D,1,N)$, which induces the 
automorphism $\Sigma$
of the geometric crystal $\cB_l$ of type $\TY(D,1,N)$, 
where $\cB_l\cong (\bbC^\times)^{2N-2}$ as 
algebraic varieties.
Let $X_l$ be the fixed-point variety for 
$\Sigma$, that is, $X_l\seteq
\{x\in\cB_l\,|\,\Sigma(x)=x\}$.
There is an invertible matrix $J$
such that $M_l(\Sigma(x),z)=JM_l(x,z)J^{-1}$
($x\in \cB_l$).
By this formula and the uniqueness of the 
solution for the equation $M_l(x,z)M_m(y,z)
=M_m(x',z)M_l(y',z)$ $(l,m\in\bbC^\times)$, 
we deduce that the \tropical R map
sends $X_l\times X_m\to X_m\times X_l$.
Furthermore, this fixed-point variety 
$X_l$ is 
equipped with the $\ge^\sigma$-geometric 
crystal structure, 
where $\ge^\sigma$ is the affine 
Lie algebra obtained from $\TY(D,1,N)$ 
by the folding associated 
with $\sigma$ and it 
is isomorphic to $\cV_l$ for $\ge^\sigma$.
Hence, we get the \tropical R map for 
$\{\cV_l\}$ in the case of $\ge^\sigma$.

In the article, we also discuss the 
uniqueness of the \tropical R maps.
As for usual R-matrices, it is unique (up to constant) by Schur's lemma,
if the corresponding modules are irreducible.
Similar situation occurs in the 
geometric crystals. 
A geometric crystal is {\it prehomogeneous}  if there exists
an open dense orbit by 
the actions of the $e_i^c$'s.
For $\phi,\phi'$ two
isomorphisms of prehomogeneous
geometric crystals, if there exists a point $p$
in the open dense orbit
which is sent to 
the same point by $\phi$ and $\phi'$,
it is shown that $\phi=\phi'$ as rational 
maps.
A crucial result in this article is that 
a positive 
geometric crystal $\bbX$ is prehomogeneous if 
the crystal ultra-discretized from 
$\bbX$ is connected (see Theorem \ref{conn-preh}).
This uniqueness implies the following fact:
For geometric crystals $X,Y,Z$, suppose that 
the product $X\times Y\times Z$ is prehomogeneous
and that there exist \tropical R maps for any 
two of $X,Y,Z$, then it follows from the uniqueness that we have 
the Yang-Baxter equation:
\[
 R^{(12)}R^{(23)}R^{(12)}=R^{(23)}R^{(12)}R^{(23)}.
\]
If $X\times Y$ is prehomogeneous, we also obtain
the inversion formula:
\[
 R_{YX}R_{XY}={\rm id},
\]
where $R^{(ij)}$ means $R$ acting on the $i$-th and the
$j$-th components of the product.
Namely, important properties of the \tropical R 
maps are deduced from the uniqueness.
The \tropical R map 
on the affine geometric crystal $\cV_l$ 
introduced in this article
is unique, since the ultra-discretized
crystal ${\mathcal UD}(\cV_l)$ 
and its tensor products are connected,
which implies that $\cV_l$ and its products
are prehomogeneous.

%%%%%%%%%%%%%%% Section 2 %%%%%%%%%%%%%%%%
\section{\bf Geometric Crystals and Crystals}
\label{sec2}

\nd
In this section, 
we review Kac-Moody groups and geometric crystals
following 
\cite{BK}, \cite{Kac}, \cite{N}, 
\cite{N2}, \cite{PK}.
\subsection{Kac-Moody algebras and Kac-Moody groups}
\label{KM}
Fix a symmetrizable generalized Cartan matrix
 $A=(a_{ij})_{i,j\in I}$ with a finite index set $I$.
Let $(\tt,\{\al_i\}_{i\in I},\{\al^\vee_i\}_{i\in I})$ 
be the associated
root data, where ${\tt}$ is a vector space 
over $\bbC$ and
$\{\al_i\}_{i\in I}\subset\tt^*$ and 
$\{\al^\vee_i\}_{i\in I}\subset\tt$
are linearly independent 
satisfying $\al_j(\al^\vee_i)=a_{ij}$.

The Kac-Moody Lie algebra $\ge=\ge(A)$ associated with $A$
is the Lie algebra over $\bbC$ generated by $\tt$, the 
Chevalley generators $e_i$ and $f_i$ $(i\in I)$
with the usual defining relations (\cite{Kac}).
There is the root space decomposition 
$\ge=\sbigoplus_{\al\in \tt^*}\ge_{\al}$.
Denote the set of roots by 
$\Delta\seteq\set{\al\in \tt^*}{\al\ne0,\,\,\ge_{\al}\ne(0)}$.
Set $Q=\sum_i\bbZ \al_i$, $Q_+=\sum_i\bbZ_{\geq0} \al_i$,
$Q^\vee=\sum_i\bbZ \al^\vee_i$
and $\Delta_+=\Delta\cap Q_+$.
An element of $\Delta_+$ is called 
a {\it positive root}.
Let $P\subset \tt^*$ be a weight lattice such that 
$\bbC\ot P=\tt^*$ and 
$Q\subset P\subset 
\{\lm\,|\, \lm(Q^\vee)\subset \bbZ\}$, 
whose element is called a weight.

Define the simple reflections $s_i\in{\rm Aut}(\tt)$ $(i\in I)$ by
$s_i(h)\seteq h-\al_i(h)\al^\vee_i$, which generate the Weyl group $W$.
It induces the action of $W$ on $\tt^*$ by
$s_i(\lm)\seteq\lm-\lm(\al^\vee_i)\al_i$.
Set $\Delre\seteq\set{w(\al_i)}{w\in W,\,\,i\in I}$, whose element 
is called a real root.

Let $G$ be the Kac-Moody group associated 
with $(\ge,P)$ (\cite{PK}).
Let $U_{\al}\seteq\exp\ge_{\al}$ $(\al\in \Delre)$
be the one-parameter subgroup of $G$.
The group $G$ is generated by $U_{\al}$ $(\al\in \Delre)$.
Let $U^{\pm}$ be the subgroup generated by $U_{\pm\al}$
($\al\in \Delre_+=\Delre\cap Q_+$), {\it i.e.,}
$U^{\pm}\seteq\lan U_{\pm\al}|\al\in\Del^{\rm re}_+\ran$.

For any $i\in I$, there exists a unique group homomorphism
$\phi_i\cl SL_2(\bbC)\rightarrow G$ such that
\[
\phi_i\left(
\left(
\begin{array}{cc}
1&t\\
0&1
\end{array}
\right)\right)=\exp(t e_i),\,\,
 \phi_i\left(
\left(
\begin{array}{cc}
1&0\\
t&1
\end{array}
\right)\right)=\exp(t f_i)\qquad(t\in\bbC).
\]
Set $\al^\vee_i(c)\seteq
\phi_i\left(\left(
\begin{smallmatrix}
c&0\\
0&c^{-1}\end{smallmatrix}\right)\right)$,
$x_i(t)\seteq\exp{(t e_i)}$, $y_i(t)\seteq\exp{(t f_i)}$, 
$G_i\seteq\phi_i(SL_2(\bbC))$,
$T_i\seteq \alpha_i^\vee(\bbC^\times)$ 
and 
$N_i\seteq N_{G_i}(T_i)$. Let
$T$ be the subgroup of $G$ 
with $P$ as its weight lattice 
which is called a {\it maximal torus} in $G$, and let
$B^{\pm}=U^{\pm}T$ be the Borel subgroup of $G$.
Let $N$ be the subgroup of $G$ generated by
the $N_i$'s. Then we have the isomorphism
$\phi\cl W\mapright{\sim}N/T$ defined by $\phi(s_i)=N_iT/T$.
An element $\ovl s_i\seteq x_i(-1)y_i(1)x_i(-1)
=\phi_i\left(
\left(\begin{smallmatrix}
0&-1\\
1&0
\end{smallmatrix}
\right)\right)$ in 
$N_G(T)$ is a representative of 
$s_i\in W=N_G(T)/T$. 

%%%%%%%%%%%%%%%%%%%%%%%%%%%%%%%%%%%%%%%%%%%%%
\subsection{Geometric crystals}\label{sec:gc}
Let $W$ be the  Weyl group associated with $\ge$. 
Define $R(w)$ for $w\in W$ by
\[
 R(w)\seteq \set{(i_1,i_2,\ld,i_l)\in I^l}{w=s_{i_1}s_{i_2}\cd s_{i_l}},
\]
where $l$ is the length of $w$, {\em i.e.},
$R(w)$ is the set of reduced expressions of $w$.

Let $X$ be a variety, 
{$\gamma_i\cl X\rightarrow \bbC$} and 
$\vep_i\cl X\longrightarrow \bbC$ ($i\in I$) 
rational functions on $X$, and
{$e_i\cl\bbC^\times \times X\longrightarrow X$}
$((c,x)\mapsto e^c_i(x))$ a
rational $\bbC^\times$-action.
For $w\in W$ and ${\bf i}=(i_1,\ld,i_l)\in R(w)$, set 
$\al^{(j)}\seteq s_{i_l}\cd s_{i_{j+1}}(\al_{i_j})$ 
$(1\leq j\leq l)$ and 
\begin{eqnarray*}
e_{\bf i}\cl T\times X&\rightarrow &X\\
(t,x)&\mapsto &e_{\bf i}^t(x)\seteq e_{i_1}^{\al^{(1)}(t)}
e_{i_2}^{\al^{(2)}(t)}\cd e_{i_l}^{\al^{(l)}(t)}(x).
\label{tx}
\end{eqnarray*}
\begin{df}
\label{def-gc}
A quadruple $(X,\{e_i\}_{i\in I},\{\gamma_i,\}_{i\in I},
\{\vep_i\}_{i\in I})$ is a 
$G$ (or $\ge$)-{\it geometric crystal} 
if
\bnum
\item
$\{1\}\times X\cap dom(e_i)$ 
is open dense in $\{1\}\times X$ for any $i\in I$.
Here $dom(e_i)$ is the domain of definition of
$e_i\cl\C^\times\times X\to X$.
\cmt{changed}
\item
$\gamma_j(e^c_i(x))=c^{a_{ij}}\gamma_j(x)$.
\item
{$e_{\bf i}=e_{\bf i'}$}
for any 
$w\in W$ and ${\bf i}$.
${\bf i'}\in R(w)$.
\item
$\vep_i(e_i^c(x))=c^{-1}\vep_i(x)$.
\ee
\end{df}
Note that the condition (iii) is 
equivalent to the following so-called 
{\it Verma relations}:
\[
 \begin{array}{lll}
&\hspace{-20pt}e^{c_1}_{i}e^{c_2}_{j}
=e^{c_2}_{j}e^{c_1}_{i}&
{\rm if }\,\,a_{ij}=a_{ji}=0,\\
&\hspace{-20pt} e^{c_1}_{i}e^{c_1c_2}_{j}e^{c_2}_{i}
=e^{c_2}_{j}e^{c_1c_2}_{i}e^{c_1}_{j}&
{\rm if }\,\,a_{ij}=a_{ji}=-1,\\
&\hspace{-20pt}
e^{c_1}_{i}e^{c^2_1c_2}_{j}e^{c_1c_2}_{i}e^{c_2}_{j}
=e^{c_2}_{j}e^{c_1c_2}_{i}e^{c^2_1c_2}_{j}e^{c_1}_{i}&
{\rm if }\,\,a_{ij}=-2,\,
a_{ji}=-1,\\
&\hspace{-20pt}
e^{c_1}_{i}e^{c^3_1c_2}_{j}e^{c^2_1c_2}_{i}
e^{c^3_1c^2_2}_{j}e^{c_1c_2}_{i}e^{c_2}_{j}
=e^{c_2}_{j}e^{c_1c_2}_{i}e^{c^3_1c^2_2}_{j}e^{c^2_1c_2}_{i}
e^{c^3_1c_2}_je^{c_1}_i&
{\rm if }\,\,a_{ij}=-3,\,
a_{ji}=-1,
\end{array}
\]
Note that the last formula is different from the one in 
\cite{BK}, \cite{N}, \cite{N2} which seems to be
incorrect. %The formula here is correct.
If $\bbX=(X,\{e_i\},\{\gamma_i\},\{\vep_i\})$
satisfies the conditions (i), (ii) and (iv), 
we call $\bbX$ a {\it pre-geometric crystal}.
%%%%%%%%%%%%%%%%%%%%%%%%%%%%%%%%%%%%
\subsection{Geometric crystal on Schubert cell}
\label{schubert}

Let $X\seteq G/B$ be the flag 
variety, which is the inductive limit 
of finite-dimensional projective varieties.
For $w\in W$, let $X_w\seteq BwB/B\subset X$ be the
Schubert cell associated with $w$, which has 
a natural geometric crystal structure
(\cite{BK}, \cite{N}).
For ${\bf i}=(i_1,\ld,i_k)\in R(w)$, set 
\begin{equation}
B_{\bf i}^-
\seteq \{Y_{\bf i}(c_1,\ld,c_k)
\seteq Y_{i_1}(c_1)\cd Y_{i_l}(c_k)
\,\vert\, c_1\cd,c_k\in\bbC^\times\}\subset B^-
\label{bw1}
\end{equation}
where $Y_i(c)\seteq y_i(c^{-1})\al_i^\vee(c)
=\al_i^\vee(c)y_i(c)$.
Then $B_{\bf i}^-$ 
is birationally isomorphic to $X_w$
and endowed with the induced geometric crystal structure.
Let $\xi$ be an element in the torus $T$. 
Then we can define the geometric 
crystal structure on $B_{\bf i}^-\cdot\xi$ and 
we shall describe 
its explicit form: The action $e^c_i$ on 
$B_{\bf i}^-\cdot\xi$ is given by
\[
e_i^c(Y_{i_1}(c_1)\cd Y_{i_l}(c_k)\xi)
=Y_{i_1}({\mathcal C}_1)\cd Y_{i_l}({\mathcal C}_k)\xi,
\]
where
\begin{equation}
{\mathcal C}_j\seteq 
c_j\cdot \frac{\displaystyle \sum_{1\leq m\leq j,\,i_m=i}
 \frac{c}
{c_1^{a_{i_1,i}}\cd c_{m-1}^{a_{i_{m-1},i}}c_m}
+\sum_{j< m\leq k,\,i_m=i} \frac{1}
{c_1^{a_{i_1,i}}\cd c_{m-1}^{a_{i_{m-1},i}}c_m}}
{\displaystyle\sum_{1\leq m<j,\,i_m=i} 
 \frac{c}
{c_1^{a_{i_1,i}}\cd c_{m-1}^{a_{i_{m-1},i}}c_m}+
\mathop\sum_{j\leq m\leq k,\,i_m=i}  \frac{1}
{c_1^{a_{i_1,i}}\cd c_{m-1}^{a_{i_{m-1},i}}c_m}}.
\label{eici}
\end{equation}
The explicit forms of 
rational functions $\vep_i$ and $\gamma_i$ are:
\[
 \vep_i(Y_{i_1}(c_1)\cd Y_{i_l}(c_k)\xi)=
\hs{-3ex}\sum_{1\leq m\leq k,\,i_m=i} \frac{1}
{c_1^{a_{i_1,i}}\cd c_{m-1}^{a_{i_{m-1},i}}c_m},\q
\gamma_i(Y_{i_1}(c_1)\cd Y_{i_l}(c_k)\xi)
=c_1^{a_{i_1,i}}\cd c_k^{a_{i_k,i}}\al_i(\xi).
\]
These will be needed in Sect.5 
to construct affine geometric crystals.
%%%%%%%%%%%%%%%%%%%%%%%%%%%%%%%%%%
\subsection{Crystals}

We recall the notion of crystals,
which is obtained by
abstracting the combinatorial 
properties of crystal bases.
\begin{df}
A {\it crystal} $B$ is a set endowed with the following maps:
\begin{eqnarray*}
&& {\rm wt}\cl B\lar P,\\
&&\vep_i\cl B\lar\ZZ\sqcup\{-\infty\},\q
  \vp_i\cl B\lar\ZZ\sqcup\{-\infty\} \q{\hbox{for}}\q i\in I,\\
&&\eit\cl B\sqcup\{0\}\lar B\sqcup\{0\},
\q\fit\cl B\sqcup\{0\}\lar B\sqcup\{0\}\q{\hbox{for}}\q i\in I,\\
&&\eit(0)=\fit(0)=0.
\end{eqnarray*}
Those maps satisfy the following axioms: for
 all $b,b_1,b_2 \in B$, we have
\begin{eqnarray*}
&&\vp_i(b)=\vep_i(b)+\lan \al^\vee_i,{\rm wt}
(b)\ran,\\
&&\wt(\eit b)=\wt(b)+\al_i{\hbox{ if  }}\eit b\in B,\\
&&\wt(\fit b)=\wt(b)-\al_i{\hbox{ if  }}\fit b\in B,\\
&&\eit b_2=b_1 \Longleftrightarrow \fit b_1=b_2,\\
&&\vep_i(b)=-\ify
   \Longrightarrow \eit b=\fit b=0.
\end{eqnarray*}
\end{df}
\begin{ex}
\label{ex-tlm}
\bnum
\item
If $(L,B)$ is a crystal base, then $B$ is a crystal.
\item For the crystal base $(L(\ify),B(\ify))$
of the subalgebra $\uqm$ of 
the quantum algebra $\uq$, 
$B(\ify)$ is a crystal. \label{exam2}
\item
\label{tlm}
For $\lm\in P$, set $T_\lm\seteq \{t_\lm\}$. We define a crystal
structure on $T_\lm$ by 
\[
 \eit(t_\lm)=\fit(t_\lm)=0,\q\vep_i(t_\lm)=
\vp_i(t_\lm)=-\ify,\q \wt(t_\lm)=\lm.
\]
\ee
\end{ex}
\begin{df}
\begin{enumerate}
\item
To a crystal $B$, a colored oriented graph
is associated by 
\[
 b_1\mapright{i}b_2\Longleftrightarrow 
\fit b_1=b_2.
\]
We call this graph the {\em crystal graph}
of $B$. 
\item
A crystal $B$ is said to be {\it connected}, 
if its crystal graph is connected as a graph.
\item
A crystal $B$ is {\it free} 
if $\eit^n(b)\ne0$ and $\fit^n(b)\ne0$
for any $b\in B$, $i\in I$ 
and $n>0$.
\end{enumerate}
\end{df}
%%%%%%%%%%%%%%%%%%%%%%%%%%%%%

%%%%%%%%%%%%%%%%%%%%%%%%%%%%%%%
\subsection{Positive structure,
Ultra-discretization and Tropicalization}
\label{positive-str}

Let us recall the notions of 
positive structure and ultra-discretization/tropicalization.

Set $R\seteq \bbC(c)$ and define
$$
\begin{array}{cccc}
v\cl&R\setminus\{0\}&\longrightarrow &\ZZ\\
&f(c)&\mapsto
&{\rm deg}(f(c)).
\end{array}
$$
Here ${\rm deg}$ is the degree of poles at $c=\infty$.
%For any $\phi\in L(T)$, 
%set 
%${\rm deg}_T(\phi)\seteq v\circ\phi|_{X^*(T)}$. Since 
Note that for $f_1,f_2\in R\setminus\{0\}$, we have
\begin{equation}
v(f_1 f_2)=v(f_1)+v(f_2),\q
v\left(\frac{f_1}{f_2}\right)=v(f_1)-v(f_2).
\label{ff=f+f}
\end{equation}
We say that a non-zero rational function $f(c)\in \bbC(c)$ is 
{\it positive} if $f$ can be expressed as a ratio
of polynomials with positive coefficients.
Note that $f\in\C(c)$ is positive if and only if
any pole of $f$ is not a positive number and
$f(x)>0$ for any $x>0$.

If $f_1,\,\,f_2\in R$ 
are positive, then we have 
\begin{equation}
v(f_1+f_2)={\rm max}(v(f_1),v(f_2)).
\label{+max}
\end{equation}

Let $T\simeq(\bbC^\times)^l$ be an algebraic torus over $\bbC$ and 
$X^*(T)\seteq {\rm Hom}(T,\bbC^\times)$ 
(resp.~$X_*(T)\seteq {\rm Hom}(\bbC^\times,T)$) 
be the lattice of characters
(resp.~co-characters)
of $T$. 
We denote by $T_+$ the set of points $x$ in $T$
such that
$\chi(x)>0$ for any character $\chi$.
Then $\bigl((\C^\times)^n\bigr)_+=(\bbR_{>0})^n$.

A non-zero rational function on
an algebraic torus $T$ is called {\em positive} if
it is written as $g/h$ where
$g$ and $h$ are a positive linear combination of
characters of $T$.
\begin{df}
Let 
$f\cl T\rightarrow T'$ be 
a rational mapping between
two algebraic tori $T$ and 
$T'$.
We say that $f$ is {\em positive},
if $\bbX\circ f$ is positive
for any character $\bbX\cl T'\to \C$.
%for each $f_i$, there exist 
%polynomials
%$g_i(x_1,\ld,x_m)$, $h_i(x_1,\ld,x_m)$
%with positive coefficients such that $f_i=g_i/h_i
%\ne0$.
\end{df}
Denote by ${\rm Mor}^+(T,T')$ the set of 
positive rational mappings from $T$ to $T'$.

Note that any $f\in {\rm Mor}^+(T,T')$ 
induces a real analytic map
$f_+\cl T_+\to T'_+$.

\begin{lem}[\cite{BK}]
\label{TTT}
For any $f\in {\rm Mor}^+(T_1,T_2)$             
and $g\in {\rm Mor}^+(T_2,T_3)$, 
the composition $g\circ f$
is well-defined and belongs to ${\rm Mor}^+(T_1,T_3)$.
\end{lem}

By Lemma \ref{TTT}, we can define a category ${\mathcal T}_+$
whose objects are algebraic tori over $\bbC$ and arrows
are positive rational mappings.
The category ${\mathcal T}_+$ admits products.
For two algebraic tori $T$ and $T'$, their product in
${\mathcal T}_+$ coincides with 
the usual product of $T$ and $T'$.

Note that $T\mapsto T_+$ gives a functor from ${\mathcal T}_+$
to the category of real analytic manifolds.

Let $f\cl T\rightarrow T'$ be a 
positive rational mapping
of algebraic tori $T$ and 
$T'$.
We define a map $\what f\cl X_*(T)\rightarrow X_*(T')$ by 
\[
\langle\chi,\what f(\xi)\rangle
=v(\chi\circ f\circ \xi),
\]
where $\chi\in X^*(T')$ and $\xi\in X_*(T)$.
Note that $\chi\circ f\circ \xi$ is a 
rational map from $\bbC^\times$ to itself.
%%%%%%%%%%%%%%%%%%%%%%%%%%%%%%%%%%%%%%%%%%%%%%%%%%%%%

\begin{lem}[\cite{BK}]
For any algebraic tori $T_1$, $T_2$, $T_3$, 
and positive rational mappings
$f\in {\rm Mor}^+(T_1,T_2)$, 
$g\in {\rm Mor}^+(T_2,T_3)$, we have
$\what{g\circ f}=\what g\circ\what f.$
\end{lem}
%{\sl Proof.}
%For $\mu\in X_*(T_1)$, 
%Let ${\mathcal B}$ be the category of free $\ZZ$-modules,
%whose arrows are piece-wise linear maps.
By this lemma, we obtain a functor 
\[
\begin{array}{cccc}
{\mathcal UD}\cl&{\mathcal T}_+&\longrightarrow &{{\mathrm{Set}}}\\[2ex]
&T&\mapsto& X_*(T)\\
&(f\cl T\rightarrow T')&\mapsto& 
(\what f\cl X_*(T)\rightarrow X_*(T')).
\end{array}
\]

Let us come back to the situation in \S \ref{sec:gc}.
Hence $G$ is a Kac-Moody group and $T$ is its Cartan subgroup.
\begin{df}[\cite{BK}]
Let $\bbX=(X,\{e_i\}_{i\in I},\{\gamma_i\}_{i\in I},
\{\vep_i\}_{i\in I})$ be a 
%\cmt{$wt_i$ or $\gamma_i$?}
$G$ (or $\ge$)-geometric crystal, $T'$ an algebraic torus
and $\theta\cl T'\rightarrow X$ 
a birational mapping.
The mapping $\theta$ is called 
a {\em positive structure} on
$\bbX$ if it satisfies
\bnum
\item for any $i\in I$ the rational functions
$\gamma_i\circ \theta\cl T'\rightarrow \bbC$ and 
$\vep_i\circ \theta\cl T'\rightarrow \bbC$ 
are positive,
\item
for any $i\in I$, the rational mapping
$e_{i,\theta}\cl\bbC^\tm \tm T'\rightarrow T'$ defined by
$e_{i,\theta}(c,t)
\seteq \theta^{-1}\circ e_i^c\circ \theta(t)$
is positive.
\end{enumerate}
\end{df}
Let $\theta\cl T'\rightarrow X$ be a positive structure on 
a geometric crystal $\bbX=(X,\{e_i\}_{i\in I},\{\gamma_i\}_{i\in I},
\{\vep_i\}_{i\in I})$.
Applying the functor ${\mathcal UD}$ 
to the positive rational mappings
$e_{i,\theta}\cl\bbC^\tm \tm T'\rightarrow T'$ and
$\gamma_i\circ \theta,\,\,
\vep_i\circ\theta\cl T'\ra \bbC^\times$,
%(the notations are as above), 
we obtain
\begin{eqnarray*}
&&\til e_i\seteq {\mathcal UD}(e_{i,\theta})\cl
\ZZ\tm X_*(T') \rightarrow X_*(T')\\
&&{\rm wt}_i\seteq 
{\mathcal UD}(\gamma_i\circ\theta),\,\,
\vep_i\seteq
{\mathcal UD}(\vep_i\circ\theta)\cl
X_*(T')\rightarrow \bbZ.
\end{eqnarray*}
Hence
the quadruple $(X_*(T'),\{\til e_i\}_{i\in I},
\{{\rm wt}_i\}_{i\in I},\{\vep_i\}_{i\in I})$
is a free pre-crystal structure (see \cite[2.2]{BK}) 
and we denote it by ${\mathcal UD}_{\theta}(\bbX)=
{\mathcal UD}_{\theta,T'}(\bbX)$.
As for the definition of crystal, 
see \ref{limit} or \cite{KKM},\cite{K0},
\cite{K1}.
We have thus the following theorem:

\begin{thm}[\cite{BK}\cite{N}]
For any geometric crystal 
$\bbX=(X,\{e_i\}_{i\in I},\{\gamma_i\}_{i\in I},
\{\vep_i\}_{i\in I})$ and a positive structure
$\theta\cl T'\rightarrow X$, the associated pre-crystal 
${\mathcal UD}_{\theta,T'}(\bbX)$
%(X_*(T'),\{e_i\}_{i\in I},\{{\rm wt}_i\}_{i\in I},
%\{\vep_i\}_{i\in I})$ 
is a crystal {\rm (see \cite[2.2]{BK})}.
\end{thm}

Now, let ${\mathcal GC}^+(\ge)$ be the category whose 
object is a triplet
$(\bbX,T',\theta)$ where 
$\bbX=(X,\{e_i\},\{\gamma_i\},\{\vep_i\})$ 
is a $\ge$-geometric crystal and $\theta\cl T'\rightarrow X$ 
is a positive structure on $\bbX$, and morphism
$f\cl(\bbX_1,T'_1,\theta_1)$ $\rightarrow 
(\bbX_2,T'_2,\theta_2)$ is given by a morphism 
$\vp\cl X_1\longrightarrow X_2$  of geometric crystals
such that 
\[
f\seteq \theta_2^{-1}\circ\vp\circ\theta_1\cl T'_1\longrightarrow T'_2,
\]
is a positive rational mapping. Let ${\mathcal CR}(\ge)$
be the category of $\ge$-crystals. 
Then by the theorem above, we have
\begin{cor}
$\mathcal UD$ defines a functor
\begin{eqnarray*}
&&\ba{ccc}
 {\mathcal UD}:{\mathcal GC}^+(\ge)
&\longrightarrow &{\mathcal CR}(\ge^L)\\[2ex]
(\bbX,T',\theta)&\mapsto& X_*(T')\\
(f\cl(\bbX_1,T'_1,\theta_1)\rightarrow 
(\bbX_2,T'_2,\theta_2))&\mapsto&
(\what f\cl X_*(T'_1)\rightarrow X_*(T'_2)),
\ea
\end{eqnarray*}
where $\ge^L$ is the Langlands dual for $\ge$.
\end{cor}
%\cmt{Caution}
We call the functor $\mathcal UD$
{\it ``ultra-discretization''} as in \cite{N},\cite{N2}.
While 
for a crystal $B\in {\mathcal CR}(\ge)$, if there
exists an object $(\bbX,T',\theta)$ in
${\mathcal GC}^+(\ge^L)$, we cal $(\bbX,\theta)$ a 
{\it tropicalization }
of $B$.
%%%%%%%%%%%%%%%%%%%%%%%%%%%%%%%%%%%%%
%%%%%%%%%%% Section 3
\section{\bf Prehomogeneous
Geometric Crystals}\label{sec3}
\setcounter{equation}{0}
\renewcommand{\theequation}{\thesection.\arabic{equation}}

%%%%%%%%%%%%%%%%%%%%%%%%%%%%%%%%%%%%
\begin{df} Let $\bbX=(X,\{e^c_i\},
\{\gamma_i\},\{\vep_i\})$ be 
a geometric crystal. We say that 
$\bbX$ is {\it prehomogeneous} 
if there exists
a Zariski open dense subset $\Omega\subset X$
which is an orbit by the actions of the
$e_i^c$'s.
\end{df}

The following lemma is obvious.
\begin{lem}
\label{uniq}
Let $\bbX_j=(X_j,\{e^c_i\},
\{\gamma_i\},\{\vep_i\})$ $(j=1,2)$
be  geometric crystals and $\bbX_1$ prehomogeneous.
Let 
$\Omega_1\subset X_1$ be an open dense 
orbit in $X_1$. 
For isomorphisms of geometric
crystals $\phi$, $\phi'\cl\bbX_1\to\bbX_2$, 
suppose that there exists $p_1\in \Omega_1$ 
such that 
$\phi(p_1)=\phi'(p_1)\in X_2$.
Then, we have $\phi=\phi'$ as rational maps.
\end{lem}

The following is the criterion for the prehomogeneity of 
a geometric crystal.
\begin{thm}
\label{conn-preh}
Let $\bbX=(X,\{e^c_i\},
\{\gamma_i\},\{\vep_i\})$ be 
a finite-dimensional positive 
geometric crystal with 
the positive structure $\theta\cl T\to X$ and 
$B\seteq
 {\mathcal UD}_\theta(\bbX)$  the 
crystal obtained as the ultra-discretization
of $\bbX$. 
If $B$ is a connected crystal, then 
$\bbX$ is prehomogeneous.
\end{thm}
\Proof 
Set $m\seteq \dim T=\dim X$. 
We identify $T$ and $X$ by $\theta$, 
and take ${\bf a}\in T_+$
(see \ref{positive-str}).
Assume that $B$ is connected and $\bbX(=T)$ 
is not prehomogeneous.
Then there exists a nowhere dense closed subset $Z$ of $T$
such that $Z$ contains ${\bf a}$ and
is invariant by the actions of the $e_i^c$'s.

Suppose that $Z$ is contained in a variety 
$\set{x\in T}{\phi(x)=0}$ where 
$\phi(x)=\sum_{\al\in\bbZ_{\geq0}^m}
 a_\al x_1^{\al_1}\cd x_m^{\al_m}\in \bbC[x_1,\ld,x_m]$
is a non-zero polynomial.
Set $A\seteq\set{\al=(\al_1,\ld,\al_m)
\in \bbZ_{\geq0}^m}{a_\al\ne0}\not=\emptyset$ and let $H$ 
be the convex hull of $A$ in $\bbZ^m$. 
Let us take an end point $\til\al$ of $H$. Then 
$\til\al\in A$ and there exists 
$c=(c_1,\cd,c_m)\in B=\bbZ^m$ such that 
$\sum_i\til\al_i c_i
>\sum_i\al_i c_i$
for any $\al\in A\setminus\{\til\al\}$.
Since the crystal $B$ is connected,
there exist $i_1,\ld,i_l\in I$ and 
$k_1,\ld,k_l\in \bbZ$ such that 
\[
 c=(c_1,\ld,c_m)=\til e_{i_1}^{k_1}\cd 
\til e_{i_l}^{k_l}(0,0,\ld,0).
\]
We define rational functions 
$f_1(t),\ld,f_m(t)$ by 
\[
 (f_1(t),\ld, f_m(t))\seteq 
e_{i_1}^{t^{k_1}}\cd
e_{i_l}^{t^{k_l}}{\bf a}\in Z.
\]
Since $\til e_i
={\mathcal UD}(e_{i,\theta})$, we have
$v(f_i(t))=c_i$,
where $v \colon\bbC(t)\setminus\{0\}\to \bbZ$ is 
as in \ref{positive-str}.
Thus we have
\[
 v(\phi(f_1(t),\cd,f_m(t)))
={\sum_i\til\al_i c_i},
\]
which implies $ \phi(f_1(t),\cd,f_m(t))\ne0$.
This is a contradiction. 
\QED

By the above proof, we obtain the following:
\begin{cor}
In the setting of Theorem \ref{conn-preh}, let $\Omega$ be the 
open dense orbit in $\bbX$ and 
identify $T$ and $X$ by $\theta$. 
Then we have $T_+\subset \Omega$.
\end{cor}
%%%%%%%%%%%%%% Section 4 %%%%%%%%%%%%%%%%
%\renewcommand{\thesection}{\arabic{section}}
\section{\bf Fundamental Representations and 
Perfect Crystals}\label{sec4}
%\setcounter{equation}{0}
%\renewcommand{\theequation}{\thesection.\arabic{equation}}

%%%%%%%%%%%%%%%%%%%%%%%%%%%%%%%%
\subsection{Affine weights}
\label{aff-wt}

Let $\ge$ be an affine Lie algebra and 
the sets $\mathfrak t$, 
$\{\al_i\}_{i\in I}$ 
and $\{\al^\vee_i\}_{i\in I}$ as in \ref{KM}. 
We take $\mathfrak t$ so that $\dim\mathfrak t=\sharp I+1$.
Let $\del\in Q_+$ be a unique element 
satisfying $$\set{\lm\in Q}{\lan \al^\vee_i,\lm\ran=0
\text{ for any $i\in I$}}=\bbZ\del,$$
and let ${\bf c}\in \sum_i\bbZ_{\geq0}\alpha_i^\vee\subset \ge$ 
be a unique central element
satisfying $$\set{h\in Q^\vee}{\lan h,\al_i\ran=0
\text{ for any }i\in I}=\bbZ{\bf c}.$$
We write (\cite[6.1]{Kac})
\begin{eqnarray}
{\bf c}=\sum_i a_i^\vee \al^\vee_i,\qq
\del=\sum_i a_i\al_i.\label{eq:cdel}
\end{eqnarray}
Let $(\ ,\ )$ be the non-degenerate
$W$-invariant symmetric bilinear form on $\mathfrak t^*$
normalized by $(\del,\lm)=\lan {\bf c},\lm\ran$
for $\lm\in\mathfrak{t}^*$.
Let us set $\tt^*_{\rm cl}\seteq \tt^*/\bbC\del$ and let
${\rm cl}\cl\tt^*\longrightarrow \tt^*_{\rm cl}$
be the canonical projection. 
Then we have 
$\tt^*_{\rm cl}\cong \sbigoplus_i(\bbC \al^\vee_i)^*$.
Set $\tt^*_0\seteq \set{\lm\in\tt^*}{\lan {\bf c},\lm\ran=0}$,
$(\tt^*_{\rm cl})_0\seteq {\rm cl}(\tt^*_0)$. 
Then we have a positive-definite
symmetric form on $(\tt^*_{\rm cl})_0$ 
induced by the one on 
$\tt^*$. 
Let $\Lm_i\in \tt^*_{\rm cl}$ $(i\in I)$ be a %classical 
weight such that $\lan \al^\vee_i,
\Lm_j\ran=\del_{i,j}$, which 
is called a {\em fundamental weight}.
We choose 
$P$ so that $P_{\rm cl}\seteq {\rm cl}(P)$ 
coincides with 
$\oplus_{i\in I}\bbZ\Lm_i$ and 
we call $P_{\rm cl}$ the
{\it classical weight lattice}.
%%%%%%%%%%%%%%%%%%%%%%
\subsection{Affinization}
Let $\uq=\lan e_i,f_i,q^h|i\in I,\,\,h\in P\ran$ 
be the quantum affine algebra associated with $P$ and 
$\uqp=\lan e_i,f_i,q^h|i\in I,\,\,h\in (P_{\rm cl})^*
\ran$ its subalgebra associated with $(P_{\rm cl})^*$.
Set ${\rm Mod}^f(\ge,P_{\rm cl})$ the category of 
a finite dimensional 
$\uqp$-module $M$ with a weight decomposition 
$M=\oplus_{\lm\in P_{\rm cl}}M_\lm$.

Let $M$ be an object in ${\rm Mod}^f(\ge,P_{\rm cl})$.
For $l\in\bbC^\times$, define the $\uqp$-module $M_l$
as follows:
There exists a $\bbC$-linear bijective homomorphism
$ \Phi_l:M\lar M,$
such that 
\begin{eqnarray*}
&&q^h\Phi_l(u)=\Phi_l(q^h u)\q
\text{for }h\in P_{\rm cl}^*,\\
&&e_i\Phi_l(u)=l^{\del_{i,0}}\Phi_l(e_i u),\\
&&f_i\Phi_l(u)=l^{-\del_{i,0}}\Phi_l(f_i u).
\end{eqnarray*}
The module $M_l:=\Phi_l(M)$ 
is said to be the {\it affinization} of $M$
(\cite{KMN1},\cite{K0}).
%%%%%%%%%% 
\subsection{Fundamental representation 
$W(\varpi_1)_l$}
\label{fundamental}

%Let $c=\sum_{i}a_i^\vee \al^\vee_i$ be the canonical
%central element of 
Here we consider the following affine 
Lie algebras $\ge=A_n^{(1)}$, $B_n^{(1)}$,
$D_n^{(1)}$, $A_{2n-1}^{(2)}$,
$D_{n+1}^{(2)}$, $A_{2n}^{(2)}$.
%(as for $a^\vee_i$, see \cite[6.1]{Kac}), 
Let $\set{\Lm_i}{i\in I}$ be the set of fundamental 
weights as in \S \ref{aff-wt}.
Let $\varpi_1\seteq \Lm_1-a^\vee_1\Lm_0$ be the
(level 0) fundamental weight, where
$i=1$ is the node of the Dynkin diagram as below
and $a^\vee_i$ is given in \eqref{eq:cdel}.

%Let $V(\varpi_1)$ be the extremal weight module
%over $\uq$
%associated with $\varpi_1$ (\cite{K0}) and 

Let $W(\varpi_1)$ be the 
fundamental representation of $\uqp$
 (see \cite[Sect 5.]{K0}). 
By \cite[Theorem 5.17]{K0}, $W(\varpi_1)$ is a
finite-dimensional irreducible integrable 
$\uqp$-module, an object in 
${\rm Mod}^f(\ge,P_{\rm lc})$ and has a global basis
with a simple crystal. Thus, we can consider its
affinization $W(\varpi_1)_l$ ($l\in \bbC^\times$)
specialized at $q=1$. Then we obtain a
finite-dimensional $\ge$-module 
denoted also by $W(\varpi_1)_l$
in which we shall construct affine geometric crystals
in the next section.

Note that $W(\varpi_1)_l$ is an irreducible 
$\ge$-module for any $\ge$ as above. But, 
for some $i$ and $\ge$, $W(\varpi_i)_l$ 
is not necessarily an irreducible $\ge$-module.

Let us see the explicit description of 
$W(\varpi_1)_l$ for 
$\ge=A_n^{(1)}$, $B_n^{(1)}$,
$D_n^{(1)}$, $A_{2n-1}^{(2)}$,
$D_{n+1}^{(2)}$, $A_{2n}^{(2)}$.

%%%%%%%%%%%%%
\subsection{$\TY(A,1,n)$ $(n\geq2)$}
\label{a-inf}

The Cartan matrix $(a_{ij})_{i,j\in I}$ 
$(I\seteq \{0,1,\ld,n\})$ of type $\TY(A,1,n)$ is
\[
 a_{ij}=\begin{cases}
2&\text{if $i=j$,}\\
-1&\text{if $i\equiv j\pm1\mod n+1$,}\\
0&\text{otherwise.}
\end{cases}
\]
We have ${\bf c}=\sum_{i\in I}\al^\vee_i$ and 
$\del=\sum_{i\in I}\al_i$.

The basis of $W(\varpi_1)_l$ is
$\{v_1,\,\,\fsquare(0.5cm,2),\cd,
\,\,\fsquare(0.5cm,n+1)\},$
and we have
\[
 \wt(\fsquare(0.5cm,i))=\Lm_i-\Lm_{i-1} \q(i=1\cd,n+1),
\]
where we understand $\Lm_{n+1}=\Lm_0$.
The explicit actions of $f_i$'s are  
\begin{eqnarray*}
&&f_i\fsquare(0.5cm,i)
=\fsquare(0.5cm,i+1)\q(1\leq i\leq n),\q
f_0\fsquare(0.5cm,n+1)=l^{-1}\fsquare(0.5cm,1),\\
&&f_i\fsquare(0.5cm,j)=0 \qq \text{otherwise,}
\end{eqnarray*}
and the explicit actions of $e_i$'s are  
\begin{eqnarray*}
&&e_i\fsquare(0.5cm,i+1)
=\fsquare(0.5cm,i)\q(1\leq i\leq n),\q
e_0\fsquare(0.5cm,1)=l\fsquare(0.5cm,n+1),\\
&&e_i\fsquare(0.5cm,j)=0 \qq \text{otherwise.}
\end{eqnarray*}
%%%%%%%%%%%%%%%%
\subsection{$\TY(B,1,n)$ $(n\geq2)$}
\label{bn-w1}

The Cartan matrix $(a_{ij})_{i,j\in I}$ 
$(I\seteq \{0,1,\ld,n\})$ of type $\TY(B,1,n)$ is
\[
 a_{ij}=\begin{cases}
2&i=j,\\
-1&|i-j|=1\text{ and }(i,j)\ne(0,1),(1,0),
(n,n-1)\text{ or }(i,j)=(0,2),(2,0)\\
-2&(i,j)=(n,n-1),\\
0&\text{otherwise.}
\end{cases}
\]
The Dynkin diagram is 
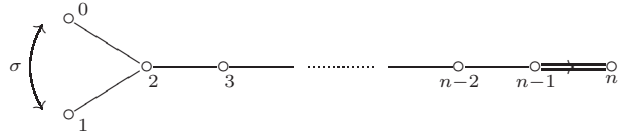
\[\SelectTips{cm}{}
\xymatrix@R=3ex{
*{\circ}<3pt> \ar@{-}[dr]^<{0} \ar@{<->}@/_/@<-2ex>[dd]_{\sigma}\\
&*{\circ}<3pt> \ar@{-}[r]_<{2} & *{\circ}<3pt> \ar@{-}[r]_<{3}
& {} \ar@{.}[r]&{} \ar@{-}[r]_>{\,\,\,\,n-2}
& *{\circ}<3pt> \ar@{-}[r]_>{\,\,\,\,n-1} &
*{\circ}<3pt> \ar@{=}[r] |-{\object@{>}}& *{\circ}<3pt>\ar@{}_<{n}\\
*{\circ}<3pt> \ar@{-}[ur]_<{1}
}
\]
where $\sigma$ is the Dynkin diagram 
automorphism $\sigma:\al_0\leftrightarrow\al_1$
which will be needed later. We have
\[
 {\bf c}=\al_0^\vee+\al_1^\vee+2\sum_{i=2}^{n-1}
{\al^\vee_i}+\al_n^\vee,\qq\q
\del=\al_0+\al_1+2\sum_{i=2}^n\al_i.
\]
The basis of $\W1_l$ is 
$\{\fsquare(0.5cm,1),\fsquare(0.5cm,2),\,\,\cd,
\,\,\fsquare(0.5cm,n),\,\,\fsquare(0.5cm,0),
\,\,\fsquare(0.5cm,\ovl n),\cd,
\,\,\fsquare(0.5cm,\ovl 2),
\,\,\fsquare(0.5cm,\ovl 1)\},$
and we have
\begin{eqnarray*}
&&\wt(\fsquare(5mm,i))=\Lm_i-\Lm_{i-1},\qq
\wt(\fsquare(5mm,\ovl i))=\Lm_{i-1}-\Lm_i
\qq(i\ne 0,2,n),\\
&&\wt(\fsquare(5mm,2))=-\Lm_0-\Lm_1+\Lm_2,\qq
\wt(\fsquare(5mm,\ovl 2))=\Lm_0+\Lm_1-\Lm_2,\\
&&\wt(\fsquare(5mm,n))=2\Lm_n-\Lm_{n-1},\qq
\wt(\fsquare(5mm,\ovl n))=\Lm_{n-1}-2\Lm_n,\qq
\wt(\fsquare(5mm,0))=0.
\end{eqnarray*}
The explicit forms of the actions by $f_i$'s and $e_i$'s are
\begin{eqnarray*}
&f_i\fsquare(5mm,i)=\fsquare(5mm,i+1),\qq
f_i\fsquare(5mm,\ovl {i+1})=\fsquare(5mm,\ovl i)\q
(1\leq i<n),
&e_i\fsquare(5mm,i+1)=\fsquare(5mm,i),\qq
e_i\fsquare(5mm,\ovl {i})=\fsquare(5mm,\ovl {i+1})\q
(1\leq i<n),\\
&f_{n}\fsquare(5mm,n)=\fsquare(5mm,0),\qq
f_n\fsquare(5mm,0)=2\fsquare(5mm,\ovl n),\qq
&e_{n}\fsquare(5mm,0)=2\fsquare(5mm,n),\qq
e_n\fsquare(5mm,\ovl n)=\fsquare(5mm,0),\\
&f_0\fsquare(5mm,\ovl 2)=l^{-1}\fsquare(5mm,1),\qq
f_0\fsquare(5mm,\ovl 1)=l^{-}\fsquare(5mm,2),\qq
&e_0\fsquare(5mm,1)=l\fsquare(5mm,\ovl 2),\qq
e_0\fsquare(5mm,2)=l\fsquare(5mm,\ovl 1),\\
&f_i\fsquare(5mm,j)=0\qq\text{otherwise},\qq
&e_i\fsquare(5mm,j)=0\qq\text{otherwise},
\end{eqnarray*}

\subsection{$\TY(D,1,n)$ $(n\geq 4)$}
\label{dn-w1}

The Cartan matrix $(a_{ij})_{i,j\in I}$ 
$(I\seteq \{0,1,\ld,n\})$ of type $\TY(D,1,n)$ is
\[
 a_{ij}=\begin{cases}
2&i=j,\\
-1&|i-j|=1 \text{ and }1\leq i,j\leq n-1
\text{ or }(i,j)=(0,2),(2,0),
(n-2,n), (n,n-2),\\
0&\text{otherwise.}
\end{cases}
\]
The Dynkin diagram is 
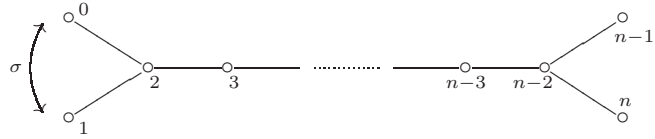
\[\SelectTips{cm}{}
\xymatrix@R=3ex{
*{\circ}<3pt> \ar@{-}[dr]^<{0} \ar@{<->}@/_/@<-2ex>[dd]_{\sigma}&&&&&&&*{\circ}<3pt> \ar@{-}[dl]^<{n-1}\\
&*{\circ}<3pt> \ar@{-}[r]_<{2} & *{\circ}<3pt> \ar@{-}[r]_<{3}
& {} \ar@{.}[r]&{} \ar@{-}[r]_>{\,\,\,\,n-3}
& *{\circ}<3pt> \ar@{-}[r]_>{n-2\,\,\,} &*{\circ}<3pt>\\
*{\circ}<3pt> \ar@{-}[ur]_<{1}&&&&&&&
*{\circ}<3pt> \ar@{-}[ul]_<{n}
}
\]
%\input{Dn-Dynkin.tex}
%\vskip3mm
where $\sigma$ is the Dynkin diagram automorphism
$\sigma\cl\al_0\leftrightarrow \al_1$ and 
$\sigma\alpha_i=\alpha_i$ for $i\not=0,1$. We have
\[
 {\bf c}=\al_0^\vee+\al_1^\vee+2\sum_{i=2}^{n-2}
{\al^\vee_i}+\al^\vee_{n-1}+\al_n^\vee,\qq\q
\del=\al_0+\al_1+2\sum_{i=2}^{n-2}\al_i
+\al_{n-1}+\al_n.
\]
The basis of $\W1_l$ is 
$\{\fsquare(0.5cm,1),\,\,
\fsquare(0.5cm,2),\cd,
\,\,\fsquare(0.5cm,n),
\,\,\fsquare(0.5cm,\ovl n),\cd,
\,\,\fsquare(0.5cm,\ovl 2),
\,\,\fsquare(0.5cm,\ovl 1)\},$
and we have
\begin{eqnarray*}
&&\wt(\fsquare(5mm,i))=\Lm_i-\Lm_{i-1},\qq
\wt(\fsquare(5mm,\ovl i))=\Lm_{i-1}-\Lm_i
\qq(i\ne 2,n-1),\\
&&\wt(\fsquare(5mm,2))=-\Lm_0-\Lm_1+\Lm_2,\qq
\wt(\fsquare(5mm,\ovl 2))=\Lm_0+\Lm_1-\Lm_2,\\
&&\wt(\fsquare(5mm,n-1))=\Lm_{n-1}+
\Lm_n-\Lm_{n-1},\qq
\wt(\fsquare(6mm,\ovl{n-1}))
=\Lm_{n-2}-\Lm_{n-1}-\Lm_n.
\end{eqnarray*}
The explicit forms of the actions by $f_i$'s and $e_i$'s are
\begin{eqnarray*}
&f_i\fsquare(5mm,i)=\fsquare(5mm,i+1),\qq
f_i\fsquare(5mm,\ovl {i+1})=\fsquare(5mm,\ovl i)
\q(1\leq i<n),
&e_i\fsquare(5mm,i+1)=\fsquare(5mm,i),\qq
e_i\fsquare(5mm,\ovl {i})=\fsquare(5mm,\ovl{i+1})
\q(1\leq i<n),\\
&f_{n}\fsquare(5mm,n)=\fsquare(6mm,\ovl{n-1}),
\qq
f_n\fsquare(5mm,n-1)=\fsquare(5mm,\ovl n),\qq\qq
&e_{n}\fsquare(5mm,\ovl n)=\fsquare(6mm,{n-1}),
\qq
e_n\fsquare(5mm,\ovl{n-1})=\fsquare(5mm,n),\\
&f_0\fsquare(5mm,\ovl 2)=l^{-1}\fsquare(5mm,1),\qq
f_0\fsquare(5mm,\ovl 1)=l^{-1}\fsquare(5mm,2),\qq\qq
&e_0\fsquare(5mm,1)=l\fsquare(5mm,\ovl 2),\qq
e_0\fsquare(5mm,2)=l\fsquare(5mm,\ovl 1),\\
&f_i\fsquare(5mm,j)=0\qq\text{otherwise},\qq\qq
&e_i\fsquare(5mm,j)=0\qq\text{otherwise}.
\end{eqnarray*}

%%%%%%%%%%%%%%%%%%%%%%
%\clearpage

\subsection{$\TY(A,2,2n-1)$ $(n\geq 3)$}
\label{aon-w1}

The Cartan matrix $(a_{ij})_{i,j\in I}$ 
$(I\seteq \{0,1,\ld,n\})$ of type $\TY(A,2,2n-1)$ is
\[
 a_{ij}=\begin{cases}
2&i=j,\\
-1&|i-j|=1\text{ and }1\leq i,j\leq n-1
\text{ or }(i,j)=(0,2),(2,0),(n,n-1)\\
-2&(i,j)=(n-1,n),\\
0&\text{otherwise.}
\end{cases}
\]
The Dynkin diagram is 
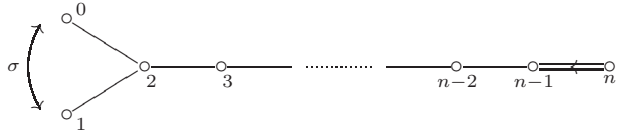
\[\SelectTips{cm}{}
\xymatrix@R=3ex{
*{\circ}<3pt> \ar@{-}[dr]^<{0} \ar@{<->}@/_/@<-2ex>[dd]_{\sigma}\\
&*{\circ}<3pt> \ar@{-}[r]_<{2} & *{\circ}<3pt> \ar@{-}[r]_<{3}
& {} \ar@{.}[r]&{} \ar@{-}[r]_>{\,\,\,\,n-2}
& *{\circ}<3pt> \ar@{-}[r]_>{\,\,\,\,n-1} &
*{\circ}<3pt> \ar@{=}[r] |-{\object@{<}}& *{\circ}<3pt>\ar@{}_<{n}\\
*{\circ}<3pt> \ar@{-}[ur]_<{1}
}
\]
where $\sigma$ is the Dynkin diagram 
automorphism $\sigma\cl\al_0\leftrightarrow \al_1$. 
We have
\[
 {\bf c}=\al_0^\vee+\al_1^\vee+2\sum_{i=2}^{n}
{\al^\vee_i},\qq\q
\del=\al_0+\al_1+2\sum_{i=2}^{n-1}\al_i
+\al_n.
\]
The basis of $\W1_l$ is 
$\{\fsquare(0.5cm,1),
\,\,\fsquare(0.5cm,2),\cd,
\,\,\fsquare(0.5cm,n),
\,\,\fsquare(0.5cm,\ovl n),\cd,
\,\,\fsquare(0.5cm,\ovl 2),
\,\,\fsquare(0.5cm,\ovl 1)\},$
and we have
\begin{eqnarray*}
&&\wt(\fsquare(5mm,i))=\Lm_i-\Lm_{i-1},\qq
\wt(\fsquare(5mm,\ovl i))=\Lm_{i-1}-\Lm_i
\qq(i\ne 2)\\
&&\wt(\fsquare(5mm,2))=-\Lm_0-\Lm_1+\Lm_2,\qq
\wt(\fsquare(5mm,\ovl 2))=\Lm_0+\Lm_1-\Lm_2.
\end{eqnarray*}
The explicit forms of the actions by $f_i$'s 
and $e_i$'s are
\begin{eqnarray*}
&f_i\fsquare(5mm,i)=\fsquare(5mm,i+1),\qq
f_i\fsquare(5mm,\ovl {i+1})=\fsquare(5mm,\ovl i)
\q(1\leq i<n),
&e_i\fsquare(5mm,i+1)=\fsquare(5mm,i),\qq
e_i\fsquare(5mm,\ovl {i})=\fsquare(5mm,\ovl{i+1})
\q(1\leq i<n),\\
&f_{n}\fsquare(5mm,n)=\fsquare(5mm,\ovl{n}),\qq\qq\qq\qq\qq
&e_{n}\fsquare(5mm,\ovl n)=\fsquare(5mm,{n}),\\
&f_0\fsquare(5mm,\ovl 2)=l^{-1}\fsquare(5mm,1),\qq
f_0\fsquare(5mm,\ovl 1)=l^{-1}\fsquare(5mm,2),\qq\qq
&e_0\fsquare(5mm,1)=l\fsquare(5mm,\ovl 2),\qq
e_0\fsquare(5mm,2)=l\fsquare(5mm,\ovl 1),\\
&f_i\fsquare(5mm,j)=0\qq\text{otherwise},\qq\qq
&e_i\fsquare(5mm,j)=0\qq\text{otherwise}.
\end{eqnarray*}

%%%%%%%%%%%%%%%%%%%%%
\subsection{$\TY(D,2,n+1)$ $(n\geq 2)$}
\label{d2n-w1}

The Cartan matrix $(a_{ij})_{i,j\in I}$ 
$(I\seteq \{0,1,\ld,n\})$ of type $\TY(D,2,n+1)$ is
\[
 a_{ij}=\begin{cases}
2&i=j,\\
-1&|i-j|=1\text{ and }(i,j)\ne(0,1),
(n,n-1),\\
-2&(i,j)=(0,1),(n,n-1),\\
0&\text{otherwise.}
\end{cases}
\]
The Dynkin diagram is
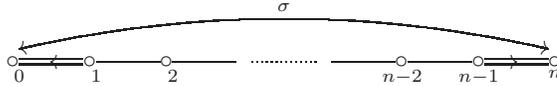
\[\SelectTips{cm}{}
\xymatrix{
*{\circ}<3pt> \ar@{=}[r] |-{\object@{<}}_<{0} 
\ar@{<->}@/^1pc/@<1ex>[rrrrrrr]^{\sigma}
&*{\circ}<3pt> \ar@{-}[r]_<{1} & *{\circ}<3pt> \ar@{-}[r]_<{2}
& {} \ar@{.}[r]&{} \ar@{-}[r]_>{\,\,\,\,n-2}
& *{\circ}<3pt> \ar@{-}[r]_>{\,\,\,\,n-1} &
*{\circ}<3pt> \ar@{=}[r] |-{\object@{>}}
& *{\circ}<3pt>\ar@{}_<{n}
}
\]
%\input{D^2_n+1-Dynkin.tex}
%\vskip3mm
where $\sigma$ is the Dynkin diagram 
automorphism $\sigma\cl\al_i
\leftrightarrow \al_{n-i}$ $(i=0,1\cd,n)$. We have
\[
 {\bf c}=\al_0^\vee+\al_1^\vee+2\sum_{i=2}^{n-1}
{\al^\vee_i}+\al_n^\vee,\qq\q
\del=\sum_{i=0}^n\al_i.
\]
The basis of $\W1_l$ is 
$\{\fsquare(0.5cm,1),
\,\,\fsquare(0.5cm,2),\cd,
\,\,\fsquare(0.5cm,n),
\,\,\fsquare(0.5cm,0),
\,\,\fsquare(0.5cm,\ovl n),\cd,
\,\,\fsquare(0.5cm,\ovl 2),
\,\,\fsquare(0.5cm,\ovl 1),
\,\,\phi\},$
and we have
\begin{eqnarray*}
&&\wt(\fsquare(5mm,i))=\Lm_i-\Lm_{i-1},\qq
\wt(\fsquare(5mm,\ovl i))=\Lm_{i-1}-\Lm_i
\qq(i\ne 0,1,n),\\
&&\wt(\fsquare(5mm,1))=\Lm_1-2\Lm_0,\qq
\wt(\fsquare(5mm,\ovl 1))=2\Lm_0-\Lm_1,\\
&&\wt(\fsquare(5mm,n))=2\Lm_n-\Lm_{n-1},\qq
\wt(\fsquare(5mm,\ovl n))=\Lm_{n-1}-2\Lm_n,\\
&&\wt(\fsquare(5mm,0))=0,\qq
\wt(\phi)=0
\end{eqnarray*}
The explicit forms of the actions 
by $f_i$'s and $e_i$'s are
\begin{eqnarray*}
&f_i\fsquare(5mm,i)=\fsquare(5mm,i+1),\qq
f_i\fsquare(5mm,\ovl {i+1})=\fsquare(5mm,\ovl i)
\q(1\leq i<n),
&e_i\fsquare(5mm,i+1)=\fsquare(5mm,i),\qq
e_i\fsquare(5mm,\ovl {i})=\fsquare(5mm,\ovl{i+1})
\q(1\leq i<n),\\
&f_{n}\fsquare(5mm,n)=\fsquare(5mm,0),\qq
f_n\fsquare(5mm,0)=2\fsquare(5mm,\ovl n),\qq\qq
&e_{n}\fsquare(5mm,0)=2\fsquare(5mm,n),\qq
e_n\fsquare(5mm,\ovl n)=\fsquare(5mm,0),\\
&f_0\fsquare(5mm,\ovl 1)=l^{-1}\phi,\qq
f_0\phi=2l^{-1}\fsquare(5mm,1),\qq\qq
&e_0\fsquare(5mm,1)=l\phi,\qq
e_0\phi=2l\fsquare(5mm,\ovl 1),\\
&f_i\fsquare(5mm,j)=0,\q
f_i\phi=0\qq\text{otherwise},\qq\qq
&e_i\fsquare(5mm,j)=0,\q
e_i\phi=0\qq\text{otherwise}.
\end{eqnarray*}

\subsection{$A_{2n}^{(2)}$ $(n\geq 2)$}
\label{aen-v1}

The Cartan matrix $(a_{ij})_{i,j\in I}$ 
$(I\seteq \{0,1,\ld,n\})$ of type 
${\TY(A,2,2n)}$ is
\[
 a_{ij}=\begin{cases}
2&i=j,\\
-1&|i-j|=1\text{ and }(i,j)\ne(0,1),
(n-1,n),\\
-2&(i,j)=(0,1),(n-1,n),\\
0&\text{otherwise.}
\end{cases}
\]
Then the Dynkin diagram is 
\[\SelectTips{cm}{}
\xymatrix{
*{\circ}<3pt> \ar@{=}[r] |-{\object@{<}}_<{0} 
%\ar@{<->}@/^1pc/@<1ex>[rrrrrrr]^{\sigma}
&*{\circ}<3pt> \ar@{-}[r]_<{1} & *{\circ}<3pt> \ar@{-}[r]_<{2}
& {} \ar@{.}[r]&{} \ar@{-}[r]_>{\,\,\,\,n-2}
& *{\circ}<3pt> \ar@{-}[r]_>{\,\,\,\,n-1} &
*{\circ}<3pt> \ar@{=}[r] |-{\object@{<}}
& *{\circ}<3pt>\ar@{}_<{n}
}
\]
Note that there exists no Dynkin diagram 
automorphism in this case. We have
\[
 {\bf c}=\al^\vee_0+2\sum_{i=1}^{n}\al^\vee_i,\qq\q
\del=2\sum_{i=0}^{n-1}\al_i+\al_n.
\]
In this case, we denote this type by 
$(A_{2n}^{(2)},\varpi_1)$ in order to distinguish it 
with the type $({\TY(A,2,2n)}\dagger,\varpi_1)$.
Then the basis of $\W1_l$ is 
$\{\fsquare(0.5cm,1),
\,\,\fsquare(0.5cm,2),\cd,
\,\,\fsquare(0.5cm,n),
\,\,\fsquare(0.5cm,\ovl n),\cd,
\,\,\fsquare(0.5cm,\ovl 2),
\,\,\fsquare(0.5cm,\ovl 1),\,\,\emptyset
\},$
and we have
\begin{eqnarray*}
&&\wt(\fsquare(5mm,1))=\Lm_1-2\Lm_0,\qq
\wt(\fsquare(5mm,i))=\Lm_i-\Lm_{i-1},\qq
\wt(\fsquare(5mm,\ovl i))=\Lm_{i-1}-\Lm_i
\qq(i=2,\ld,n),\\
&&\wt(\fsquare(5mm,\ovl 1))=2\Lm_0-\Lm_1,\qq
\wt(\emptyset)=0.
\end{eqnarray*}
The explicit forms of the actions by 
$f_i$'s and $e_i$'s are
\begin{eqnarray*}
&f_i\fsquare(5mm,i)=\fsquare(5mm,i+1),\q
f_i\fsquare(5mm,\ovl {i+1})=\fsquare(5mm,\ovl i)
\q(1\leq i<n),
&e_i\fsquare(5mm,i+1)=\fsquare(5mm,i),\qq
e_i\fsquare(5mm,\ovl {i})=\fsquare(5mm,\ovl{i+1})
\q(1\leq i<n),
\\
&f_{n}\fsquare(5mm,n)=\fsquare(5mm,\ovl n),\qq\qq\qq\qq\qq
&e_{n}\fsquare(5mm,n)=\fsquare(5mm,n),\\
&f_0\fsquare(5mm,\ovl 1)=l^{-1}\emptyset,\q
f_0\emptyset=2l^{-1}\fsquare(5mm,1),\qq\qq
&e_0\fsquare(5mm,1)=l\emptyset,\qq
e_0\emptyset=2lv_{\ovl 1},
\\
&f_iv_j=0\qq
\text{otherwise},\qq\qq
&e_i\fsquare(5mm,j)=0,\q
\text{otherwise}.
\end{eqnarray*}

%%%%%%%%%%%%%%%%%%%%%%%%%%%%%%%%%%%%%
\subsection{$A_{2n}^{(2)\,\dagger}$ $(n\geq 2)$}
\label{aen-v1-dag}

The Cartan matrix $(a_{ij})_{i,j\in I}$ 
$(I\seteq \{0,1,\ld,n\})$ of type 
${\TY(A,2,2n)}^\dagger$ is
\[
 a_{ij}=\begin{cases}
2&i=j,\\
-1&|i-j|=1\text{ and }(i,j)\ne(1,0),
(n,n-1),\\
-2&(i,j)=(1,0),(n,n-1),\\
0&\text{otherwise.}
\end{cases}
\]
Then the Dynkin diagram is
\[\SelectTips{cm}{}
\xymatrix{
*{\circ}<3pt> \ar@{=}[r] |-{\object@{>}}_<{0} 
%\ar@{<->}@/^1pc/@<1ex>[rrrrrrr]^{\sigma}
&*{\circ}<3pt> \ar@{-}[r]_<{1} & *{\circ}<3pt> \ar@{-}[r]_<{2}
& {} \ar@{.}[r]&{} \ar@{-}[r]_>{\,\,\,\,n-2}
& *{\circ}<3pt> \ar@{-}[r]_>{\,\,\,\,n-1} &
*{\circ}<3pt> \ar@{=}[r] |-{\object@{>}}
& *{\circ}<3pt>\ar@{}_<{n}
}
\]
Note that there exists no Dynkin diagram 
automorphism in this case. We have
\[
 {\bf c}=2\sum_{i=0}^{n-1}\al^\vee_i+\al_n^\vee
,\qq\q
\del=\al_0+2\sum_{i=1}^{n}\al_i.
\]
In this case, we denote this type by 
$(A_{2n}^{(2)\,\dagger},\varpi_1)$ in order to distinguish it 
with the type $(\TY(A,2,2n),\varpi_1)$.
Then the basis of $\W1_l$ is 
$\{\fsquare(0.5cm,1),
\,\,\fsquare(0.5cm,2),\cd,
\,\,\fsquare(0.5cm,n),
\,\,\fsquare(0.5cm,0),
\,\,\fsquare(0.5cm,\ovl n),\cd,
\,\,\fsquare(0.5cm,\ovl 2),
\,\,\fsquare(0.5cm,\ovl 1)
\},$
and we have
\begin{eqnarray*}
&&\wt(\fsquare(5mm,i))=\Lm_i-\Lm_{i-1},\qq
\wt(\fsquare(5mm,\ovl i))=\Lm_{i-1}-\Lm_i
\qq(i=1,\ld,n-1),\\
&&\wt(\fsquare(5mm,n))=2\Lm_n-\Lm_{n-1},\qq
\wt(\fsquare(5mm,0))=0,\qq
\wt(\fsquare(5mm,\ovl 1))=\Lm_{n-1}-2\Lm_n.
\end{eqnarray*}
The explicit forms of the actions by 
$f_i$'s and $e_i$'s are
\begin{eqnarray*}
&f_i\fsquare(5mm,i)=\fsquare(5mm,i+1),\qq
f_i\fsquare(5mm,\ovl {i+1})=\fsquare(5mm,\ovl i)
\q(1\leq i<n),
&e_i\fsquare(5mm,i+1)=\fsquare(5mm,i),\qq
e_i\fsquare(5mm,\ovl {i})=\fsquare(5mm,\ovl{i+1})
\q(1\leq i<n),
\\
&f_{n}\fsquare(5mm,n)=\fsquare(5mm,\ovl 0),
\qq f_{n}\fsquare(5mm,0)=
2\fsquare(5mm,\ovl n),\qq\qq
&e_{n}\fsquare(5mm,0)=2\fsquare(5mm,n),
\qq 
e_{n}\fsquare(5mm,\ovl n)=\fsquare(5mm,0),
\\
&f_0\fsquare(5mm,\ovl 1)=l^{-1}\fsquare(5mm,1),
\qq\qq\qq\qq\qq
&e_0\fsquare(5mm,1)=l\fsquare(5mm,\ovl 1),\\
&f_i\fsquare(5mm,j)=0,\q
\text{otherwise},\qq\qq
&e_i\fsquare(5mm,j)=0,\q
\text{otherwise}.
\end{eqnarray*}

%%%%%%%%%%%%%%%%%%%%%%%%%%%%
\subsection{Limit of perfect crystals}
\label{limit}

We review the limits of perfect crystals following \cite{KKM}.
(See also \cite{KMN1},\cite{KMN2}.)

Let $\ge$ be an affine Lie algebra, $P_{\mathrm{cl}}$ 
the classical weight lattice as above and for $l\in\ZZ_{>0}$
set 
$(P_{\mathrm{cl}})^+_l\seteq \set{\lm\in P_{\mathrm{cl}}}{
\lan {\bf c},\lm\ran=l,\,\,\lan \al^\vee_i,\lm\ran\geq0}$.
\begin{df}
\label{perfect-def}
We say that a crystal $B$ is {\it perfect} of level $l$ if 
\bnum
\item
$B\ot B$ is connected as a crystal graph.
\item
There exists $\lm_0\in P_{\rm cl}$ such that 
\[
 \wt(B)\subset \lm_0+\sum_{i\ne0}\ZZ_{\leq0}
{\rm cl}(\al_i),\qq
\sharp B_{\lm_0}=1
\]
\item There exists a finite-dimensional 
$U'_q(\ge)$-module $V$ with a
crystal pseudo-base $B_{ps}$ 
such that $B\cong B_{ps}/{\pm1}$
\item
Set $\vep(b)\seteq \sum_i\vep_i(b)\Lm_i$ and 
$\vp(b)\seteq \sum_i\vp_i(b)\Lm_i$.
For any $b\in B$, $\lan{\bf c},\vep(b)\ran\geq l$
and the maps 
$\vep,\vp\cl B^{min}\seteq \set{b\in B}{\lan c,\vep(b)\ran=l}
\mapright{}(P_{\rm cl}^+)_l$ are bijective.
\enum
\end{df}

For an affine Lie algebra $\ge$, 
let $\{B_l(\ge)\}_{l\geq1}$ be a family of 
perfect crystals of level $l$ and set 
$J\seteq \set{(l,b)}{l>0,\,b\in B^{min}_l}$.
\begin{df}
\label{def-limit}
A crystal $B_\ify=B_\ify(\ge)$ 
with an element $b_\ify$ is called the
{\it limit of $\{B_l\}_{l\geq1}$}
if 
\bnum
\item
$\wt(b_\ify)=\vep(b_\ify)=\vp(b_\ify)=0$.
\item
For any $(l,b)\in J$, there exists an
embedding of 
crystals:
\begin{eqnarray*}
 f_{(l,b)}\cl&
T_{\vep(b)}\ot B_l\ot T_{-\vp(b)}\hookrightarrow
B_\ify\\
&t_{\vep(b)}\ot b\ot t_{-\vp(b)}\mapsto b_\ify
\end{eqnarray*}
\item
$B_\ify=\bigcup_{(l,b)\in J} {\rm Im}f_{(l,b)}$.
\end{enumerate}
\end{df}
\noindent
As for the crystal $T_\lm$, see Example \ref{ex-tlm} \eqref{tlm}.
If the limit of a family $\{B_l\}$ exists, 
we say that $\{B_l\}$
is a {\it coherent family} of perfect crystals.

{\sl Remark. } By the definition of perfect crystals,
any perfect crystal is connected and then 
the limit of a coherent family of perfect crystals 
is also connected.

%%%%%%%%%%% Section 5 %%%%%%%%
\section{\bf 
Affine Geometric Crystals}\label{sec5}

Following the method in \cite{KNO},
we shall construct the  
affine geometric crystal $\cV(\ge)_l$ 
$(l\in\bbC^\times)$ in the $\ge$-module 
$W(\varpi_1)_l$
the affinization of 
the fundamental representation $W(\varpi_1)$. 

%%%%%%%%%%%%%%%%%%%%%%%%%
\subsection{Translation $t(\til\varpi_1)$}
\label{shift}%5.1
For $\xi_0\in (\frt^*_{\rm cl})_0$, let $t(\xi_0)$ be 
as in \cite[Sect.4]{K0}, that is, 
\[
 t(\xi_0)(\lm)\seteq \lm+(\del,\lm)
\xi-
(\xi,\lm)\del
-\frac{(\xi,\xi)}{2}
(\del,\lm)\del
\]
for $\xi\in \mathfrak{t}^*$ such that $\mathrm{cl}(\xi)=\xi_0$.
Then $t(\xi_0)$ does not depend on the choice of $\xi$,
and it is well-defined.

Let  $c_i^\vee$ be as follows(\cite{K1}):
\begin{eqnarray}
&& c_i^\vee\seteq%\begin{cases}
\mathrm{max }(1,\frac{2}{(\al_i,\al_i)}).%&\text{if }\ge
%\text{ is untwisted,}\\
%1&\text{if }\ge\text{ is twisted.}
%\end{cases}
\label{eq:ci}
\end{eqnarray}

Then $t(m\varpi_i)$ belongs to the extended Weyl group
$\widetilde W$ if and only if 
$m\in c_i^\vee\bbZ$.
Setting $\wtil\varpi_i\seteq c_i^\vee\varpi_i$ $(i\in I)$,
$t(\wtil\varpi_1)$ 
is expressed as follows (see {\it e.g.} \cite{KMOTU}):
\[
 t(\wtil\varpi_1)=
\begin{cases}
\io(s_ns_{n-1}\cd s_2s_1)&\TY(A,1,n)\text{ case},\\
\io(s_1\cd s_n)(s_{n-1}\cd s_2s_1)
&\TY(B,1,n),\,\TY(A,2,2n-1)\text{ cases},\\
(s_0s_1\cd s_n)(s_{n-1}\cd s_2s_1)
\quad&\TY(C,1,n),\,\TY(D,2,n+1)\text{ cases},\\
\io(s_1\cd s_n)(s_{n-2}\cd s_2s_1)
&\TY(D,1,n)\text{ case},\\
(s_0s_1\cd s_n)(s_{n-1}\cd s_2s_1)
%&\TY(A,2,2n),\,A_{2n}^{(2)\,\dagger}\text{ cases},
&\TY(A,2,2n),\,\,
A_{2n}^{(2)\,\dagger}\text{ cases},
\end{cases}
\]
where $\io$ is the Dynkin diagram automorphism
\[
\io=
 \begin{cases}\sigma&\ge=\TY(A,1,n), \TY(B,1,n),
\TY(A,2,2n-1),\\
\al_0\leftrightarrow \al_1\text{ and }
\al_{n-1}\leftrightarrow \al_{n}&\ge=\TY(D,1,n).
\end{cases}
\]
Now, we know that each $t(\wtil\varpi_1)$ is 
in the form $w_1$ or $\io\cdot w_1$ for some $w_1\in W$,
%and denote this $w$ by $w_1$, 
{\it e.g.,} $w_1=s_n\cd s_1$ for $\TY(A,1,n)$,
$w_1=(s_1\cd s_n)(s_{n-1}\cd s_1)$ for 
$\TY(B,1,n)$, {\it etc.}, \ldots.

In the case $\TY(A,2,2n)$
(resp.~$\ge=A_{2n}^{(2)\,\dagger}$), 
$\eta\seteq 
\on{wt}(\fsquare(5mm,\ovl n))
=\Lm_{n-1}-\Lm_n$ 
(resp.~$\Lm_{n-1}-2\Lm_n$)
is a unique weight of $W(\varpi_1)_l$
which satisfies $\lan \al^\vee_i,\eta\ran\geq0$
for $i\ne n$.
For this $\eta$ we have
\begin{equation}
t(\eta)=(s_ns_{n-1}\cd s_1)(s_0s_1\cd s_{n-1})
=:w_2,
\label{t-pi-n}
\end{equation}
which will be  used later.
%%%%%%%%%%%%%%%%%%%%%%%%%%%%%%%%%%%%%
\subsection{Affine geometric crystals 
in $W(\varpi_1)_l$}%5.2

Let $\sigma$ be the Dynkin diagram automorphism
as in \ref{a-inf}--\ref{d2n-w1} and 
$w_1=s_{i_1}\cd s_{i_k}$ be as 
in the previous subsection.
Let $H$ be an element in $\tt$ such that 
\[
 \al_i(H)=\begin{cases}1&\text{ if }i=1 
\text{ and }\ge\ne\TY(D,2,n+1),\TY(A,2,2n),
{\TY(A,2,2n)}^\dagger,\\
2&\text{ if }i=1 
\text{ and }\ge=\TY(D,2,n+1),\TY(A,2,2n),
{\TY(A,2,2n)}^\dagger,\\
-1&\text{ if }i=0\text{ and }\ge\ne
\TY(D,2,n+1),\TY(A,2,2n),
{\TY(A,2,2n)}^\dagger,\\
-2&\text{ if }i=0\text{ and }\ge=\TY(D,2,n+1),
\TY(A,2,2n), {\TY(A,2,2n)}^\dagger,\\
0&\text{ otherwise.}
\end{cases}
\]
Set
\begin{equation}
{\cV}(\ge)_l\seteq 
\{v(x_1,\ldots,x_k)\seteq 
Y_{i_1}(x_1)\cd Y_{i_k}(x_k)l^H\fsquare(5mm,1)\,
 \big\vert \,
x_1,\ld,x_k\in \bbC^\times\}\subset W(\varpi_1)_l
\label{Vxrn}
\end{equation}
Let $\ge_0\subset \ge$ (resp. $G_0\subset G$) 
be a simple Lie algebra
(resp. simple algebraic group) corresponding to 
the index set $I_0\seteq I\setminus\{0\}$.
Since the vector $\fsquare(5mm,1)$ 
satisfies $x_i(c)v_1=v_1$ for any $i\in I_0$,
the actions of $e_i^c$ ($i\in I_0$)
on $v(x)$ and $Y_{i_1}(x_1)\cd Y_{i_k}(x_k)
\cdot l^H$ coincide each other. Therefore, 
$\cV(\ge)_l$ has a 
$G_0$-geometric crystal structure same as 
that of 
$B^-_{i_1\cd i_k}\cdot l^H$ 
(see \ref{schubert}), 
Moreover $(\bbC^\times)^k\to {\cV}(\ge)_l$ 
given by $(x_1,\ldots,x_k)\mapsto 
v(x_1,\ld,x_k)$ is a birational map.
We shall define a $G$-geometric crystal structure 
on $\cV(\ge)_l$ by using the Dynkin diagram 
automorphism $\sigma$ except for $\TY(A,2,2n)$. 
This $\sigma$ induces an automorphism 
of $W(\varpi_1)_l$, which is denoted by 
$\sigma_l\cl W(\varpi_1)_l\longrightarrow 
W(\varpi_1)_l$.
The following theorems 
are analogous results to 
Theorem 5.1 and 5.2 in \cite{KNO}.
\begin{thm}
\label{birat}
\bnum
\item
Case $\ge\ne A_{2n}^{(2)},\,
{\TY(A,2,2n)}^\dagger$.
For $x=(x_1\cd,x_k)\in (\bbC^\times)^k$, there 
exist a unique $y=(y_1,\ld,y_k)$
$\in (\bbC^\times)^k$
and a positive rational function $a(x)$ such that
\begin{equation}
v(y)=a(x)\sigma_l(v(x)),\q
\vep_{\sigma(i)}(v(y))=\vep_i(v(x))\quad
\text{if $i,\sigma(i)\not=0$.}
\end{equation}
\item
Case $\ge=A_{2n}^{(2)}$ 
(resp. ${\TY(A,2,2n)}^\dagger$). 
Associated with 
$w_1$ and $w_2$ as in the previous section, we define
\begin{eqnarray*}
 \cV(\ge)_l&\seteq &\{v_1(x)
=Y_0(x_0)Y_1(x_1)\cd Y_n(x_n)
Y_{n-1}(\bar x_{n-1})\cd Y_1(\bar x_1)l^H
\fsquare(5mm,1)\,\,\big\vert\,\,x_i,\bar x_i\in\bbC^\times\},\\
\cV_2(\ge)_l&\seteq &\{v_2(y)
=Y_n(y_n)\cd Y_1(y_1)Y_0(y_0)Y_{1}(\ovl y_1)\cd
Y_{n-1}(\ovl y_{n-1})l^{H'}\fsquare(5mm,\ovl n)
\,\big\vert\,y_i,\ovl y_i\in\bbC^\times\},
\end{eqnarray*}
where $\al_0(H')=2$ (resp. $\al_0(H')=2$),
$\al_n(H')=-4$ (resp. $\al_n(H')=-1$) and 
$\al_i(H')=0$ otherwise. (Note that 
$\wt(v_1)(H)=\wt(v_{\ovl n})(H')$.)
%%% see note IX page 98. %%%%
For any $x\in (\bbC^\times)^{2n}$ 
there exist a unique $y\in (\bbC^\times)^{2n}$ and 
a rational function $a(x)$ 
such that $v_2(y)=a(x)v_1(x)$.
\end{enumerate}
\end{thm}

Now, using this theorem, 
we define the rational mapping
\begin{equation}
\begin{array}{cccccccccc}
\ovl\sigma\cl \cV(\ge)_l
&\longrightarrow &\cV(\ge)_l,&&&
\ovl\sigma\cl\cV(\ge)_l
&\longrightarrow &\cV_2(\ge)_l,&\\
v(x)&\mapsto &v(y)&(\ge\ne A_{2n}^{(2)},
\,\,{\TY(A,2,2n)}^\dagger),&\quad&
v_1(x)&\mapsto &v_2(y)&
(\ge=A_{2n}^{(2)},\,\,
{\TY(A,2,2n)}^\dagger),
\end{array}
\end{equation}
%Put $0^*\seteq \sigma(0)\in I\setminus\{0\}$. 

\begin{thm}
\label{aff-geo}
The rational mapping $\ovl\sigma$ is birational. 
If we define a $\ge_0$-geometric 
crystal structure on $\cV(\ge)_l$ by the one 
on $B^-_{i_1\cd i_k}\cdot l^H$ ($w_1=s_{i_1}\cd s_{i_k}$)
as in \ref{schubert}
and a rational 
$\bbC^\times$-action $e_0:\bbC^\times
\times\cV(\ge)_l\to\cV(\ge)_l$ and rational 
functions $\wt_0$ and $\vep_0$ 
on $\cV(\ge)_l$ by
\begin{equation}
\begin{cases}
e_0^c\seteq \ovl\sigma^{-1}\circ 
e_{\sigma(0)}^c\circ\ovl\sigma
,\q \vep_0\seteq \vep_{\sigma(0)}\circ\ovl\sigma,\q
\gamma_0\seteq \gamma_{\sigma(0)}\circ\ovl\sigma,&
\text{for $\ge\ne A_{2n}^{(2)},\,\,
{\TY(A,2,2n)}^\dagger$},\\
e_0^c\seteq \ovl\sigma^{-1}\circ 
e_0^c\circ\ovl\sigma
,\q \vep_0\seteq \vep_0\circ\ovl\sigma,\q
\gamma_0\seteq \gamma_0\circ\ovl\sigma,&
\text{for $\ge=A_{2n}^{(2)},\,\,
{\TY(A,2,2n)}^\dagger$}.
\end{cases}
\label{ese}
\end{equation}
then $(\cV(\ge)_l,\{e_i\}_{i\in I},
\{\gamma_i\}_{i\in I}, \{\vep_i\}_{i\in I})$
turns out to be an affine $\ge$-geometric crystal.
\end{thm}
{\sl Remark.}
In the case $\ge=A_{2n}^{(2)}$ and 
${\TY(A,2,2n)}^\dagger$, 
%let $\ge_n\subset \ge$ be a simple Lie algebra 
%corresponding to the index $i=0,1,\ld,n-1$.
$\cV_2(\ge)_l$ 
has a $\ge_{I\setminus\{n\}}$-geometric crystal structure. 
Thus, $e_0,\gamma_0,\vep_0$ 
are well-defined on $\cV_2(\ge)_l$.

\medskip
The following lemma shows 
Theorem \ref{aff-geo} partially.
\begin{lem}
\label{io}
Suppose that $\ge\ne
A_{2n}^{(2)}$, ${\TY(A,2,2n)}^\dagger$.
If there exists 
$\osigma$ as above and 
\begin{equation}
e_{\sigma(i)}^c=\osigma\circ e_i^c\circ
\osigma^{-1}, \q\gamma_{i}
=\gamma_{\sigma(i)}\circ\osigma,\q
\vep_{i}=\vep_{\sigma(i)}\circ
\osigma, 
\label{cond-i}
\end{equation}
for $i\ne \sigma^{-1}(0),0$,  then we obtain
\bnum
\item
\begin{eqnarray*}
&&e_0^{c_1}e_i^{c_2}=e_i^{c_2}e_0^{c_1}\q
\text{if }a_{0i}=a_{i0}=0,\\
&&e^{c_1}_{0}e^{c_1c_2}_{i}e^{c_2}_{0}
=e^{c_2}_{i}e^{c_1c_2}_{0}e^{c_2}_{i}
\q\text{if }a_{0i}a_{i0}=1,\\
&&e^{c_1}_{0}e^{c^2_1c_2}_{i}
e^{c_1c_2}_{0}e^{c_2}_{i}
=e^{c_2}_{i}e^{c_1c_2}_{0}e^{c^2_1c_2}_{i}
e^{c_1}_{0}\q
{\rm if }\,\,a_{0i}=-2,\,a_{i0}=-1,\\
&&
e^{c_2}_{0}e^{c_1c_2}_{i}e^{c^2_1c_2}_{0}
e^{c_1}_{i}
=e^{c_1}_{i}e^{c^2_1c_2}_{0}
e^{c_1c_2}_{i}e^{c_2}_{0}
\q
{\rm if }\,\,a_{0i}=-1,\,a_{i0}=-2.
\end{eqnarray*}
\item
$\gamma_0(e_i^c(v(x)))=c^{a_{i0}}
\gamma_0(v(x))$ and 
$\gamma_i(e_0^c(v(x)))=c^{a_{0i}}
\gamma_i(v(x))$.
\item
$\vep_0(e_0^c(v(x)))=c^{-1}\vep_0(v(x))$.
\enum
\end{lem}
{\sl Proof.} 
For example, we have
\begin{eqnarray*}
\gamma_0(e_i^c(v(x)))
&=&\gamma_{\sigma(0)}
(\osigma e_i^c\osigma^{-1}
(\osigma(v(x))))\\
&=&\gamma_{\sigma(0)}
(e_{\sigma(i)}^c(\osigma(v(x))))
=c^{a_{\sigma(i),\sigma(0)}}
\gamma_{\sigma(0)}(\osigma(v(x)))\\
&=&c^{a_{i,0}}\gamma_0(v(x)),
\end{eqnarray*}
where we use
$a_{\sigma(i),\sigma(0)}=a_{i0}$ 
in the last equality. 
The other assertions are obtained similarly.\qed

In the rest of this section, we shall prove
Theorem~\ref{birat} and Theorem~\ref{aff-geo}
in case-by-case methods.

%%%%%%%%%%%%%%%%%%%%%%%%%%%%%%%%%%%%%%%%%%%%%%
\subsection{$A^{(1)}_n$-case $(n\geq2)$}%5.3
\label{geo-a}

We have $w_1\seteq s_ns_{n-1}\cd s_2s_1$, and
\[
{\mathcal V}(\TY(A,1,n))_l
\seteq \{Y_n(x_n)\cd Y_2(x_2)Y_1(x_1)l^H
\fsquare(5mm,1)
\,\big\vert\,x_i\in \bbC^\times\}
\subset W(\varpi_1)_l.
\]
Since $y_i(\frac{1}{c})=
\exp(\frac{f_i}{c})=1+c^{-1}f_i$ on 
$W(\varpi_1)$,  
$v(x)=Y_n(x_n)\cd Y_2(x_2)Y_1(x_1)l^H
\fsquare(5mm,1)$
is explicitly written as
\[
 v(x)=v(x_1,\ld,x_n)=
l^m\left(\sum_{i=1}^nx_i\fsquare(5mm,i)
+\fsquare(6mm,n+1)\right),
\]
where $m=\varpi_1(H)$.
Let  $\sigma\cl\al_k\mapsto \al_{k+1}$ 
$(k\in I)$ be 
the  Dynkin diagram automorphism for $\TY(A,1,n)$, 
which gives rise to the automorphism 
$\sigma_l\cl W(\varpi_1)_l
\rightarrow W(\varpi_1)_l$.
We have
\[
\sigma_l(v_i)=
\begin{cases}v_{i+1}&i\ne n+1,\\
l^{-1}v_1&i=n+1.
\end{cases}
\]
Then, we obtain
\[
 \sigma_l(v(x))=l^m(l^{-1}\fsquare(5mm,1)+
\sum_{i=1}^nx_i\fsquare(5mm,i+1)).
\]
Then the equation $v(y)=a(x)\sigma_l(v(x))$, 
{\em i.e.}
\[
 \sum_{i=1}^ny_i\fsquare(5mm,i)
+\fsquare(6mm,n+1)
=a(x)(l^{-1}\fsquare(5mm,1)
+\sum_{i=1}^nx_i\fsquare(5mm,i+1))
\]
is solved by
\begin{equation}
a(x)=\frac{1}{x_n}, \q y_1=\frac{1}{lx_n},\q
y_i=\frac{x_{i-1}}{x_n} \,\,(i=2,\ld,n), 
\end{equation}
that is
\begin{equation}
\ovl\sigma(v(x_1\cd,x_n))
=v(\frac{1}{lx_n},\frac{x_1}{x_n},\ld,
\frac{x_{n-1}}{x_n}).
\end{equation}
The $A_n$-geometric crystal structure on 
$\cV(\TY(A,1,n))_l$ 
induced from the one on $B^-_{w_1}\cdot l^H$
is given by:
\begin{eqnarray}
&&e_i^c(v(x_1,\ld,x_n))=v(x_1,\ld, cx_i,\ld,x_n)
\q(i=1,\ld,n), \label{ei-A}\\
&&\gamma_1(v(x))=\frac{lx_1^2}{x_2},\q
\gamma_i(v(x))=\frac{x_i^2}{x_{i-1}x_{i+1}}\q
(i=2,\ld,n-1),\q
\gamma_n(v(x))=\frac{x_n^2}{x_{n-1}},\\
&&\vep_i(v(x))=\frac{x_{i+1}}{x_i}\q
(i=1,\ld, n-1),\q
\vep_n(v(x))=\frac{1}{x_n}.
\end{eqnarray}
Then we have 
\[
 \vep_{i+1}(\osigma(v(x)))=
\begin{cases}
\frac{x_{i+1}}{x_{i}}&\text{if }i=1\cd,n-2,\\
\frac{x_n}{x_{n-1}}&\text{if }i=n-1,
\end{cases}
\]
which implies 
$\vep_{\sigma(i)}(\osigma(v(x)))=\vep_i(v(x))$,
and then we completed the 
proof of Theorem \ref{birat} for $\TY(A,1,n)$.

\medskip
Now, we define $e_0^c,\gamma_0$ and $\vep_0$ by 
\begin{equation}
e_0^c\seteq \ovl\sigma^{-1}\circ e_1^c
\circ \ovl\sigma,\q
\gamma_0\seteq \gamma_1\circ\ovl\sigma,\q
\vep_0\seteq \vep_1\circ\ovl\sigma.
\end{equation}
Their explicit forms are
\begin{eqnarray}
&&e_0^c(v(x))=v(\frac{x_1}{c},\frac{x_2}{c},\ld,
\frac{x_n}{c}),\\
&&\gamma_0(v(x))=\frac{1}{lx_1x_n},\q
\vep_0(v(x))=lx_1.
\label{ge0}
\end{eqnarray}
Thus, we can check \eqref{cond-i} easily,
and then Lemma \ref{io} 
reduces the proof of Theorem \ref{aff-geo} to
the statements:
\begin{eqnarray}
&&e^{c_1}_{0}e^{c_1c_2}_{n}e^{c_2}_{0}
=e^{c_2}_{n}e^{c_1c_2}_{0}e^{c_1}_{n},\\
&&\gamma_0(e_n^c(v(x)))=c^{-1}
\gamma_0(v(x)),\q
\gamma_n(e_0^c(v(x)))=c^{-1}
\gamma_n(v(x)).
\end{eqnarray}
These are immediate from (\ref{ei-A})--
(\ref{ge0}). Thus, we obtain 
Theorem \ref{aff-geo} for $\TY(A,1,n)$.

Let us introduce the following 
$\TY(A,1,n)$-geometric crystal 
$\cB_L(\TY(A,1,n))$ $(L\in \bbC^\times)$ 
(\cite{KOTY}):
\begin{eqnarray*}
&&\cB_L(\TY(A,1,n))
\seteq\{l=(l_1,\ld,l_{n},l_{n+1})
\in (\bbC^\times)^{n+1}|
l_1\cd l_{n+1}=L\},\\
&&e_i^c(l)=(\cd,cl_i,c^{-1}l_{i+1},\ld),
\qq\gamma_i(l)=\frac{l_i}{l_{i+1}},\qq
\vep_i(l)=l_{i+1}\q
(i=0,1,\ld,n),
\end{eqnarray*}
where 
we understand $l_0=l_{n+1}$.

We have $\cB_L(\TY(A,1,n))
\cong\cV(\TY(A,1,n))_L$. 
Indeed, defining $\phi:\cB_L(\TY(A,1,n))\to
\cV(\TY(A,1,n))_L$ by 
$\phi(l_1,\cd,l_{n+1})=
v(\frac{l_1}{L},\frac{l_1l_2}{L},\cd,\frac{l_1\cd l_{n}}{L})$,
it is easy to see that $\phi$ is an
isomorphism of geometric crystals.

%%%%%%%%%%%%%%%%%%%%%%%%%%%%%%%%%%%%%%%%%%%%%%
\subsection{$B^{(1)}_n$-case $(n\geq2)$}%5.4
\label{b1n}
We have $w_1=s_1\cd s_{n-1}s_ns_{n-1}\cd s_1$, and
\[
 {\cV}(\TY(B,1,n))_l
\seteq \{v(x)=Y_1(x_1)\cd Y_n(x_n)
Y_{n-1}(\bar x_{n-1})\cd Y_1(\bar x_1)l^H
\fsquare(5mm,1)\,\,\big\vert\,\,x_i,\bar x_i\in\bbC^\times\}.
\]
It follows from the explicit description of 
$W(\varpi_1)_l$
as in \ref{bn-w1}
\begin{eqnarray*}
&&v(x_1,\ld,x_n,\ovl x_{n-1},\ld,\ovl x_1)
=l^m\left\{\left(\sum_{i=1}^n\xi_i(x)\fsquare(0.5cm,i)\right)
+x_n\fsquare(0.5cm,0)
+\left(\sum_{i=2}^{n}x_{i-1}\fsquare(0.5cm,\ovl i)
\right)
+\fsquare(0.5cm,\ovl 1)\right\},\\
&&\qq\q{\rm where}\q
m:=\varpi_1(H)\q \text{ and }
\q \xi_i(x)\seteq \begin{cases}
x_1\ovl x_1&i=1\\
\dfrac{x_{i-1}\ovl x_{i-1}+x_i\ovl x_i}{x_{i-1}}&i\ne1,n\\
\dfrac{x_{n-1}\ovl x_{n-1}+x^2_n}{x_{n-1}}&i=n
\end{cases}
\end{eqnarray*}
The automorphism $\sigma_l
:W(\varpi_1)_l\to W(\varpi_1)_l$ is given as 
\begin{equation}
\sigma_l\fsquare(0.5cm,1)
=l\fsquare(0.5cm,\ovl 1),\q
\sigma_l\fsquare(0.5cm,\ovl 1)
=l^{-1}\fsquare(0.5cm,1), 
\sigma_l\fsquare(0.5cm,k)
=\fsquare(0.5cm,k)\text{ otherwise.}
\label{sigma-l-b}
\end{equation}
Then
we have 
\[
\sigma_l(v(x))
=l^m\left\{l^{-1}\fsquare(0.5cm,1)
+\left(\sum_{i=2}^n \xi_i(x)\fsquare(0.5cm,i)
\right)
+x_n\fsquare(0.5cm,0)
+\left(\sum_{i=2}^{n}x_{i-1}\fsquare(0.5cm,\ovl i)
\right)
+l x_1\ovl x_1\fsquare(0.5cm,\ovl 1)\right\},
\]
The equation $v(y)=a(x)\sigma_l(v(x))$
($x,y\in (\bbC^\times)^{2n-1}$)
has a unique solution:
\begin{equation}
a(x)=\frac{1}{lx_1\ovl x_1},\q
y_i=a(x)x_i=
\frac{x_i}{lx_1\ovl x_1}\q(1\leq i\leq n),\q
\ovl y_i=a(x)\ovl x_i
=\frac{\ovl x_i}{lx_1\ovl x_1}\q
(1\leq i<n).
\end{equation}
Hence we have the rational mapping:
\begin{equation}
\osigma(v(x))\seteq v(y)=
v\left(\frac{x_1}{lx_1\ovl x_1},
\frac{x_2}{lx_1\ovl x_1},\cd,
\frac{x_n}{lx_1\ovl x_1},
\frac{\ovl x_{n-1}}{lx_1\ovl x_1},
\cd,
\frac{\ovl x_1}{lx_1\ovl x_1}\right).
\label{bn-osigma}
\end{equation}
By the explicit form of 
$\osigma$ in (\ref{bn-osigma}), 
we have $\osigma^2={\rm id}$, which means 
that the morphism $\osigma$ is birational.
In this case, the second condition in 
Theorem \ref{birat} is trivial since 
$\sigma(i)=i$ if $i,\,\sigma(i)\ne 0$.
Thus, the proof of 
Theorem \ref{birat} in this case
is completed.

Now, we set 
$e_0^c\seteq \osigma\circ e_1^c\circ\osigma$, 
$\gamma_0\seteq \gamma_1\circ\osigma$ and 
$\vep_0\seteq \vep_1\circ\osigma$.

The explicit forms of $e_i$'s, 
$\vep_i$'s and $\gamma_i$'s are:
\begin{eqnarray*}
e_0^c\cl
&x_1\mapsto&x_1\frac{cx_1\bar x_1+x_2\bar x_2}
{c(x_1\bar x_1+x_2\bar x_2)}\qq
x_i\mapsto\frac{x_i}{c}\,\,(2\leq i\leq n)\\
&\bar x_1\mapsto&\bar x_1\frac{x_1\bar x_1+x_2\bar x_2}
{cx_1\bar x_1+x_2\bar x_2},\qq
\bar x_i\mapsto\frac{\bar x_i}{c}\,\,(2\leq i\leq n-1),\\
\hspace{-20pt}e_i^c\cl&x_i\mapsto&
x_i\frac{cx_i\ovl x_i+x_{i+1}\bar x_{i+1}}
{x_i\ovl x_i+x_{i+1}\ovl x_{i+1}},\q
\ovl x_i\mapsto
\ovl x_i\frac{c(x_i\ovl x_i+x_{i+1}\bar x_{i+1})}
{cx_i\ovl x_i+x_{i+1}\ovl x_{i+1}},\\
&x_j\mapsto& x_j,\q \ovl x_j\mapsto \ovl x_j\,\;
(j\ne i)
\qq\qq(1\leq i<n-1),\\
e_{n-1}^c\cl&x_{n-1}\mapsto&
x_{n-1}\frac{cx_{n-1}\ovl x_{n-1}+x_{n}^2}
{x_{n-1}\ovl x_{n-1}+x_{n}^2},\q
\ovl x_{n-1}\mapsto
\ovl x_{n-1}
\frac{c(x_{n-1}\ovl x_{n-1}+x_{n}^2)}
{cx_{n-1}\ovl x_{n-1}+x_{n}^2},\\
&x_j\mapsto& x_j,\q \ovl x_j\mapsto \ovl x_j\,\,
(j\ne n-1),\\
e^c_n\cl&x_n\mapsto& cx_n,\qq
x_j\mapsto x_j\q \ovl x_j\mapsto \ovl x_j\,\,
(j\ne n),
\end{eqnarray*}
\begin{eqnarray*}
&&\vep_0(v(x))=
\frac{l(x_1\ovl x_1+x_2\ovl x_2)}{x_1},\q
\vep_1(v(x))=\frac{1}{x_1}
\left(1+\frac{x_{2}\ovl x_{2}}{x_{1}\ovl x_1}\right),\\
&&\vep_i(v(x))=\frac{x_{i-1}}{x_i}
\left(1+\frac{x_{i+1}\ovl x_{i+1}}{x_{i}\ovl x_i}\right)\,\,(2\leq i\leq n-2),
\\
&&\vep_{n-1}(v(x))=\frac{x_{n-2}}{x_{n-1}}
\left(1+\frac{x_n^2}{x_{n-1}\ovl x_{n-1}}\right),\,\,
\vep_n(v(x))=\frac{x_{n-1}}{x_n},
\end{eqnarray*}
\begin{eqnarray*}
&&\hspace{-10pt}
\gamma_0(v(x))=\frac{1}{l x_2\ovl x_2},\,\,
\gamma_1(v(x))=\frac{l(x_1\ovl x_1)^2}{x_2\ovl x_2},\,\,
\gamma_i(v(x))=
\frac{(x_i\ovl x_i)^2}{x_{i-1}\ovl x_{i-1}x_{i+1}\ovl x_{i+1}}
\,(2\leq i\leq n-2),\\
&&\gamma_{n-1}(v(x))=
\frac{(x_{n-1}\ovl x_{n-1})^2}{x_{n-2}\ovl x_{n-2}x_n^2},\q
\gamma_n(v(x))=\frac{x_n^2}{x_{n-1}\ovl x_{n-1}}.
\end{eqnarray*}
%where we denote $\gamma_(v(x))$ (resp.~$\vep_i(v(x))$)
%by $\gamma_i(x)$ (resp.~$\vep_i(x)$).
Since $\sigma(i)=i$ for $i\ne0,1$,
the condition (\ref{cond-i}) in Lemma \ref{io}
can be easily seen by (\ref{bn-osigma}) and by the 
explicit form of $e_i$, $\gamma_i$ and $\vep_i$
($i\in I$).
Thus, in order to prove Theorem \ref{aff-geo}, 
it suffices to show that 
\begin{eqnarray}
&&e_0^{c_1}e_1^{c_2}=e_1^{c_2}e_0^{c_1},
\label{b1}\\
&&\gamma_0(e_1^c(v(x)))=\gamma_0(v(x)),\q
\gamma_1(e_0^c(v(x)))=\gamma_1(v(x)).
\label{b2}
\end{eqnarray}
It follows from the explicit formula above that
\[
 e_0^{c_1}e_1^{c_2}(v(x))=e_1^{c_2}e_0^{c_1}(v(x))
=v\left(
x_1\frac{c_1c_2x_1\ovl x_1+x_2\ovl x_2}
{c_1(x_1\ovl x_1+x_2\ovl x_2)},\frac{x_2}{c_1},\cd,
\frac{\ovl x_2}{c_1},
\ovl x_1\frac{c_2(x_1\ovl x_1+x_2\ovl x_2)}
{c_1c_2x_1\ovl x_1+x_2\ovl x_2}
\right),
\]
which implies (\ref{b1}). 
We get  (\ref{b2}) 
immediately from the formula above and we 
complete the proof of Theorem \ref{aff-geo}
for $\TY(B,1,n)$.

%%%%%%%%%%%%%%%%%%%%%%%%%%%%%%%%%%%%%%%%%%%%%%
%%%%%%%%%%%%%%%%%%%%%%%%%%%%%%%%%%%%%%%%%%%%%%
\subsection{$D^{(1)}_n$-case $(n\geq4)$}%5.5
\label{d1n}
We have $w_1=s_1s_2\cd s_{n-1}s_ns_{n-2}s_{n-3}\cd s_2s_1$, and
\[
 {\cV}(\TY(D,1,n))_l\seteq 
\{v(x)=Y_1(x_1)\cd Y_{n-1}(x_{n-1})Y_n(x_n)
Y_{n-2}(\bar x_{n-2})\cd Y_1(\bar x_1)
l^H\fsquare(5mm,1)
\,\,\big\vert\,
\,x_i,\bar x_i\in\bbC^\times\}.
\]
It follows from the explicit form of $\W1_l$ 
in \ref{dn-w1}
that $y_i(c^{-1})=\exp(c^{-1}f_i)=
1+c^{-1}f_i$ on $\W1$. Thus, 
we have 
\begin{eqnarray*}
&&v(x)%=v(x_1,\ld,x_{n-1},x_n,\ovl x_{n-2},\ld,\ovl x_1)
=l^m\left\{\left(\sum_{i=1}^{n-1}\xi_i(x)\fsquare(0.5cm,i)\right)
+x_n\fsquare(0.5cm,n)
+\left(\sum_{i=2}^{n}x_{i-1}
\fsquare(0.5cm,\ovl i)\right)
+\fsquare(0.5cm,\ovl 1)\right\},\\
&&\qq\qq\q{\rm where}\q 
m:=\varpi_1(H),\qq\xi_i(x)\seteq \begin{cases}
x_1\ovl x_1&i=1\\
\frac{x_{i-1}\ovl x_{i-1}+x_i\ovl x_i}{x_{i-1}}&i\ne1,n-1\\
\frac{x_{n-2}\ovl x_{n-2}+x_{n-1}x_n}{x_{n-2}}&i=n-1
\end{cases}
\end{eqnarray*}
The automorphism $\sigma_l
:W(\varpi_1)_l\to W(\varpi_1)_l$ is given as
\[
\sigma_l\fsquare(0.5cm,1)=l\fsquare(0.5cm,\ovl 1),\q
\sigma_l\fsquare(0.5cm,\ovl 1)
=l^{-1}\fsquare(0.5cm,1),\q
\sigma_l\fsquare(0.5cm,k)=\fsquare(0.5cm,k)
\q\text{ otherwise.}
\]
Then we have 
\[
\sigma_l(v(x))
=l^m\left\{l^{-1}\fsquare(0.5cm,1)+
\left(\sum_{i=2}^n \xi_i(x)\fsquare(0.5cm,i)\right)
+x_n\fsquare(0.5cm,n)
+\left(\sum_{i=1}^{n-1}x_{i-1}\fsquare(0.5cm,\ovl i)
\right)
+l \xi_1\fsquare(0.5cm,\ovl 1)\right\}.
\]
Then the equation $v(y)=a(x)\sigma_l(v(x))$
($x,y\in (\bbC^\times)^{2n-2}$)
has the following unique
solution:
\begin{equation}
a(x)=\frac{1}{lx_1\ovl x_1},\q
y_i=a(x)x_i=
\frac{x_i}{lx_1\ovl x_1}\q(1\leq i\leq n),\q
\ovl y_i=a(x)\ovl x_i
=\frac{\ovl x_i}{lx_1\ovl x_1}\q
(1\leq i\leq n-2).
\end{equation}
We define the rational mapping $\osigma\cl
\cV(\TY(D,1,n))_l\longrightarrow 
\cV(\TY(D,1,n))_l$ by 
\begin{equation}
\osigma(v(x))=
v\left(\frac{x_1}{lx_1\ovl x_1},
\frac{x_2}{lx_1\ovl x_1},\cd,
\frac{x_n}{lx_1\ovl x_1},
\frac{\ovl x_{n-2}}{lx_1\ovl x_1},
\cd,
\frac{\ovl x_1}{lx_1\ovl x_1}\right).
\label{dn-osigma}
\end{equation}
It is immediate from 
(\ref{dn-osigma}) that 
$\osigma^2={\rm id}$, which implies
the birationality of 
the morphism $\osigma$.
In this case, the second condition in 
Theorem \ref{birat} is trivial since 
$\sigma(i)=i$ if $i,\,\sigma(i)\ne 0$.
Thus, the proof of 
Theorem \ref{birat} for $\TY(D,1,n)$
is completed.

Now, we set 
$e_0^c\seteq \osigma\circ e_1^c\circ\osigma$, 
$\gamma_0\seteq \gamma_1\circ\osigma$ and 
$\vep_0\seteq \vep_1\circ\osigma$.
The explicit forms of $e_i$, $\vep_i$ and $\gamma_i$ are:
\begin{eqnarray*}
e_0^c\cl
&x_1\mapsto&x_1\frac{cx_1\bar x_1+x_2\bar x_2}
{c(x_1\bar x_1+x_2\bar x_2)}\qq
x_i\mapsto\frac{x_i}{c}\,\,(2\leq i\leq n)\\
&\bar x_1\mapsto&\bar x_1\frac{x_1\bar x_1+x_2\bar x_2}
{cx_1\bar x_1+x_2\bar x_2},\qq
\bar x_i\mapsto\frac{\bar x_i}{c}\,\,(2\leq i\leq n-2),\\
e_i^c\cl&x_i\mapsto&
x_i\frac{cx_i\ovl x_i+x_{i+1}\bar x_{i+1}}
{x_i\ovl x_i+x_{i+1}\ovl x_{i+1}}\q
\ovl x_i\mapsto
\ovl x_i\frac{c(x_i\ovl x_i+x_{i+1}\bar x_{i+1})}
{cx_i\ovl x_i+x_{i+1}\ovl x_{i+1}},\\
&x_j\mapsto& x_j,\q \ovl x_j\mapsto\ovl x_j\q(j\ne i),
\qq\q(1\leq i\leq n-3),\\
e_{n-2}^c\cl&x_{n-2}\mapsto&
x_{n-2}\frac{cx_{n-2}\ovl x_{n-2}+x_{n-1}x_{n}}
{x_{n-2}\ovl x_{n-2}+x_{n-1}x_{n}},\q
\ovl x_{n-2}\mapsto
\ovl x_{n-2}\frac{c(x_{n-2}\ovl x_{n-2}+x_{n-1}x_{n})}
{cx_{n-2}\ovl x_{n-2}+x_{n-1}x_{n}},\\
&x_j\mapsto& x_j,\q \ovl x_j\mapsto\ovl x_j\q(j\ne n-2),\\
e^c_{n-1}\cl&x_{n-1}\mapsto& cx_{n-1},\q
x_j\mapsto x_j,\q \ovl x_j\mapsto\ovl x_j\q(j\ne n-1),\\
e^c_n\cl&x_n\mapsto& cx_n,\q
x_j\mapsto x_j,\q \ovl x_j\mapsto\ovl x_j\q(j\ne n),
\end{eqnarray*}
\begin{eqnarray*}
&&\vep_0(v(x))=
\frac{l(x_1\ovl x_1+x_2\ovl x_2)}{x_1},\q
\vep_1(v(x))=\frac{1}{x_1}
\left(1+\frac{x_{2}\ovl x_{2}}{x_{1}\ovl x_1}\right),\q\\
&&\vep_i(v(x))=\frac{x_{i-1}}{x_i}
\left(1+\frac{x_{i+1}\ovl x_{i+1}}{x_{i}\ovl x_i}\right)\q
(2\leq i\leq n-3),\\
&&\vep_{n-2}(v(x))=\frac{x_{n-3}}{x_{n-2}}
\left(1+\frac{x_{n-1}x_n}{x_{n-2}\ovl x_{n-2}}\right),\q
\vep_{n-1}(v(x))=\frac{x_{n-2}}{x_{n-1}},\q
\vep_n(v(x))=\frac{x_{n-2}}{x_n},
\end{eqnarray*}
\begin{eqnarray*}
&&\hspace{-10pt}
\gamma_0(v(x))=\frac{1}{lx_2\ovl x_2},\,\,
\gamma_1(v(x))=\frac{l(x_1\ovl x_1)^2}{x_2\ovl x_2},\,\,
\gamma_i(v(x))=
\frac{(x_i\ovl x_i)^2}{x_{i-1}\ovl x_{i-1}x_{i+1}\ovl x_{i+1}}
\,(2\leq i\leq n-3),\\
&&\gamma_{n-2}(v(x))=\frac{(x_{n-2}\ovl x_{n-2})^2}
{x_{n-3}\ovl x_{n-3}x_{n-1}x_n},\q
\gamma_{n-1}(v(x))=
\frac{x_{n-1}^2}{x_{n-2}\ovl x_{n-2}},\q
\gamma_n(v(x))=\frac{x_n^2}{x_{n-2}\ovl x_{n-2}}.
\end{eqnarray*}
By these formulas, we can show Theorem \ref{aff-geo} 
for $\TY(D,1,n)$ similarly to the one for $\TY(B,1,n)$.

%%%%%%%%%%%%%%%%%%%%%%%%%%%%%%%%%%%%%%%%%%%%%%%%%%%%%%%
\subsection{$A^{(2)}_{2n-1}$-case $(n\geq3)$}%5.6
\label{aon}
We have $w_1=s_1s_2\cd s_ns_{n-1}\cd s_2s_1$, and
\[
{\cV}(\TY(A,2,2n-1))_l\seteq 
\{v(x)\seteq Y_1(x_1)\cd Y_n(x_n)
Y_{n-1}(\bar x_{n-1})\cd Y_1(\bar x_1)
l^H\fsquare(5mm,1)\,
\,\big\vert\,\,x_i,\bar x_i\in\bbC^\times\}.
\]
%It follows from the explicit form 
%of $\W1$ in \ref{aon-w1}
In this case, $y_i(c^{-1})=\exp(c^{-1}f_i)=
1+c^{-1}f_i$ on $\W1$, and
we have 
\begin{eqnarray*}
&&v(x_1,\ld,x_n,\ovl x_{n-1},\ld,\ovl x_1)\\
&&=l^m\left\{\left(\sum_{i=1}^n\xi_i\fsquare(0.5cm,i)\right)
+\left(\sum_{i=2}^{n}x_{i-1}\fsquare(0.5cm,\ovl i)
\right)
+\fsquare(0.5cm,\ovl 1)\right\},
\q{\rm where}\q m:=\varpi_1(H),\q
\xi_i\seteq \begin{cases}
x_1\ovl x_1&i=1\\
\frac{x_{i-1}\ovl x_{i-1}+x_i\ovl x_i}{x_{i-1}}&i\ne1,n\\
\frac{x_{n-1}\ovl x_{n-1}+x_n}{x_{n-1}}&i=n
\end{cases}
\end{eqnarray*}
The automorphism 
$\sigma_l:W(\varpi_1)_l\to W(\varpi_1)_l$ is given as 
\[
\sigma_l\fsquare(0.5cm,1)=l\fsquare(0.5cm,\ovl 1),\q
\sigma_l\fsquare(0.5cm,\ovl 1)=l^{-1}
\fsquare(0.5cm,1), 
\sigma_l\fsquare(0.5cm,k)=\fsquare(0.5cm,k)
\q\text{otherwise}.
\]
Then we have 
\[
\sigma_l(v(x))
=l^m\left\{l^{-1}\fsquare(0.5cm,1)+\left(
\sum_{i=2}^n \xi_i\fsquare(0.5cm,i)\right)
+\left(\sum_{i=2}^{n}x_{i-1}\fsquare(0.5cm,\ovl i)
\right)
+l x_1\ovl x_1\fsquare(0.5cm,\ovl 1)\right\},
\]
Solving $v(y)=a(x)\sigma_l(v(x))$
($x, y\in (\bbC^\times)^{2n-1}$), we obtain 
a unique solution:
\begin{equation}
a(x)=\frac{1}{l x_1\ovl x_1},\,\,
y_i=a(x)x_i=
\frac{x_i}{lx_1\ovl x_1},\,\,
\ovl y_i=a(x)\ovl x_i
=\frac{\ovl x_i}{l x_1\ovl x_1}\,\,
(1\leq i\leq n-1),\ 
y_n=a(x)^2x_n=\frac{x_n}{(lx_1\ovl x_1)^2}.
\end{equation}
Here we have 
\begin{equation}
\osigma(v(x))=
v\left(\frac{x_1}{lx_1\ovl x_1},
\frac{x_2}{lx_1\ovl x_1},\cd,
\frac{x_n}{(lx_1\ovl x_1)^2},
\frac{\ovl x_{n-1}}{x_1\ovl x_1},
\cd,
\frac{\ovl x_1}{x_1\ovl x_1}\right).
\label{a2o-osigma}
\end{equation}
By the explicit form of 
$\osigma$ in (\ref{a2o-osigma}), 
we have $\osigma^2={\rm id}$, which means 
that the morphism $\osigma$ is birational.
In this case, the second condition in 
Theorem \ref{birat} is trivial since 
$\sigma(i)=i$ if $i,\,\sigma(i)\ne 0$.
Thus, the proof of 
Theorem \ref{birat} for $\TY(A,2,2n-1)$
is completed.

Now, we set 
$e_0^c\seteq \osigma\circ e_1^c\circ\osigma$, 
$\gamma_0\seteq \gamma_1\circ\osigma$ and 
$\vep_0\seteq \vep_1\circ\osigma$.
The explicit forms of $e_i$, $\vep_i$ and $\gamma_i$ are:
\begin{eqnarray*}
e_0^c\cl
&x_1\mapsto&x_1
\frac{cx_1\bar x_1+x_2\bar x_2}
{c(x_1\bar x_1+x_2\bar x_2)}\qq
x_i\mapsto\frac{x_i}{c}\,\,
(2\leq i\leq n-1),\q
x_n\mapsto \frac{x_n}{c^2}\\
&\bar x_1\mapsto&\bar x_1\frac{x_1\bar x_1+x_2\bar x_2}
{cx_1\bar x_1+x_2\bar x_2},\qq
\bar x_i\mapsto\frac{\bar x_i}{c}\,\,(2\leq i\leq n-1),\\
e_i^c\cl&x_i\mapsto&
x_i\frac{cx_i\ovl x_i+x_{i+1}\bar x_{i+1}}
{x_i\ovl x_i+x_{i+1}\ovl x_{i+1}}\qq
\ovl x_i\mapsto
\ovl x_i\frac{c(x_i\ovl x_i+x_{i+1}\bar x_{i+1})}
{cx_i\ovl x_i+x_{i+1}\ovl x_{i+1}},\\
&x_j\mapsto& x_j,\q \ovl x_j\mapsto \ovl x_j\q(j\ne i)
\qq(1\leq i<n-1),\\
e_{n-1}^c\cl&x_{n-1}\mapsto&
x_{n-1}\frac{cx_{n-1}\ovl x_{n-1}+x_{n}}
{x_{n-1}\ovl x_{n-1}+x_{n}}\qq
\ovl x_{n-1}\mapsto
\ovl x_{n-1}\frac{c(x_{n-1}\ovl x_{n-1}+x_{n})}
{cx_{n-1}\ovl x_{n-1}+x_{n}},\\
&x_j\mapsto& x_j,\q \ovl x_j\mapsto \ovl x_j\q(j\ne n-1),\\
e^c_n\cl&x_n\mapsto& cx_n,\q
x_j\mapsto x_j,\q \ovl x_j\mapsto \ovl x_j\q(j\ne n).
\end{eqnarray*}
\begin{eqnarray*}
&&\vep_0(v(x))=
\frac{l(x_1\ovl x_1+x_2\ovl x_2)}{x_1},\q
\vep_1(v(x))=\frac{1}{x_1}
\left(1+\frac{x_{2}\ovl x_{2}}{x_{1}\ovl x_1}\right),\\
&&\vep_i(v(x))=\frac{x_{i-1}}{x_i}
\left(1+\frac{x_{i+1}\ovl x_{i+1}}{x_{i}\ovl x_i}\right)\q
(2\leq i\leq n-2),\\
&& \vep_{n-1}(v(x))=\frac{x_{n-2}}{x_{n-1}}
\left(1+\frac{x_n}{x_{n-1}\ovl x_{n-1}}\right),\q
\vep_n(v(x))=\frac{x_{n-1}^2}{x_n},
\end{eqnarray*}
\begin{eqnarray*}
&&\gamma_0(v(x))=\frac{1}{l x_2\ovl x_2},\q
\gamma_1(v(x))
=\frac{l(x_1\ovl x_1)^2}{x_2\ovl x_2},\q
\gamma_i=
\frac{(x_i\ovl x_i)^2}{x_{i-1}\ovl x_{i-1}x_{i+1}\ovl x_{i+1}}
\q(2\leq i\leq n-2),\\
&&\gamma_{n-1}(v(x))=
\frac{(x_{n-1}\ovl x_{n-1})^2}{x_{n-2}\ovl x_{n-2}x_n},\q
\gamma_n(v(x))=\frac{x_n^2}{(x_{n-1}\ovl x_{n-1})^2}.
\end{eqnarray*}
We can show Theorem \ref{aff-geo} for $\TY(A,2,2n-1)$ 
similarly to the one for $\TY(B,1,n)$.
%%%%%%%%%%%%%%%%%%%%%%%%%%%%%%%%%%%%%%%%%%%%%%%
\subsection{$D^{(2)}_{n+1}$-case $(n\geq2)$}%5.7
\label{d2n}
%see note \Tropical R pp22-25
We have $w_1=s_0s_1\cd s_ns_{n-1}\cd s_2s_1$, and
\[
\hspace{-20pt}
 {\cV}(\TY(D,2,n+1))_l
\seteq \{v(x)\seteq Y_0(x_0)Y_1(x_1)\cd Y_n(x_n)
Y_{n-1}(\bar x_{n-1})\cd Y_1(\bar x_1)l^H
\fsquare(5mm,1)\,\,\big\vert\,\,x_i,\bar x_i
\in\bbC^\times\}.
\]
It follows from the explicit 
description of $W(\varpi_1)_l$
as in \ref{d2n-w1}
that 
on $W(\varpi_1)_l$:
\[
y_i(c^{-1})= \exp(c^{-1}f_i)=\begin{cases}
1+c^{-1}f_i&i\ne 0,n,\\
1+c^{-1}f_i+\frac{1}{2c^2}f_i^2&i=0, n.
\end{cases}
\]
Then we have 
\[
v(x)=l^m\left\{
\left(\sum_{i=1}^n\xi_i(x)\fsquare(0.5cm,i)\right)
+x_n\fsquare(0.5cm,0)+l^{-1}x_0\phi
+\left(\sum_{i=2}^{n}x_{i-1}\fsquare(0.5cm,\ovl i)
\right)
+x_0^2\fsquare(0.5cm,\ovl 1)\right\}
\]
where 
\[
m\seteq \varpi_1(H),\q\qq
\xi_i(x)\seteq \begin{cases}
\frac{l^{-2}x_{0}^2+x_1\ovl x_1}{x_{0}^2}&i=1\\
\frac{x_{i-1}\ovl x_{i-1}+x_i\ovl x_i}{x_{i-1}}&i\ne1,n\\
\frac{x_{n-1}\ovl x_{n-1}+x_n^2}{x_{n-1}}&i=n.
\end{cases}
\]
The automorphism $\sigma_l:W(\varpi_1)_l\to
W(\varpi_1)_l$ is given as 
\begin{eqnarray*}
&&\sigma_l\fsquare(5mm,1)=
l\fsquare(5mm,\ovl n), \q
\sigma_l\fsquare(5mm,\ovl 1)=
l^{-1}\fsquare(5mm,n),\q
\sigma_l\fsquare(5mm,n)=
l\fsquare(5mm,\ovl 1),\q
\sigma_l\fsquare(5mm,\ovl n)=
l^{-1}\fsquare(5mm,1),\\
&&\sigma_l\fsquare(5mm,0)=\phi,\q
\sigma_l\phi=\fsquare(5mm,0),\q
\sigma_l\fsquare(5mm,i+1)
=\fsquare(5mm,\ovl{n-i}),\q
\sigma_l\fsquare(5mm,\ovl{i+1})
=\fsquare(5mm,{n-i})\q
(1\leq i<n-1).
\end{eqnarray*}
Then we have 
\begin{eqnarray*}
&&\sigma_l(v(x))\\
&&\q=l^m\left\{l^{-1}x_{n-1}v_1+
\left(\sum_{i=2}^{n-1} x_{n-i}\fsquare(0.5cm,i)
\right)+l^{-1}x_0^2\fsquare(0.5cm,n)
+l \xi_n v_{\ovl n}
+\left(\sum_{i=2}^{n-1}
\xi_{n-i+1}\fsquare(0.5cm,\ovl i)
\right)+l\xi_1v_{\ovl n}
+x_n\phi+x_0\fsquare(0.5cm,0)\right\},
\end{eqnarray*}
Solving $v(y)=a(x)\sigma_l(v(x))$
($x, y\in (\bbC^\times)^{2n}$), we get a unique solution:
\begin{eqnarray*}
&&a(x)=\frac{x_{n-1}\ovl x_{n-1}+x_n^2}
{lx_{n-1}x_n^2},\\
&&y_0=la(x)x_n=
\frac{x_{n-1}\ovl x_{n-1}+x_n^2}{x_{n-1}x_n},\\
&&y_i=
\frac{(x_{n-i-1}\ovl x_{n-i-1}+x_{n-i}\ovl x_{n-i})
(x_{n-1}\ovl x_{n-1}+x_n^2)}
{lx_{n-i-1}x_{n-1}x_n^2}
\q(1\leq i<n),\\
&&y_{n-1}=
\frac{(l^{-2}x_0^2+x_1\ovl x_1)(x_{n-1}\ovl x_{n-1}+x_n^2)}
{x_0^2x_{n-1}x_n^2}
,\\
&&y_n=
\frac{x_0(x_{n-1}\ovl x_{n-1}+x_n^2)}
{l^2x_{n-1}x_n^2},\\
&&\ovl y_i
=\frac{(x_{n-1}\ovl x_{n-1}+x_n^2)x_{n-i-1}x_{n-i}\ovl x_{n-i}}
{l(x_{n-i-1}\ovl x_{n-i-1}+x_{n-i}\ovl x_{n-i})x_{n-1}x_n^2}\q
(1\leq i\leq n-2),\\
&&\ovl y_{n-1}
=\frac{(x_{n-1}\ovl x_{n-1}+x_n^2)x_0^2x_{1}\ovl x_{1}}
{(x_0^2+l^2x_{1}\ovl x_{1})x_{n-1}x_n^2}.
\end{eqnarray*}
Then we have the rational mapping 
$\osigma\cl\cV(\TY(D,2,n+1))_l
\longrightarrow \cV(\TY(D,2,n+1))_l$
defined by $v(x)\mapsto v(y)$.
% Indeed, this is a birational 
%isomorphism since $\osigma^2={\rm id}$.
The explicit forms of 
$\vep_i$ ($1\leq i\leq n$) are 
as follows:
\begin{eqnarray*}
&&\hspace{-15pt}\vep_1(v(x))=\frac{x_0^2}{x_1}
\left(1+\frac{x_2\ovl x_2}{x_{1}\ovl x_1}\right),\q
\vep_n(v(x))=\frac{x_{n-1}}{x_n},\\
&&\hspace{-15pt}
\vep_i(v(x))=\frac{x_{i-1}}{x_i}
\left(1+\frac{x_{i+1}\ovl x_{i+1}}{x_{i}\ovl x_i}\right)\,\,
(2\leq i\leq n-2),\,\,
\vep_{n-1}(v(x))=\frac{x_{n-2}}{x_{n-1}}
\left(1+\frac{x_n^2}{x_{n-1}\ovl x_{n-1}}\right).
\end{eqnarray*}
Then we get easily that 
$\vep_{n-i}(v(y))=\vep_i(v(x))$ ($1\leq i\leq n-1$), which 
finishes the proof of Theorem \ref{birat} for $\TY(D,2,n+1)$.

Let us define $e^c_0\seteq \osigma\circ e^c_n\circ\osigma$
($\osigma^2={\rm id}$), 
$\gamma_0\seteq \gamma_n\circ\osigma$
and $\vep_0\seteq \vep_n\circ\osigma$. 
The explicit forms of $e_i$, $\gamma_i$ and $\vep_0$ are 
\begin{eqnarray*}
e_0^c\cl&x_0\mapsto&
x_0\frac{c^2x_0^2+l^2x_1\bar x_1}
{c(x_0^2+l^2x_1\bar x_1)},\qq
x_i\mapsto{x_i}\frac{c^2x_0^2+l^2x_1\bar x_1}
{c^2(x_0^2+l^2x_1\bar x_1)}\,\,(1\leq i\leq n),\\
&
\bar x_i\mapsto&
{\bar x_i}\frac{c^2x_0^2+l^2x_1\bar x_1}
{c^2(x_0^2+l^2x_1\bar x_1)}\,\,(1\leq i\leq n-1),\\
e_i^c\cl&x_i\mapsto&
x_i\frac{cx_i\ovl x_i+x_{i+1}\bar x_{i+1}}
{x_i\ovl x_i+x_{i+1}\ovl x_{i+1}},\qq
\ovl x_i\mapsto
\ovl x_i\frac{c(x_i\ovl x_i+x_{i+1}\bar x_{i+1})}
{cx_i\ovl x_i+x_{i+1}\ovl x_{i+1}},\\
&x_j\mapsto& x_j,\q \ovl x_j\mapsto \ovl x_j\q(j\ne i),
\qq\qq\qq\q(1\leq i<n-1),\\
e_{n-1}^c\cl&x_{n-1}\mapsto&
x_{n-1}\frac{cx_{n-1}\ovl x_{n-1}+x_{n}^2}
{x_{n-1}\ovl x_{n-1}+x_{n}^2},\q
\ovl x_{n-1}\mapsto
\ovl x_{n-1}\frac{c(x_{n-1}\ovl x_{n-1}+x_{n}^2)}
{cx_{n-1}\ovl x_{n-1}+x_{n}^2},\\
&x_j\mapsto& x_j,\q \ovl x_j\mapsto \ovl x_j\q(j\ne n-1),\\
e^c_n\cl&x_n\mapsto& cx_n,\q
x_j\mapsto x_j,\q \ovl x_j\mapsto \ovl x_j\q(j\ne n),
\end{eqnarray*}
\begin{eqnarray*}
&&\hspace{-15pt}
\gamma_0(v(x))=\frac{x_0^2}{l^2x_1\ovl x_1},\q
\gamma_1(v(x))=\frac{(lx_1\ovl x_1)^2}
{x_0^2x_2\ovl x_2},\,\,
\gamma_i(v(x))=
\frac{(x_i\ovl x_i)^2}{x_{i-1}\ovl x_{i-1}x_{i+1}\ovl x_{i+1}}\,\,
(2\leq i\leq n-2),\\
&&\hspace{-15pt}\gamma_{n-1}(v(x))=
\frac{(x_{n-1}\ovl x_{n-1})^2}{x_{n-2}\ovl x_{n-2}x_n^2},\,\,
\gamma_n(v(x))=\frac{x_n^2}{x_{n-1}\ovl x_{n-1}},\\
&&\hspace{-15pt}\vep_0(v(x))=
\frac{x_0^2+l^2x_1\ovl x_1}{x_0^3}.
\end{eqnarray*}
Let us check the condition (\ref{cond-i}) 
in Lemma \ref{io}. 
The following are useful for this purpose:
\begin{eqnarray}
&&y_i\ovl y_i=a(x)^2x_{n-i}\ovl x_{n-i},\q
y_0=a(x)x_n,\label{u1}\\
&&a(v(y))=a(\osigma(v(x)))=\frac{1}{a(v(x))}.
\label{u2}
\end{eqnarray}
Using these we can easily check 
the two conditions 
$\gamma_{i}
=\gamma_{\sigma(i)}\circ\osigma$ and 
$\vep_{i}=\vep_{\sigma(i)}\circ
\osigma$.
The condition 
$e_{\sigma(i)}^c=\osigma\circ e_i^c\circ
\osigma^{-1}$ for $i=2,\cd,n-2$ is also 
immediate from (\ref{u1}) and (\ref{u2}).
Next let us see the case $i=1,n-1$.
We have 
\[
a(e_{n-1}^c(v(y)))=
\frac{y_n^2+cy_{n-1}\ovl y_{n-1}}{y^2_ny_{n-1}
\frac{cy_{n-1}\ovl y_{n-1}{y}+y_n^2}
{y_{n-1}\ovl y_{n-1}{y}+y_n^2}}
=\frac{y_{n-1}\ovl y_{n-1}{y}+y_n^2}{y_{n-1}y_n^2}
=\frac{1}{a(v(x))}
\]
Using this, we can get 
$e_{n-i}^c=\osigma\circ e_1^c\circ
\osigma^{-1}$ and then 
$e_{1}^c=\osigma\circ e_{n-1}^c\circ
\osigma^{-1}$ since $\osigma^2={\rm id}$.
Now, it remains to show that 
\begin{eqnarray*}
&&e_0^{c_1}e_n^{c_2}=e_n^{c_2}e_0^{c_1},\q
\vep_0(e_n^c(v(x)))=\vep_0(v(x)),\q
\vep_n(e_0^c(v(x)))=\vep_n(v(x)).
\end{eqnarray*}
They easily follow from the explicit form of 
$e_0^c$. Thus, the proof of Theorem \ref{aff-geo}
in this case is completed.

%%%%%%%%%%%%%%%%%%%%%%%%%%%%%%%
\subsection{$A^{(2)}_{2n}$-case $(n\geq2)$}%5.8
\label{aen-1}

As in the beginning of this section, 
we have $w_1=s_0s_1\cd s_ns_{n-1}\cd s_2s_1$
and
\[
 {\mathcal V}(A_{2n}^{(2)})_l
\seteq \{v_1(x)
=Y_0(x_0)Y_1(x_1)\cd Y_n(x_n)
Y_{n-1}(\bar x_{n-1})\cd Y_1(\bar x_1)
l^H\fsquare(5mm,1)
\,\,\big\vert\,\,x_i,\bar x_i\in\bbC^\times\}.
\]
By the explicit description of 
$W(\varpi_1)_l$
as in \ref{aen-v1}, 
on $W(\varpi_1)_l$ we have:
\[
y_i(c^{-1})=\exp(c^{-1}f_i)=\begin{cases}
1+c^{-1}f_i&i\ne 0,\\
1+c^{-1}f_0+\frac{1}{2c^2}f_0^2&i= 0.
\end{cases}
\]
Then we have 
\begin{eqnarray*}
&&v_1(x)
=l^m\left\{\left(\sum_{i=1}^{n}\xi_i(x)\fsquare(0.5cm,i)\right)
+x_0^2v_{\ovl 1}
+\left(\sum_{i=2}^{n}x_{i-1}
\fsquare(0.5cm,\ovl i)+
\right)+\frac{x_0}{l}\emptyset\right\},\\
&&{\rm where}\q 
m:=\varpi_1(H), \q 
\xi_i(x)\seteq \begin{cases}
\frac{l^{-2}x_{0}^2+x_1\ovl x_1}{x_{0}^2}
&i=1\\
\frac{x_{i-1}\ovl x_{i-1}+x_i\ovl x_i}{x_{i-1}}&i\ne1,n\\
\frac{x_{n-1}\ovl x_{n-1}+x_n}{x_{n-1}}&i=n.
\end{cases}
\end{eqnarray*}

Next, for $w_2=s_ns_{n-1}\cd s_1
s_0s_1\cd s_{n-1}$, we set 
\[
\cV_2(A_{2n}^{(2)})_l\seteq \{v_2(y)
=Y_n(y_n)\cd Y_1(y_1)Y_0(y_0)Y_{1}(\ovl y_1)\cd
Y_{n-1}(\ovl y_{n-1})
l^{H'}\fsquare(5mm,\ovl n)
\,\big\vert\,y_i,\ovl y_i\in\bbC^\times\}.
\]
Then we have 
\[
v_2(y)
=l^m\left\{\left(\sum_{i=1}^{n} 
\frac{y_i}{l^2}
\fsquare(0.5cm,i)\right)
+\left(\sum_{i=1}^{n}
\eta_i(y)\fsquare(0.5cm,\ovl i)
\right)+\frac{y_0}{l}\emptyset\right\}
\ {\rm where}\ \eta_i(y)\seteq \begin{cases}
\frac{y_{0}^2+y_1\ovl y_1}{y_{1}}&i=1\\
\frac{y_{i-1}\ovl y_{i-1}+y_i\ovl y_i}{y_{i}}&i\ne1,n\\
\frac{y_{n-1}\ovl y_{n-1}+l^{-2}y_n}
{y_{n}}&i=n.
\end{cases}
\]
Note that 
$\varpi_1(H)=
m=\wt(v_{\ovl n})(H')=(\varpi_n(H'))$.
For $x\in(\bbC^\times)^{2n}$ 
there exist a unique
$y=(y_0,\ld,\ovl y_1)\in (\bbC^\times)^{2n}$ 
and $a(x)$ such 
that $v_2(y)=a(x)v_1(x)$.
They are given by
\begin{eqnarray*}
&&a(x)=\frac{x_{n-1}\ovl x_{n-1}+x_n}
{l^2 x_{n-1}x_n},\\
&&y_0=a(x)x_0=
\frac{x_0(x_{n-1}\ovl x_{n-1}
+x_n)}{l^2x_{n-1}x_n},\\
&&y_1=l^2a(x)\xi_1(x)=
\frac{(x_{n-1}\ovl x_{n-1}
+x_{n})(x_0+l^2x_1\ovl x_1)}
{x_0^2x_{n-1}x_n},\\
&&y_i=l^2a(x)\xi_i(x)=
\frac{(x_{i-1}\ovl x_{i-1}
+x_{i}\ovl x_{i})
(x_{n-1}\ovl x_{n-1}+x_n)}
{x_{i-1}x_{n-1}x_n}
\q(1< i<n),\\
&&y_n=l^2a(x)\xi_n=
\frac{(x_{n-1}\ovl x_{n-1}+x_n)^2}
{x_{n-1}^2x_n},\\
&&\ovl y_{1}=a(x)\frac{l^2x_0^2x_1\ovl x_1}
{x_0+l^2x_1\ovl x_1}
=\frac{(x_{n-1}\ovl x_{n-1}
+x_n)x_0^2x_{1}\ovl x_{1}}
{(x_0+l^2 x_{1}\ovl x_{1})x_{n-1}x_n},\\
&&\ovl y_i=a(x)\frac{x_{i-1}x_i\ovl x_i}
{x_{i-1}\ovl x_{i-1}+x_i\ovl x_i}
=\frac{(x_{n-1}\ovl x_{n-1}+x_n^2)
x_{i-1}x_{i}\ovl x_{i}}
{l^2(x_{i-1}\ovl x_{i-1}+x_{i}\ovl x_{i})x_{n-1}x_n^2}\q
(1<i\leq n-1).
\end{eqnarray*}
It defines a rational mapping
$\osigma\cl\cV(A_{2n}^{(2)})_l
\longrightarrow 
\cV_2(A_{2n}^{(2)})_l$
($v_1(x)\mapsto v_2(y)$).\\
The inverse 
$\osigma^{-1}\cl\cV_2(A_{2n}^{(2)})_l
\longrightarrow \cV(A_{2n}^{(2)})_l$
$(v_2(y)\mapsto v_1(x))$ is given by 
\begin{eqnarray*}
&&a(y)\seteq \frac{y_0^2y_1}
{y_0^2+y_1\ovl y_1}(=a(x)),\\
&&x_0=\frac{y_0}{a(y)},\\
&&x_i=a(y)^{-1}\frac{y_i\ovl y_i
+y_{i+1}\ovl y_{i+1}}{y_{i+1}}\q(1\leq i\leq n-2),
\\
&& x_{n-1}=a(y)^{-1}
\frac{y_{n-1}\ovl y_{n-1}+l^{-2}y_n}{y_n},\\
&&x_n=\frac{y_n}{l^4a(y)^2},\\
&&\ovl x_i=(l^2a(y))^{-1}\frac{y_i\ovl y_iy_{i+1}}
{y_i\ovl y_i+y_{i+1}\ovl y_{i+1}}\q
(1\leq i\leq n-2),\\
&&\ovl x_{n-1}=\frac{y_{n-1}\ovl y_{n-1}
y_n}{a(y)(l^2y_{n-1}\ovl y_{n-1}+y_n)},
\end{eqnarray*}
which means that the morphism $\osigma\cl
\cV(A_{2n}^{(2)})_l\longrightarrow 
\cV_2(A_{2n}^{(2)})_l$ is birational.
Thus, we obtain Theorem \ref{birat}(ii). 

The actions
of $e_i$ $(0\le i<n)$ on $v_2(y)$
are induced from the ones on 
$Y_{\bf i_2}(y)\cdot l^{H'}\seteq Y_n(y_n)
\cd Y_{n-1}(\ovl y_{n-1})\cdot l^{H'}$
(${\bf i_2}=(n,\ld,1,0,1,$ $\cd,n-1)$)
since $e_i\fsquare(0.5cm,\ovl n)=0$ for $i=0,1,\ld,n-1$.
We also get $\gamma_i(v_2(y))$ and 
$\vep_i(v_2(y))$ from the ones for 
$Y_{\bf i_2}(y)\cdot l^{H'}$ 
where $v_2(y)=\osigma(v_1(x))$:
%on $\cV_2(A_{2n}^{(2)\,\dagger})$, 
\begin{eqnarray*}
&&e_0^c\cl y_0\mapsto cy_0,
\q y_i\mapsto y_i,\q
\ovl y_i\mapsto\ovl y_i\,\,(i\ne0),\\
&&e_1^c\cl y_1\mapsto 
y_1\frac{cy_1\ovl y_1+y_0^2}
{y_1\ovl y_1+y_0^2}, \q
\ovl y_1\mapsto \ovl y_1
\frac{c(y_1\ovl y_1+y_0^2)}
{cy_1\ovl y_1+y_0^2},
\q y_i\mapsto y_i,\q
\ovl y_i\mapsto\ovl y_i\,\,(i\ne1),
\\
&&e_{i}^c\cl y_i\mapsto y_i
\frac{cy_i\ovl y_i+y_{i-1}\ovl y_{i-1}}
{y_i\ovl y_i+y_{i-1}\ovl y_{i-1}},\q
\ovl y_i\mapsto \ovl y_i
\frac{c(y_i\ovl y_i+y_{i-1}\ovl y_{i-1})}
{cy_i\ovl y_i+y_{i-1}\ovl y_{i-1}},\\
&&\q y_j\mapsto y_j,\q
\ovl y_j\mapsto\ovl y_j\,\,(j\ne i),
\q(i=2,\ld,n-1),\\
&&\gamma_0(v_2(y))
=\frac{y_0^2}{y_1\ovl y_1},\q
\gamma_1(v_2(y))=
\frac{(y_1\ovl y_1)^2}{y_0^2y_2\ovl y_2},\\
&&\gamma_i(v_2(y))=
\frac{(y_i\ovl y_i)^2}
{y_{i-1}\ovl y_{i-1}y_{i+1}\ovl y_{i+1}}
\q(i=2,\ld, n-2),\q
\gamma_{n-1}(v_2(y))
=\frac{(y_{n-1}\ovl y_{n-1})^2}
{y_{n-2}\ovl y_{n-2}y_{n}},\\
&&\vep_0(v_2(y))=\frac{y_1}{y_0},\,
\vep_1(v_2(y))=\frac{y_2}{y_1}
\left(1+\frac{y_0^2}{y_1\ovl y_1}\right),\,
\vep_{n-1}(v_2(y))=
\frac{y_n}{y_{n-1}}\left(
1+\frac{y_{n-2}\ovl y_{n-2}}
{y_{n-1}\ovl y_{n-1}}\right),\\
&&\vep_i(v_2(y))=\frac{y_{i+1}}{y_i}
\left(1+\frac{y_{i-1}\ovl y_{i-1}}
{y_{i}\ovl y_i}\right)\q
(i=2\cd,n-1).
\end{eqnarray*}

The explicit forms of $\vep_i(v_1(x))$and 
$\gamma_i(v_1(x))$
$(1\leq i\leq n)$ are also induced from 
the ones
for $Y_{\bf i_1}(x)\cdot l^H\seteq 
Y_0(x_0)\cd Y_1(\ovl x_1)\cdot l^H$ and, we define
$\vep_0(v_1(x))\seteq \vep_0(v_2(y))$ and 
$\gamma_0(v_1(x))\seteq \gamma_0(v_2(y))$
$(v_2(y)\seteq \osigma(v_1(x)))$:
\begin{eqnarray*}
&&\vep_0(v_1(x))=
\frac{1}{x_0}
\left(1+\frac{l^2x_1\ovl x_1}{x_0^2}
\right)^2,\q
\vep_i(v_1(x))=\frac{x_{i-1}}{x_i}
\left(1+\frac{x_{i+1}\ovl x_{i+1}}
{x_{i}\ovl x_i}\right)\q
(1<i\leq n-2),\\
&&
\vep_1(v_1(x))=\frac{x_0^2}{x_1}
\left(1+\frac{x_{2}\ovl x_{2}}
{x_{1}\ovl x_1}\right),\q
\vep_{n-1}(v_1(x))=\frac{x_{n-2}}{x_{n-1}}
\left(1+\frac{x_n}{x_{n-1}\ovl x_{n-1}}
\right),\q
\vep_n(v_1(x))=\frac{x_{n-1}^2}{x_n}.
\end{eqnarray*}
\begin{eqnarray*}
&&\hspace{-23pt}\gamma_0(v_1(x))
=\frac{x_0^2}{l^4 x_1\ovl x_1},\,\,
\gamma_1(v_1(x))=\frac{(lx_1\ovl x_1)^2}
{x_0^2x_2\ovl x_2},\,\,
\gamma_i(v_1(x))=
\frac{(x_i\ovl x_i)^2}{x_{i-1}\ovl x_{i-1}x_{i+1}\ovl x_{i+1}}\,\,
(2\leq i\leq n-2),\\
&&\hspace{-23pt}\gamma_{n-1}(v_1(x))=
\frac{(x_{n-1}\ovl x_{n-1})^2}
{x_{n-2}\ovl x_{n-2}x_n},\q
\gamma_n(v_1(x))=\frac{x_n^2}
{(x_{n-1}\ovl x_{n-1})^2}
\end{eqnarray*}
For $i=0$, we define 
$e^c_0(v_1(x))=
\osigma^{-1}\circ e^c_0\circ\osigma(v_1(x))
=\osigma^{-1}\circ e^c_0(v_2(y))$.
Then we get
\begin{eqnarray*}
e_0^c\cl&x_0\mapsto&x_0
\frac{c^2x_0+l^2x_1\bar x_1}
{c(x_0^2+l^2x_1\bar x_1)},\qq
x_i\mapsto{x_i}\frac{c^2x_0+l^2x_1\bar x_1}
{c^2(x_0^2+l^2x_1\bar x_1)}
\,\,(1\leq i<n),
\\
&\bar x_i\mapsto&{\bar x_i}
\frac{c^2x_0+l^2x_1\bar x_1}
{c^2(x_0^2+l^2x_1\bar x_1)}
\,\,(1\leq i\leq n-1),\,\,
x_n\mapsto{x_n}
\frac{(c^2x_0+l^2x_1\bar x_1)^2}
{c^4(x_0^2+l^2x_1\bar x_1)^2},\\
e_i^c\cl&x_i\mapsto&
x_i\frac{cx_i\ovl x_i+x_{i+1}\bar x_{i+1}}
{x_i\ovl x_i+x_{i+1}\ovl x_{i+1}},\q
\ovl x_i\mapsto
\ovl x_i\frac{c(x_i\ovl x_i+x_{i+1}\bar x_{i+1})}
{cx_i\ovl x_i+x_{i+1}\ovl x_{i+1}},\\
&&\q\q x_j\mapsto x_j,\q
\ovl x_j\mapsto\ovl x_j\,\,(j\ne i)\q
(1\leq i<n-1)\\
e_{n-1}^c\cl&x_{n-1}\mapsto&
x_{n-1}\frac{cx_{n-1}\ovl x_{n-1}+x_{n}}
{x_{n-1}\ovl x_{n-1}+x_{n}},\qq
\ovl x_{n-1}\mapsto
\ovl x_{n-1}
\frac{c(x_{n-1}\ovl x_{n-1}+x_{n})}
{cx_{n-1}\ovl x_{n-1}+x_{n}}\q,\\
&&\q\q x_i\mapsto x_i,\q
\ovl x_i\mapsto\ovl x_i\,\,(i\ne  n-1),\\
e^c_n\cl&x_n\mapsto& cx_n, 
\q\q x_i\mapsto x_i,\q
\ovl x_i\mapsto\ovl x_i\,\,(i\ne n).
\end{eqnarray*}
In order to prove Theorem \ref{aff-geo}, it 
suffices to show the following:
\begin{eqnarray}
&&
e_i^c=\osigma^{-1}\circ 
e_{i}^c\circ \osigma,\,\,
\gamma_i=\gamma_{i}\circ \osigma,\,\,
\vep_i=\vep_{i}\circ \osigma
\,\,(i\ne0,n),
\label{0i}\\
&&e_0^{c_1}e_n^{c_2}=e_n^{c_2}e_0^{c_1},
\\
&&\gamma_n(e_0^c(v_1(x)))=\gamma_n(v_1(x)),\q
\gamma_0(e_n^c(v_1(x)))=\gamma_0(v_1(x)),
\\
&&\vep_0(e_0^c(v_1(x)))=c^{-1}\vep_0(v_1(x)),
\end{eqnarray}
which are immediate from the above formulae.
Let us show \eqref{0i}.
Set $v_2(y)\seteq \osigma(v_1(x))$ and $v_1(x')
\seteq \osigma^{-1}(e_i^c(v_2(y)))$
for $i=2,\ld,n-2$. Then
$x_j'=x_j$ and $\ovl x_j'=\ovl x_j$ for $j\ne i-1,i$, and we have 
\begin{eqnarray*}
&&a(v_2(y))=a(v_1(x)),\\
&&x_i'=\frac{1}{a(v_2(y))}\left(
\ovl y_{i+1}+\frac{cy_i\ovl y_i}{y_{i+1}}\right)
=x_i\frac{cx_i\ovl x_i+x_{i+1}\ovl x_{i+1}}
{x_i\ovl x_i+x_{i+1}\ovl x_{i+1}},\\
&&\ovl x_i'=\frac{1}{a(v_2(y))}
\frac{cy_i\ovl y_iy_{i+1}}
{cy_i\ovl y_i+y_{i+1}\ovl y_{i+1}}
=\ovl x_i\frac{c(x_i\ovl x_i
+x_{i+1}\ovl x_{i+1})}
{cx_i\ovl x_i+x_{i+1}\ovl x_{i+1}},\\
&&x_{i-1}'=\frac{1}{a(v_2(y))}\left(
\ovl y_{i}+\frac{cy_{i-1}\ovl y_{i-1}}
{y_{i}}\right)
=x_{i-1},\\
&&\ovl x_{i-1}'=\frac{1}{a(v_2(y))}
\frac{cy_{i-1}\ovl y_{i-1}y_{i}}
{cy_i\ovl y_i+y_{i-1}\ovl y_{i-1}}
=\ovl x_{i-1},
\end{eqnarray*}
where the formula
$y_i\ovl y_i=a(v_1(x))x_i\ovl x_i$
is useful to obtain these results.
Therefore we have $e_i^c=\osigma^{-1}\circ 
e_i^c\circ\osigma$ for $i=2,\ld,n-2$. 
Others are obtained similarly.

%%%%%%%%%%%%%%%%%%%%%%%%%%%%%%%%%%%%%%%
\subsection{$A^{(2)\,\dagger}_{2n}$-case
$(n\geq2)$}%5.9
\label{aend-1}
%\cmt{\bf Caution!!Fix here!}
%In this subsection, we shall realize
%$\cV(A_{2n}^{(2)\,\dagger})$ in $W(\varpi_1)$ as in
%\ref{aen-v1}. 
%In this case, 
%$W(\varpi_1)$ is isomorphic to 
%$V(\Lm_1)$ as a $U_q(B_n)$-module.

As in the beginning of this section, 
we have $w_1=s_0s_1\cd s_ns_{n-1}\cd s_2s_1$
and
\[
 {\mathcal V}(A_{2n}^{(2)\,\dagger})_l
\seteq \{v_1(x)
=Y_0(x_0)Y_1(x_1)\cd Y_n(x_n)
Y_{n-1}(\bar x_{n-1})\cd Y_1(\bar x_1)
l^H\fsquare(5mm,1)
\,\,\big\vert\,\,x_i,\bar x_i\in\bbC^\times\}.
\]
By the explicit description of $W(\varpi_1)_l$
as in \ref{aen-v1}, 
on $W(\varpi_1)_l$ we have:
\[
y_i(c^{-1})=\exp(c^{-1}f_i)=\begin{cases}
1+c^{-1}f_i&i\ne n,\\
1+c^{-1}f_n+\frac{1}{2c^2}f_n^2&i= n.
\end{cases}
\]
Then we have 
\[
v_1(x)
=l^m\left\{\left(\sum_{i=1}^{n}\xi_i(x)\fsquare(0.5cm,i)\right)
+x_n\fsquare(0.5cm,0)
+\left(\sum_{i=1}^{n}x_{i-1}\fsquare(0.5cm,\ovl i)
\right)\right\}
\q{\rm where}\q \xi_i(x)\seteq \begin{cases}
\frac{l^{-1}x_{0}+x_1\ovl x_1}{x_{0}}&i=1\\
\frac{x_{i-1}\ovl x_{i-1}+x_i\ovl x_i}{x_{i-1}}&i\ne1,n\\
\frac{x_{n-1}\ovl x_{n-1}+x_n^2}{x_{n-1}}&i=n,
\end{cases}
\]
and $m:=\varpi_1(H)$. 
Next, for $w_2=s_ns_{n-1}\cd s_1
s_0s_1\cd s_{n-1}$ and $H'$ 
such that $m=\varpi_n(H')$, we set 
\[
\cV_2(A_{2n}^{(2)\,\dagger})_l\seteq \{v_2(y)
=Y_n(y_n)\cd Y_1(y_1)Y_0(y_0)Y_{1}(\ovl y_1)\cd
Y_{n-1}(\ovl y_{n-1})
l^{H'}\fsquare(5mm,\ovl n)
\,\big\vert\,y_i,\ovl y_i\in\bbC^\times\}.
\]
Then we have 
\[
v_2(y)
=l^m\left\{\left(\sum_{i=1}^{n-1} \frac{y_i}{l}
\fsquare(0.5cm,i)\right)
+\frac{y_n^2}{l}\fsquare(0.5cm,n)
+\frac{y_n}{l}\fsquare(0.5cm,0)
+\left(\sum_{i=1}^{n}
\eta_i(y)\fsquare(0.5cm,\ovl i)
\right)\right\}
\ {\rm where}\ \eta_i(y)\seteq \begin{cases}
\frac{y_{0}+y_1\ovl y_1}{y_{1}}&i=1\\
\frac{y_{i-1}\ovl y_{i-1}+y_i\ovl y_i}{y_{i}}&i\ne1,n\\
\frac{y_{n-1}\ovl y_{n-1}+l^{-1}y_n^2}
{y_{n}^2}&i=n.
\end{cases}
\]
For $x\in(\bbC^\times)^{2n}$ 
there exist a unique
$y=(y_0,\ld,\ovl y_1)\in (\bbC^\times)^{2n}$ 
and $a(x)$ such 
that $v_2(y)=a(x)v_1(x)$.
They are given by
\begin{eqnarray*}
&&a(x)=\frac{x_{n-1}\ovl x_{n-1}+x_n^2}
{l x_{n-1}x_n^2},\\
&&y_0=a(x)^2x_0=
\frac{x_0(x_{n-1}\ovl x_{n-1}
+x_n^2)^2}{(lx_{n-1}x_n^2)^2},\\
&&y_1=la(x)\xi_1(x)=
\frac{(x_{n-1}\ovl x_{n-1}
+x_{n}^2)(l^{-1} x_0+x_1\ovl x_1)}
{x_0x_{n-1}x_n^2},\\
&&y_i=la(x)\xi_i(x)=
\frac{(x_{i-1}\ovl x_{i-1}
+x_{i}\ovl x_{i})
(x_{n-1}\ovl x_{n-1}+x_n^2)}
{x_{i-1}x_{n-1}x_n^2}
\q(1< i<n),\\
&&y_n=la(x)x_n=
\frac{x_{n-1}\ovl x_{n-1}
+x_n^2}{x_{n-1}x_n},\\
&&\ovl y_{1}=a(x)\frac{x_0x_1\ovl x_1}
{l^{-1}x_0+x_1\ovl x_1}
=\frac{(x_{n-1}\ovl x_{n-1}
+x_n^2)x_0x_{1}\ovl x_{1}}
{(x_0+l x_{1}\ovl x_{1})x_{n-1}x_n^2},\\
&&\ovl y_i=a(x)\frac{x_{i-1}x_i\ovl x_i}
{x_{i-1}\ovl x_{i-1}+x_i\ovl x_i}
=\frac{(x_{n-1}\ovl x_{n-1}+x_n^2)
x_{i-1}x_{i}\ovl x_{i}}
{l(x_{i-1}\ovl x_{i-1}+x_{i}\ovl x_{i})x_{n-1}x_n^2}\q
(1<i\leq n-1).
\end{eqnarray*}
It defines a rational mapping
$\osigma\cl\cV(A_{2n}^{(2)\,\dagger})
\longrightarrow 
\cV_2(A_{2n}^{(2)\,\dagger})$
($v_1(x)\mapsto v_2(y)$).\\
The inverse 
$\osigma^{-1}\cl\cV_2(A_{2n}^{(2)\,\dagger})
\longrightarrow \cV(A_{2n}^{(2)\,\dagger})$
$(v_2(y)\mapsto v_1(x))$ is given by 
\begin{eqnarray*}
&&a(y)\seteq \frac{y_0y_1}{y_0+y_1\ovl y_1}(=a(x)),\\
&&x_0=a(y)^{-1}\frac{y_0+y_1\ovl y_1}{y_1},\\
&&x_i=a(y)^{-1}\frac{y_i\ovl y_i
+y_{i+1}\ovl y_{i+1}}{y_{i+1}}\q(1\leq i\leq n-2),
\\
&& x_{n-1}=a(y)^{-1}
\frac{y_{n-1}\ovl y_{n-1}+l^{-1}y_n^2}{y_n^2},\\
&&x_n=\frac{y_n}{la(y)},\\
&&\ovl x_i=(la(y))^{-1}\frac{y_i\ovl y_iy_{i+1}}
{y_i\ovl y_i+y_{i+1}\ovl y_{i+1}}\q
(1\leq i\leq n-2),\\
&&\ovl x_{n-1}=a(y)^{-1}\frac{y_{n-1}\ovl y_{n-1}
y_n^2}{ly_{n-1}\ovl y_{n-1}+y_n^2},
\end{eqnarray*}
which means that the morphism $\osigma\cl
\cV(A_{2n}^{(2)\,\dagger})\longrightarrow 
\cV_2(A_{2n}^{(2)\,\dagger})$ is birational.
Thus, we obtain Theorem \ref{birat}(ii). 

The actions
of $e_i$ $(0\leq i<n)$ on $v_2(y)$
are induced from the ones on 
$Y_{\bf i_2}(y)\cdot l^{H'}\seteq Y_n(y_n)\cd
Y_{n-1}(\ovl y_{n-1})\cdot l^{H'}$
(${\bf i_2}=(n,\ld,1,0,1,$ $\cd,n-1)$)
since $e_i\fsquare(0.5cm,\ovl n)=0$ for $i=0,1,\ld,n-1$.
We also get $\gamma_i(v_2(y))$ and 
$\vep_i(v_2(y))$ from the ones for 
$Y_{\bf i_2}(y)\cdot l^{H'}$ 
where $v_2(y)=\osigma(v_1(x))$:
%on $\cV_2(A_{2n}^{(2)\,\dagger})$, 
\begin{eqnarray*}
&&e_0^c\cl y_0\mapsto cy_0,\q y_i\mapsto y_i,\q
\ovl y_i\mapsto\ovl y_i\,\,(i\ne0),\\
&&e_1^c\cl y_1\mapsto y_1
\frac{cy_1\ovl y_1+y_0}
{y_1\ovl y_1+y_0}, \q
\ovl y_1\mapsto \ovl y_1
\frac{c(y_1\ovl y_1+y_0)}{cy_1\ovl y_1+y_0},
\q y_i\mapsto y_i,\q
\ovl y_i\mapsto\ovl y_i\,\,(i\ne1),
\\
&&e_{i}^c\cl y_i\mapsto y_i
\frac{cy_i\ovl y_i+y_{i-1}\ovl y_{i-1}}
{y_i\ovl y_i+y_{i-1}\ovl y_{i-1}},\q
\ovl y_i\mapsto \ovl y_i
\frac{c(y_i\ovl y_i+y_{i-1}\ovl y_{i-1})}
{cy_i\ovl y_i+y_{i-1}\ovl y_{i-1}},\q\\
&&\q\q y_j\mapsto y_j,\q
\ovl y_j\mapsto\ovl y_j\,\,(j\ne i),
\q(i=2,\ld,n-1),\\
&&\gamma_0(v_2(y))
=\frac{y_0^2}{(y_1\ovl y_1)^2},\q
\gamma_1(v_2(y))=
\frac{(y_1\ovl y_1)^2}{y_0y_2\ovl y_2},\\
&&\gamma_i(v_2(y))=
\frac{(y_i\ovl y_i)^2}  
{y_{i-1}\ovl y_{i-1}y_{i+1}\ovl y_{i+1}}
\q(i=2,\ld, n-2),\q
\gamma_{n-1}(v_2(y))
=\frac{(y_{n-1}\ovl y_{n-1})^2}
{y_{n-2}\ovl y_{n-2}y_{n}^2},\\
&&\vep_0(v_2(y))=\frac{y_1^2}{y_0},\,
\vep_1(v_2(y))=\frac{y_2}{y_1}
\left(1+\frac{y_0}{y_1\ovl y_1}\right),\,
\vep_{n-1}(v_2(y))=
\frac{y_n^2}{y_{n-1}}\left(
1+\frac{y_{n-2}\ovl y_{n-2}}
{y_{n-1}\ovl y_{n-1}}\right),\\
&&\vep_i(v_2(y))=\frac{y_{i+1}}{y_i}
\left(1+\frac{y_{i-1}\ovl y_{i-1}}
{y_{i}\ovl y_i}\right)\q
(i=2\cd,n-1).
\end{eqnarray*}
The explicit forms of $\vep_i(v_1(x))$and 
$\gamma_i(v_1(x))$
$(1\leq i\leq n)$ are also induced from 
the ones
for $Y_{\bf i_1}(x)\cdot l^H\seteq 
Y_0(x_0)\cd Y_1(\ovl x_1)\cdot l^H$ and, we define
$\vep_0(v_1(x))\seteq \vep_0(v_2(y))$ and 
$\gamma_0(v_1(x))\seteq \gamma_0(v_2(y))$
$(v_2(y)\seteq \osigma(v_1(x)))$:
\begin{eqnarray*}
&&\vep_0(v_1(x))=
\frac{1}{x_0}\left(1+\frac{lx_1\ovl x_1}{x_0}
\right)^2,\q
\vep_i(v_1(x))=\frac{x_{i-1}}{x_i}
\left(1+\frac{x_{i+1}\ovl x_{i+1}}
{x_{i}\ovl x_i}\right)\q
(1\leq i\leq n-2),\\
&&\vep_{n-1}(v_1(x))=\frac{x_{n-2}}{x_{n-1}}
\left(1+\frac{x_n^2}{x_{n-1}\ovl x_{n-1}}
\right),\q
\vep_n(v_1(x))=\frac{x_{n-1}}{x_n}.
\end{eqnarray*}
\begin{eqnarray*}
&&\hspace{-23pt}\gamma_0(v_1(x))
=\frac{x_0^2}{(lx_1\ovl x_1)^2},\,\,
\gamma_1(v_1(x))=\frac{(lx_1\ovl x_1)^2}
{x_0x_2\ovl x_2},\,\,
\gamma_i(v_1(x))=
\frac{(x_i\ovl x_i)^2}{x_{i-1}\ovl x_{i-1}x_{i+1}\ovl x_{i+1}}\,\,
(2\leq i\leq n-2),\\
&&\hspace{-23pt}\gamma_{n-1}(v_1(x))=
\frac{(x_{n-1}\ovl x_{n-1})^2}
{x_{n-2}\ovl x_{n-2}x_n^2},\q
\gamma_n(v_1(x))=\frac{x_n^2}{x_{n-1}\ovl x_{n-1}}
\end{eqnarray*}
For $i=0$, we define 
$e^c_0(v_1(x))=
\osigma^{-1}\circ e^c_0\circ\osigma(v_1(x))
=\osigma^{-1}\circ e^c_0(v_2(y))$.
Then we get
\begin{eqnarray*}
e_0^c\cl&x_0\mapsto&x_0
\frac{(cx_0+lx_1\bar x_1)^2}
{c(x_0+lx_1\bar x_1)^2},\qq
x_i\mapsto{x_i}\frac{cx_0+lx_1\bar x_1}
{c(x_0+lx_1\bar x_1)}\,\,(1\leq i\leq n),\\
&\bar x_i\mapsto&{\bar x_i}
\frac{cx_0+lx_1\bar x_1}
{c(x_0+lx_1\bar x_1)}\,\,(1\leq i\leq n-1),\\
e_i^c\cl&x_i\mapsto&
x_i\frac{cx_i\ovl x_i+x_{i+1}\bar x_{i+1}}
{x_i\ovl x_i+x_{i+1}\ovl x_{i+1}},\q
\ovl x_i\mapsto
\ovl x_i\frac{c(x_i\ovl x_i+x_{i+1}\bar x_{i+1})}
{cx_i\ovl x_i+x_{i+1}\ovl x_{i+1}},\ \\
&&\q x_j\mapsto x_j,\q
\ovl x_j\mapsto\ovl x_j\,\,(j\ne i),\q
(1\leq i<n-1),\\
e_{n-1}^c\cl&x_{n-1}\mapsto&
x_{n-1}\frac{cx_{n-1}\ovl x_{n-1}+x_{n}^2}
{x_{n-1}\ovl x_{n-1}+x_{n}^2},\qq
\ovl x_{n-1}\mapsto
\ovl x_{n-1}
\frac{c(x_{n-1}\ovl x_{n-1}+x_{n}^2)}
{cx_{n-1}\ovl x_{n-1}+x_{n}^2}\q,\\
&&\q x_i\mapsto x_i,\q
\ovl x_i\mapsto\ovl x_i\,\,(i\ne n-1),\\
e^c_n\cl&x_n\mapsto& cx_n,
\q x_i\mapsto x_i,\q
\ovl x_i\mapsto\ovl x_i\,\,(i\ne n).
\end{eqnarray*}
In order to prove Theorem \ref{aff-geo}, it 
suffices to show the following:
\begin{eqnarray}
&&
e_i^c=\osigma^{-1}\circ 
e_{i}^c\circ \osigma,\,\,
\gamma_i=\gamma_{i}\circ \osigma,\,\,
\vep_i=\vep_{i}\circ \osigma
\,\,(i\ne0,n),
\label{0is}\\
&&e_0^{c_1}e_n^{c_2}=e_n^{c_2}e_0^{c_1},
\\
&&\gamma_n(e_0^c(v_1(x)))=\gamma_n(v_1(x)),\q
\gamma_0(e_n^c(v_1(x)))=\gamma_0(v_1(x)),
\\
&&\vep_0(e_0^c(v_1(x)))=c^{-1}\vep_0(v_1(x)),
\end{eqnarray}
which are immediate from the above formulae.
Let us show \eqref{0is}.
Set $v_2(y)\seteq \osigma(v_1(x))$ and $v_1(x')
\seteq \osigma^{-1}(e_i^c(v_2(y)))$
for $i=2,\ld,n-2$. Then
$x_j'=x_j$ and $\ovl x_j'=\ovl x_j$ for $j\ne i-1,i$, and we have 
\begin{eqnarray*}
&&a(v_2(y))=a(v_1(x)),\\
&&x_i'=\frac{1}{a(v_2(y))}\left(
\ovl y_{i+1}+\frac{cy_i\ovl y_i}{y_{i+1}}\right)
=x_i\frac{cx_i\ovl x_i+x_{i+1}\ovl x_{i+1}}
{x_i\ovl x_i+x_{i+1}\ovl x_{i+1}},\\
&&\ovl x_i'=\frac{1}{a(v_2(y))}
\frac{cy_i\ovl y_iy_{i+1}}
{cy_i\ovl y_i+y_{i+1}\ovl y_{i+1}}
=\ovl x_i\frac{c(x_i\ovl x_i
+x_{i+1}\ovl x_{i+1})}
{cx_i\ovl x_i+x_{i+1}\ovl x_{i+1}},\\
&&x_{i-1}'=\frac{1}{a(v_2(y))}\left(
\ovl y_{i}+\frac{cy_{i-1}\ovl y_{i-1}}
{y_{i}}\right)
=x_{i-1},\\
&&\ovl x_{i-1}'=\frac{1}{a(v_2(y))}
\frac{cy_{i-1}\ovl y_{i-1}y_{i}}
{cy_i\ovl y_i+y_{i-1}\ovl y_{i-1}}
=\ovl x_{i-1},
\end{eqnarray*}
where the formula
$y_i\ovl y_i=a(v_1(x))x_i\ovl x_i$
is useful to obtain these results.
Therefore we have $e_i^c=\osigma^{-1}\circ 
e_i^c\circ\osigma$ for $i=2,\ld,n-2$. 
Others are obtained similarly.
%%%%%%%%%%%%%%%%%%%%%%%%%%%%%%%%%%
\subsection{Ultra-discretization of $\cV(\ge)_l$}
\label{subsec-ud}

Let us investigate the 
ultra-discretization of $\cV(\ge)_l$.

By the explicit forms of the geometric crystal $\cV(\ge)_l$, 
if we assume that $l$ is a positive real number, 
it is clear that it has a natural 
positive structure 
$\theta_l\cl(\bbC^\times)^m\to\cV(\ge)_l$ 
($x\mapsto v(x)$)
where $m=\dim \cV(\ge)_l$.
Then we have the following theorem:
\begin{thm}
\label{ud-geo}
For $\ge=\TY(A,1,n), \TY(B,1,n),\TY(D,1,n), 
\TY(D,2,n+1), \TY(A,2,2n-1)$ and $\TY(A,2,2n)$,
suppose that $l>0$. Then the ultra-discretization 
${\mathcal UD}_{\theta_l}(\cV(\ge)_l)$ associated with the 
positive structure $\theta_l$ is isomorphic to the 
crystal $B_\ify(\ge^L)$ (\ref{limit}).
\end{thm}
{\sl Proof.}
First we consider the case $\ge\ne\TY(A,2,2n)$.
Let $\cV(\ge)$ be 
the affine geometric crystal in \cite{KNO}.
It follows from its explicit form (\cite{KNO}) 
that 
restricting $l=1$, we have the isomorphism
\begin{equation}
{\cV(\ge)_l}_{|l=1}\mapright{\sim}\cV(\ge).
\label{l1=1}
\end{equation}
The resulting crystal of the 
ultra-discretization ${\mathcal UD}_{\theta_l}$
does not depend on $l$. 
This and (\ref{l1=1}) imply 
\[
{\mathcal UD}_{\theta_l}(\cV(\ge)_l)
\cong{\mathcal UD}_{\theta}(\cV(\ge))
\]
as crystals.
We show in \cite{KNO} that 
${\mathcal UD}_{\theta}(\cV(\ge))\cong B_\ify(\ge^L)$.
Thus, we have $ {\mathcal UD}_{\theta_l}(\cV(\ge)_l)
\cong B_\ify(\ge^L)$.

% The proof of the last theorem in section 9 
% has been moved here
%%% See note X87-89 %%%
As for the case $\ge=\TY(A,2,2n)$, we consider 
as follows:
For the geometric crystal 
$\cV({\TY(A,2,2n)}^\dagger)_l$, we define 
$\ovl e_i\seteq e_{n-i}$, $\ovl\gamma_i\seteq
\gamma_{n-i}$ and $\ovl\vep_i\seteq \vep_{n-i}$
and set $x_i\mapsto y_{n-i}$ ($i=0,\cd,n$)
and $\ovl x_i\mapsto \ovl y_{n-i}$
($i=1\cd,n-1$). 
We denote by $\ovl \cV_l$ the $\TY(A,2,2n)$-geometric 
crystal thus obtained.
The explicit form of $\ovl\cV_l$
is {\it e.g.}, 
$\ovl e_0^c:y_0\mapsto cy_0$, 
$\ovl\vep_0(y)=\frac{y_1}{y_0}$ and 
$\ovl\gamma_0(y)=\frac{y_n^2}{y_1\ovl y_1}$, 
etc.
By using the birational map 
$\osigma$ for $\TY(A,2,2n)$ as in 
\ref{aen-1}, we obtain the 
isomorphism of $\TY(A,2,2n)$-geometric crystal 
\[
\osigma:
\cV(\TY(A,2,2n))_l\mapright{\sim}\ovl\cV_{l^2}.
\]
Here $\ovl\cV_l$ is isomorphic to 
$\cV({\TY(A,2,2n)}^\dagger)_l$ as an 
${\TY(A,2,2n)}^\dagger$-geometric crystal.
Due to the results in
\cite{KNO}
and the fact
$\cV({\TY(A,2,2n)}^\dagger)_l|_{l=1}=\cV({\TY(A,2,2n)}^\dagger)$,
we have
\[
 {\mathcal UD}_{\theta_{l^2}}
(\ovl\cV_{l^2})\cong
{\mathcal UD}_{\theta_{l^2}}(\cV({\TY(A,2,2n)}^\dagger)_{l^2})
={\mathcal UD}_{\theta}(\cV({\TY(A,2,2n)}^\dagger))\cong B_\ify(\TY(A,2,2n)).
\]
Let us denote by $\ovl B_l$ the perfect crystal 
$B_l(\TY(A,2,2n))$ considered as 
an ${\TY(A,2,2n)}^\dagger$-
perfect crystal.
Then the crystal $B_\ify(\TY(A,2,2n))$ 
can be regarded as  the limit 
of perfect ${\TY(A,2,2n)}^\dagger$-crystals 
$\{\ovl B_l\}$.
\qed

%%%%%%% Section 6 %%%%%%%
%%%%%%%%%%%%%%%%%%%%%%%%%%%%
\section{\bf Folding of 
$\TY(D,1,n)$-Geometric Crystal}
\label{fold}

\subsection{\bf $\TY(D,1,n)$-
geometric crystal $\cB_L(\TY(D,1,n))$}
%\hfill\break

We review the geometric crystal 
$\cB_L(\TY(D,1,n))$
($n\geq4$)
in \cite{KOTY}.

Taking $L\in\bbC^\times$, 
define the geometric crystal 
$\cB_L(\TY(D,1,n))\seteq(\cB_L(\TY(D,1,n)),
\{e_i\},\{\gamma_i\},\{\vep_i\})$ as follows:
\begin{eqnarray*}
&&\cB_L(\TY(D,1,n)):=\{l=(l_1,\ld,l_n,
\ovl l_{n-1},\ld,\ovl l_1)\in 
(\bbC^\times)^{2n-1}\,|\,l_1\cd 
l_n\ovl l_{n-1}\cd\ovl l_1=L\},\\
&&\hspace{-30pt}\vep_0(l)=l_1\left(\frac{l_2}{\ovl l_2}+1\right),\,
\vep_{n-1}(l)=l_n\ovl l_{n-1},\,
\vep_{n}(l)=\ovl l_{n-1},\,
\vep_i(l)=\ovl l_i\left(\frac{l_{i+1}}{\ovl l_{i+1}}+1\right),\\
&&\hspace{-30pt}
\gamma_0(l)=\frac{\ovl l_1\ovl l_2}{l_1l_2},\,
\gamma_{n-1}(l)=\frac{l_{n-1}}{l_n\ovl l_{n-1}},\,
\gamma_{n}(l)=\frac{l_{n-1}l_n}{\ovl l_{n-1}},\,
\gamma_i(l)=\frac{l_i\ovl l_{i+1}}{\ovl l_il_{i+1}},\\
&&\hspace{-30pt}
e_0^c(l)=(\frac{l_1}{\xi_2},\frac{\xi_2l_2}{c},\ld,\xi_2\ovl l_2,
\frac{c\ovl l_1}{\xi_2}),\,
e_{n-1}^c(l)=(\cdot\cdot,cl_{n-1},\frac{l_n}{c},\cdot\cdot),\,
e_{n}^c(l)=(\cdot\cdot,cl_{n},\frac{\ovl l_{n-1}}{c},\cdot\cdot),\\
&&\hspace{-30pt}
e_i^c(l)=(\cd,\frac{cl_i}{\xi_{i+1}},\frac{\xi_{i+1}l_{i+1}}{c},
\cd,\xi_{i+1}\ovl l_{i+1},\frac{\ovl l_i}{\xi_{i+1}},\ld),\q
(\xi_i\seteq\frac{c\ovl l_i+l_i}{\ovl l_i+l_i})
\,\,\,(i=1,\cd,n-2).
\end{eqnarray*}
%\newpage
This is isomorphic to the geometric crystal 
$\cV(\TY(D,1,n))_L$.
\begin{pro}
We have the following isomorphism of 
$\TY(D,1,n)$-geometric crystals:
\begin{equation}
\cB_L(\TY(D,1,n))\mapright{\sim}
\cV(\TY(D,1,n))_L.
\end{equation}
\end{pro}
{\sl Proof.}
Define $\Xi:\cB_L(\TY(D,1,n))\to
\cV(\TY(D,1,n))_L$ ($l\mapsto x$) to be 
\[
x_i\seteq \frac{1}{\ovl l_1\ovl l_2\cd\ovl l_i}
,\,\,
\ovl x_i\seteq\frac{l_1\cd l_i}{L}\,\,
(i=1\,\ld,n-2),\,\,
x_{n-1}\seteq \frac{1}
{\ovl l_1\ovl l_2\cd\ovl l_{n-1}l_n},\,\,
x_{n}\seteq \frac{1}
{\ovl l_1\ovl l_2\cd\ovl l_{n-1}}.
\]
The inverse $\Xi^{-1}$ is given by
\[
l_1=L\ovl x_1,\,\,\ovl l_1=\frac{1}{x_1},\,\,
l_i=\frac{\ovl x_i}{\ovl x_{i-1}},\,\,
\ovl l_i=\frac{x_{i-1}}{x_i}
\,\,
(i=2,\ld,n-2),\,\,
l_{n-1}=\frac{x_{n-1}}{\ovl x_{n-2}},\,\,
l_n=\frac{x_n}{x_{n-1}},\,\,
\ovl l_{n-1}=\frac{x_{n-2}}{x_n}.
\]
Then it is easy to see 
that it commutes with the action
of $e_i$ and preserves $\gamma_i$ and $\vep_i$.
For example, for $x=\Xi(l)$, we have 
\[
 \xi_2=\frac{c\ovl l_2+l_2}{\ovl l_2+l_2}
=\frac{c\frac{x_1}{x_2}+\frac{\ovl x_2}
{\ovl x_1}}{\frac{x_1}{x_2}
+\frac{\ovl x_2}{\ovl x_1}}
=\frac{cx_1\ovl x_1+x_2\ovl x_2}
{x_1\ovl x_1+x_2\ovl x_2}=:\Gamma_1
\]
Then we obtain
\begin{eqnarray*}
&&\Xi(e_0^c(l))=\Xi(\frac{l_1}{\xi_2},\frac{\xi_2l_2}{c},\ld,\xi_2\ovl l_2,
\frac{c\ovl l_1}{\xi_2})
=(\frac{cx_1}{\Gamma_1},\frac{x_2}{c},\ld,
\frac{\ovl x_2}{c},\Gamma_1\ovl x_1)
=e_0^c\Xi(l),\\
&& \vep_0(\Xi(l))=\vep_0(x)=
L\frac{x_1\ovl x_1+x_2\ovl x_2}{x_1}
=L\frac{\frac{l_1}{L\ovl l_1}+\frac{l_1l_2}
{L\ovl l_1\ovl l_2}}{\frac{1}{\ovl l_1}}
=l_1\left(\frac{l_2}{\ovl l_2}+1\right)=
\vep_0(l).
\end{eqnarray*}
Other cases are shown similarly.
\qed
%%%%%%%%%%%%%%%%%%%%%%%%%%%%%%%%%%%%%%%
\subsection{\bf Folding}

We introduce 
certain involutions
on $\cB_L(\TY(D,1,n))$
corresponding to 
the Dynkin diagram automorphisms
for $\TY(D,1,n)$.
Let us consider the following 
Dynkin diagram automorphisms
for $\TY(D,1,N)$:
\begin{eqnarray*}
N=n+1&&\sigma_0^{(n)}:
\al_{0}\leftrightarrow \al_1,\q
\al_i\mapsto \al_i\q(i\ne 0,1),\\
N=n+1&&\sigma_1^{(n)}:
\al_{n}\leftrightarrow \al_{n+1},\q
\al_i\mapsto \al_i\q(i\ne n,n+1),\\
N=2n&&\sigma_2^{(n)}:
\al_{i}\leftrightarrow \al_{2n-i}\q
(i=0,1,\ld, 2n).
\end{eqnarray*}
For each involution $\sigma^{(n)}_j$ 
on $\TY(D,1,N)$($j=0,1,2$),
we define the involution $\Sigma^{(n)}_j$ on
$\cB_L(\TY(D,1,N))$
($j=0,1,2$) as a unique solution $l'$
of the equations for a given $l$:
\begin{equation}
 \gamma_{\sigma^{(n)}_j(i)}(l')
=\gamma_i(l), \qq
 \vep_{\sigma^{(n)}_j(i)}(l')
=\vep_i(l)\q(i=0,1,\ld,N).
\label{Sigma-sigma}
\end{equation}
Then we get: 
($l'\seteq \Sigma_j^{(n)}(l)$)
\begin{eqnarray*}
\Sigma_0^{(n)}&\cl&\cB_L(\TY(D,1,n+1))
\longrightarrow \cB_L(\TY(D,1,n+1)),\\
&& l'_1= \ovl l_1,\q \ovl l'_1=l_1,\q
l'_i=l_i,\q \ovl l'_i=\ovl l_i\,\,
(i\ne1).\\
\Sigma_1^{(n)}&\cl&\cB_L(\TY(D,1,n+1))
\longrightarrow \cB_L(\TY(D,1,n+1)),\\
&& 
{l'_n=l_nl_{n+1},\q l'_{n+1}=\frac{1}{l_{n+1}},\q
\ovl l'_n=l_{n+1}\ovl l_n}, \q
l'_i=l_i,\q 
\ovl l'_i=\ovl l_i \q(i\ne,n,n+1).\\
\Sigma_2^{(n)}&\cl&\cB_L(\TY(D,1,2n))
\longrightarrow \cB_L(\TY(D,1,2n))\\
&&l'_1=\frac{l_{2n-1}\ovl l_{2n-1}}{l_{2n-1}+\ovl l_{2n-1}},\q
\ovl l'_1=\frac{l_{2n}l_{2n-1}\ovl l_{2n-1}}
{l_{2n-1}+\ovl l_{2n-1}},\q
l'_{2n}=\frac{\ovl l_1}{l_1},\\
&&l'_{2n-1}=\frac{l_1\ovl l_2}{l_2}
\left(\frac{l_2}{\ovl l_2}+1\right),\q
\ovl l'_{2n-1}=l_1\left(\frac{l_2}{\ovl l_2}+1\right),\\
&&
l'_{2n-i}=\frac{l_i\ovl l_i}{l_{i+1}}
\left(\frac{l_{i+1}+\ovl l_{i+1}}{l_i+\ovl l_i}\right),\q
\ovl l'_{2n-i}=\frac{l_i\ovl l_i}{\ovl l_{i+1}}
\left(\frac{l_{i+1}+\ovl l_{i+1}}{l_i+\ovl l_i}\right)\q
(2\leq i\le 2n-2).
\end{eqnarray*}
\begin{lem}
\label{a3}
Let $(X,\{e_i\}_{i=1,2,3},\{\gamma_i\}_{i=1,2,3},
\{\vep_i\}_{i=1,2,3})$ be a $A_3$-geometric
crystal ($a_{ij}=-1$ if $|i-j|=1$.) and 
set $E_1^c\seteq e_1^c\circ e_3^c
(=e_3^c\circ e_1^c)$. Then we have 
\begin{eqnarray*}
&&E_1^{c}e_2^{c^2d}E_1^{cd}e_2^{d}
=e_2^{d}E_1^{cd}e_2^{c^2d}E_1^{c} 
\q(c,d\in\bbC^\times),\q
\gamma_2(E_1^c(x))=c^{-2}\gamma_2(x).
\end{eqnarray*}
\end{lem}
{\sl Proof.} The second equation is easily 
obtained:
\[
\gamma_2(E_1^c(x))=\gamma_2(e_1^c\circ e_3^c(x))
=c^{-1}\gamma_2(e_3^c(x))=c^{-2}\gamma_2(x)\q
(x\in X)
\]
The first one is derived as follows:
\begin{eqnarray*}
&&E_1^{c}e_2^{c^2d}E_1^{cd}e_2^{d}=
(e_1^{c}e_3^c)e_2^{c^2d}(e_3^{cd}e_1^{cd})e_2^{d}=
e_1^{c}e_2^{cd}e_3^{c^2d}e_2^{c}e_1^{cd}e_2^{d}
=e_1^{c}e_2^{cd}e_3^{c^2d}e_1^{d}e_2^{cd}e_1^{c}\\
&&=e_1^{c}e_2^{cd}e_1^{d}e_3^{c^2d}e_2^{cd}e_1^{c}
=e_2^{d}e_1^{cd}e_2^{c}e_3^{c^2d}e_2^{cd}e_1^{c}
=e_2^{d}(e_1^{cd}e_3^{cd})e_2^{c^2d}
(e_3^{c}e_1^{c})
=e_2^{d}E_1^{cd}e_2^{c^2d}E_1^{c}.\qq\qq\qed
\end{eqnarray*}

\subsection{Fixed-point variety-$\TY(B,1,n)$}
\label{subsec-bn}
As for the involution $\Sigma_1^{(n)}$ 
we consider the fixed-point variety $X_1^{(n,L)}
\subset \cB_L(\TY(D,1,n+1))$:
\[
 X_1^{(n,L)}\seteq
\{l=(l_1,\ld,l_{n+1},\ovl l_n,\ld,\ovl l_1)|
\Sigma_1^{(n)}(l)=l
(\Leftrightarrow l_{n+1}=1)\}.
\]
\begin{pro}
\label{fix-b}
Let 
$\cB_L(\TY(D,1,n+1))=(\cB_L,\{e_i\},\{\gamma_i\},
\{\vep_i\})$ be the 
$\TY(D,1,n+1)$-geometric crystal as above.
Then we have a 
$\TY(B,1,n)$-geometric crystal
structure on $X_1^{(n,L)}$ 
as follows:
\[
( e^{(\TY(B,1,n))}_i)^c\seteq
\begin{cases}
e_i^c&\text{ for }i\ne n,\\
e_n^c\circ e_{n+1}^c&\text{ for }i=n,
\end{cases},\q
\gamma^{(\TY(B,1,n))}_i\seteq
\gamma_i,\qq
\vep^{(\TY(B,1,n))}_i\seteq
\vep_i\q(i=0,\ld,n).
\]
\end{pro}
{\sl Proof.}
It suffices to check the conditions 
for a geometric crystal. But, most ones are 
trivial except the cases 
related to $i=n$. Namely,
we have to see the Verma relation between 
$e_{n-1}$ and $e_n$, 
$\gamma_i^{(\TY(B,1,n))}((e_n^{\TY(B,1,n)})^c(x))$,
$\gamma_n^{(\TY(B,1,n))}((e_i^{\TY(B,1,n)})^c(x))$
and 
$\vep_n^{(\TY(B,1,n))}((e_n^{\TY(B,1,n)})^c(x))$, 
which are immediate from Lemma \ref{a3}.\qed

Moreover,  the geometric crystal $X_1^{(n,L)}$ 
induces another $\TY(B,1,n)$-
geometric crystal $\cB_L(\TY(B,1,n))$:
\begin{eqnarray*}
&&\cB_L(\TY(B,1,n))\seteq
\{m=(m_1,\ld,m_n,\ovl m_n,\ld,\ovl m_1)
(\bbC^\times)^{2n}\,|\,
m_1\cd 
m_n\ovl m_{n}\cd\ovl m_1=L\},\\
&&\hspace{-30pt}\vep_0(m)
=m_1\left(\frac{m_2}{\ovl m_2}+1\right),\,
\vep_{n}(m)=\ovl m_{n},\,
\vep_i(m)=\ovl m_i\left(\frac{m_{i+1}}
{\ovl m_{i+1}}+1\right)
\,(i=1,\ld,n-1),\\
&&\hspace{-30pt}
\gamma_0(m)=\frac{\ovl m_1\ovl m_2}{m_1m_2},\,
\gamma_{n}(m)=\frac{m_n}{\ovl m_{n}},\,
\gamma_i(m)
=\frac{m_i\ovl m_{i+1}}{\ovl m_im_{i+1}}
\,(i=1,\ld,n-1),\\
&&\hspace{-30pt}
e_0^c(m)=(\frac{m_1}{\xi_2},
\frac{\xi_2m_2}{c},\ld,\xi_2\ovl m_2,
\frac{c\ovl m_1}{\xi_2}),\,
e_{n}^c(m)=(\cdot\cdot,cm_{n},
\frac{\ovl m_{n}}{c},\cdot\cdot),\\
&&\hspace{-30pt}
e_i^c(m)=(\cd,\frac{cm_i}{\xi_{i+1}},\frac{\xi_{i+1}m_{i+1}}{c},
\cd,\xi_{i+1}\ovl m_{i+1},\frac{\ovl m_i}{\xi_{i+1}},\ld)\q(i=1,\ld,n-1)\q
(\xi_i\seteq\frac{c\ovl m_i+m_i}{\ovl m_i+m_i}).
\end{eqnarray*}
Let $\eta\cl\cB_L(\TY(B,1,n))
\longrightarrow X_1^{(n,L)}$ ($m\mapsto l$)
be the morphism defined by
\[
l_i=m_i,\q \ovl l_i=\ovl m_i \q
(i=1,\cd,n),
\]
where $(l_1,\ld,l_n,1,\ovl l_n,\ld, \ovl l_1)
\in X_1^{(n,L)}$ and 
$m=(m_1,\ld,m_n,\ovl m_n,\ovl m_1)\in 
\cB_L(\TY(B,1,n))$.
It is trivial that 
$\eta$ commutes with the actions $e_i$ and 
preserves $\gamma_i$ and $\vep_i$. 
And then $\eta$ is an isomorphism of 
$\TY(B,1,n)$-geometric crystals.
\begin{pro}
We have the following isomorphisms of 
$\TY(B,1,n)$-geometric crystals:
\begin{equation}
\cV(\TY(B,1,n))_L\mapleft{\sim}
\cB_L(\TY(B,1,n))\mapright{\sim}
X_1^{(n,L)}.
\label{3iso-bn}
\end{equation}
\end{pro}
{\sl Proof.}
The second isomorphism in (\ref{3iso-bn})
is given by $\eta$. 
Then let us see the first one.
Define $\Xi\cl\cB_L(\TY(B,1,n))
\to\cV(\TY(B,1,n))_L$ ($m\mapsto x$) to be 
\[
x_i\seteq \frac{1}{\ovl m_1
\ovl m_2\cd\ovl m_i} \,\,(i=1,\ld, n),\q
\ovl x_i=\frac{m_1\cd m_i}{L}\,\,
(i=1,\ld, n-1),
\]
and the inverse is 
\[
m_1=L\ovl x_1,\q
m_i=\frac{\ovl x_i}{\ovl x_{i-1}}\,\,
(2\le i\leq n-1),\q
m_n=\frac{x_n}{\ovl x_{n-1}},\q
\ovl m_1=\frac{1}{x_1},\q
\ovl m_i=\frac{x_{i-1}}{x_i}\,\,
(2\leq i\leq n).
\]
Then, by a direct calculation, we can check 
that $\Xi$ commutes with any $e_i^c$ and preserves
$\gamma_i$ and $\vep_i$.
Thus, $\Xi$ is an 
isomorphism of $\TY(B,1,n)$-geometric 
crystals.\qed
%%%%%%%%%%%%%%%%%%%%%%%%%%%%%%%%
\subsection{Fixed-point variety-$\TY(D,2,n+1)$}
\label{subsec-d2n}
As for the involution 
$\Sigma_3^{(n)}
\seteq\Sigma_0^{(n+1)}\circ\Sigma_1^{(n+1)}
(=\Sigma_1^{(n+1)}\circ\Sigma_0^{(n+1)})$ 
on $\cB_L(\TY(D,1,n+2))$, 
we consider the fixed-point variety $X_3^{(n,L)}
\subset \cB_L(\TY(D,1,n+2))$:
\[
 X_3^{(n,L)}\seteq
\{l=(l_1,\ld,l_{n+2},\ovl l_{n+1},\ld,\ovl l_1)|
\Sigma_3^{(n,L)}(l)=l
(\Leftrightarrow l_{n+2}=1,\,\,
l_1=\ovl l_1)\}.
\]
\begin{pro}
Let 
$\cB_L(\TY(D,1,n+2))=(\cB_L,\{e_i\},\{\gamma_i\},
\{\vep_i\})$ be the 
$\TY(D,1,n+2)$-geometric crystal as above.
Then we have a 
$\TY(D,2,n+1)$-geometric crystal
structure on $X_3^{(n,L)}$ 
as follows:
\[
( e^{(\TY(D,2,n+1))}_i)^c\seteq
\begin{cases}
e_{i+1}^c&\text{ for }i\ne 0,n,\\
e_{0}^c\circ e_{1}^c&\text{ for }i=0,\\
e_{n+1}^c\circ e_{n+2}^c&\text{ for }i=n,
\end{cases},\q
\gamma^{(\TY(D,2,n+1))}_i\seteq
\gamma_{i+1},\qq
\vep^{(\TY(D,2,n+1))}_i\seteq
\vep_{i+1}\q(i=0,\ld,n).
\]
\end{pro}
{\sl Proof.} We can prove it by the 
similar argument for
 $X^{(n,L)}_1$ using Lemma \ref{a3}.\qed

The geometric crystal $X_3^{(n,L)}$ induces 
another $\TY(D,2,n+1)$-
geometric crystal $\cB_L(\TY(D,2,n+1))$:
\begin{eqnarray*}
&&\cB_L(\TY(D,2,n+1))\seteq
\{m=(m_0,m_1,\ld,m_n,\ovl m_n,\ld,\ovl m_1)
(\bbC^\times)^{2n}\,|\,
m_0^2m_1\cd 
m_n\ovl m_{n}\cd\ovl m_1=L\},\\
&&\hspace{-30pt}\vep_0(m)
=m_0\left(\frac{m_1}{\ovl m_1}+1\right),\,
\vep_{n}(m)=\ovl m_{n},\,
\vep_i(m)=\ovl m_i\left(\frac{m_{i+1}}
{\ovl m_{i+1}}+1\right)\,\,(i=1,\ld,n-1),\\
&&\hspace{-30pt}
\gamma_0(m)=\frac{\ovl m_1}{m_1},\,
\gamma_{n-1}(m)=\frac{m_{n-1}}{m_n\ovl m_{n-1}},\,
\gamma_{n}(m)=\frac{m_n}{\ovl m_{n}},\,
\gamma_i(m)
=\frac{m_i\ovl m_{i+1}}{\ovl m_im_{i+1}}\,\,
(i=1,\ld,n-2),\\
&&\hspace{-30pt}
e_0^c(m)=(\frac{c m_0}{\xi_1},
\frac{\xi_1m_1}{c^2},\ld,\xi_1\ovl m_1)
\,\,(\xi_1(m)\seteq
\frac{\ovl m_1+m_1}{c^2\ovl m_1+m_1}),\,\,
e_{n}^c(m)=(\cdot\cdot,cm_{n},
\frac{\ovl m_{n}}{c},\cdot\cdot),\\
&&\hspace{-30pt}
e_i^c(m)=(\cd,\frac{cm_i}{\xi_{i+1}},\frac{\xi_{i+1}m_{i+1}}{c},
\cd,\xi_{i+1}\ovl m_{i+1},\frac{\ovl m_i}{\xi_{i+1}},\ld)\q(i=1,\ld,n-1)\q
(\xi_i\seteq\frac{c\ovl m_i+m_i}
{\ovl m_i+m_i}).
\end{eqnarray*}
Let $\eta\cl\cB_L(\TY(D,2,n+1))
\longrightarrow X_3^{(n,L)}$
($m\mapsto l$)
be the morphism defined by
\[
l_{i+1}=m_i\q (i=0,1,\cd,n),\q
\ovl l_{i+1}=\ovl m_i\q(i=1,\cd,n),
\]
where $(l_1,\ld,l_{n+1},1,\ovl l_{n+1},
\cd, \ovl l_2,l_1)
\in X_3^{(n,L)}$ and 
$m=(m_0,m_1,\ld,m_n,\ovl m_n,\ovl m_1)\in 
\cB_L(\TY(D,2,n+1))$.
Then, $\eta$ is an isomorphism of 
$\TY(D,2,n+1)$-
geometric crystal.
\begin{pro}
We have the following isomorphisms of 
$\TY(D,2,n+1)$-geometric crystals:
\begin{equation}
\cV(\TY(D,2,n+1))_L\mapleft{\sim}
\cB_{L^2}(\TY(D,2,n+1))\mapright{\sim}
X_3^{(n,L^2)}.
\label{3iso-d2n}
\end{equation}
\end{pro}
{\sl Proof.}
The second isomorphism in (\ref{3iso-d2n})
is given by $\eta$. 
Then let us see the first one.
Define $\Xi\cl\cB_{L^2}(\TY(D,2,n+1))
\to\cV(\TY(D,2,n+1))_L$ ($m\mapsto x$) to be 
\[
x_0\seteq \frac{1}{m_0},\q
x_i\seteq \frac{1}{m_0^2\ovl m_1
\ovl m_2\cd\ovl m_i} \,\,(i=1,\ld, n),\q
\ovl x_i=\frac{m_1\cd m_i}{L^2}\,\,
(i=1,\ld, n-1).
\]
and the inverse is 
\[
m_0=\frac{1}{x_0},\,\, m_1=L^2\ovl x_1,\,\,
m_i=\frac{\ovl x_i}{\ovl x_{i-1}}\,\,
(2\le i\leq n-1),\,\,
m_n=\frac{x_n}{\ovl x_{n-1}},\,
\ovl m_1=\frac{x_0^2}{x_1},\,\,
\ovl m_i=\frac{x_{i-1}}{x_i}\,\,
(2\leq i\leq n).
\]
Then, calculating directly, we can check 
that $\Xi$ is an isomorphism of 
geometric crystals.\qed

\subsection{Fixed-point variety-$\TY(A,2,2n-1)$}
\label{subsec-a2o}
As for the involution $\Sigma_2^{(n)}$ 
we consider the fixed-point variety $X_2^{(n,L)}
\subset \cB_L(\TY(D,1,2n))$:
\[
 X_2^{(n,L)}\seteq
\{l=(l_1,\ld,l_{2n},\ovl l_{2n-1},
\cd,\ovl l_1)|
\Sigma_2^{(n)}(l)=l
\},
\]
where the condition $\Sigma_2^{(n)}(l)=l$ is 
equivalent to 
\begin{equation}
\begin{array}{ll}
& l_1=\frac{l_{2n-1}\ovl l_{2n-1}}
{l_{2n-1}+\ovl l_{2n-1}},\q
\ovl  l_1=\frac{l_{2n}l_{2n-1}\ovl l_{2n-1}}
{l_{2n-1}+\ovl l_{2n-1}},\q
l_{2n-1}=\frac{l_1\ovl l_2}{l_2}
\left(\frac{l_2}{\ovl l_2}+1\right),\q
\ovl l_{2n-1}=l_1
\left(\frac{l_2}{\ovl l_2}+1\right),\\
&l_{2n-i}=\frac{l_i\ovl l_i}{l_{i+1}}
\frac{l_{i+1}+\ovl l_{i+1}}{l_i+\ovl l_i},\q
l_{2n-i}=\frac{l_i\ovl l_i}{\ovl l_{i+1}}
\frac{l_{i+1}+\ovl l_{i+1}}{l_i+\ovl l_i}\q
(i=2,\ld, 2n-2), \q
l_{2n}=\frac{\ovl l_1}{l_1}.
\end{array}
\label{aon-eq}
\end{equation}
Note that all the equations in (\ref{aon-eq})
are not necessarily independent.
Indeed, we only need $2n$-equations:
$l_1=\cd,l_n=\cd$ and
$\ovl l_1=\cd,\ovl l_n=\cd$ in (\ref{aon-eq}).
\begin{pro}
Let 
$\cB_L(\TY(D,1,2n))=(\cB_L,\{e_i\},\{\gamma_i\},
\{\vep_i\})$ be the 
$\TY(D,1,2n)$-geometric crystal as above.
Then we have a 
$\TY(A,2,2n-1)$-geometric crystal
structure on $X_2^{(n,L)}$ 
as follows:
\[
( e^{(\TY(A,2,2n-1))}_i)^c\seteq
\begin{cases}
e_i^c\circ e_{2n-i}^c
&\text{ for }0\leq i\leq n-1,\\
e_n^c&\text{ for }i=n,
\end{cases}\q
\begin{array}{l}
\gamma^{(\TY(A,2,2n-1))}_i\seteq
\gamma_i(=\gamma_{2n-i}),\\
\vep^{(\TY(A,2,2n-1))}_i\seteq
\vep_i(=\vep_{2n-i})\q(i=0,\ld,n).
\end{array}
\]
\end{pro}
The proof is similar to the previous cases.

The geometric crystal $X_2^{(n,L)}$ induces 
another $\TY(A,2,2n-1)$-
geometric crystal $\cB_L(\TY(A,2,2n-1))$:
\[
 \cB_L(\TY(A,2,2n-1))\seteq
\{m=(m_1,\ld,m_n,\ovl m_n,\ld,\ovl m_1)\in
(\bbC^\times)^{2n}\,|\,
m_n\ovl m_n(m_1\cd 
m_{n-1}\ovl m_{n-1}\cd\ovl m_1)^2=L^2\}.
\]
For $m\in \cB_L(\TY(A,2,2n-1))$ set 
\[
M(m)\seteq \frac{L}{m_1\cd m_{n-1}
\ovl m_n\cd \ovl m_1},\qq
\ovl M(m)\seteq \frac{L}{m_1\cd m_{n}
\ovl m_{n-1}\cd \ovl m_1}(=\frac{1}{M(m)}).
\]
\begin{eqnarray*}
&&\hspace{-30pt}\vep_0(m)
=m_1\left(\frac{m_2}{\ovl m_2}+1\right),\,
\vep_{n-1}(m)=\ovl m_{n-1}(M(m)+1),\,
\vep_{n}(m)=\ovl m_{n},\,
\vep_i(m)=\ovl m_i\left(\frac{m_{i+1}}
{\ovl m_{i+1}}+1\right),\\
&&\hspace{-30pt}
\gamma_0(m)=\frac{\ovl m_1\ovl m_2}{m_1m_2},
\,\,
\gamma_{n-1}(m)=\frac{m_{n-1}}{\ovl m_{n-1}}
\ovl M(m),\,\,
\gamma_{n}(m)=\frac{m_n}{\ovl m_{n}},\,\,
\gamma_i(m)
=\frac{m_i\ovl m_{i+1}}{\ovl m_im_{i+1}}
\,\,\,(1\leq i\leq n-2),\\
&&\hspace{-30pt}
e_0^c(m)=(\frac{m_1}{\xi_1},
\frac{\xi_1m_2}{c},\ld,\xi_1\ovl m_2,
\frac{c\ovl m_1}{\xi_1}),\,
e_{n}^c(m)=(\cdot\cdot,cm_{n},
\frac{\ovl m_{n}}{c},\cdot\cdot),\\
&&\hspace{-30pt}
e_{n-1}^c(m)=(\cd,m_{n-1}\frac{c}{\xi_{n-1}},
m_n\frac{\xi_{n-1}^2}{c^2},\ovl m_n\xi_{n-1}^2,
\frac{\ovl m_{n-1}}{\xi_{n-1}},\ld),\\
&&\hspace{-30pt}
e_i^c(m)=(\cd,\frac{cm_i}{\xi_{i}},
\frac{\xi_{i}m_{i+1}}{c},
\cd,\xi_{i}\ovl m_{i+1},\frac{\ovl m_i}
{\xi_{i}},\ld)\q(i=1,\ld,n-1),
\end{eqnarray*}
where
\[
\xi_i\seteq
\begin{cases}
\frac{c\ovl m_{i+1}+m_{i+1}}
{\ovl m_{i+1}+m_{i+1}}
&\text{ for }i\ne n-1\\
\frac{c+M(m)}{1+M(m)}&\text{ for }i=n-1.
\end{cases}
\]
Let $\eta\cl\cB_L(\TY(A,2,2n-1))
\longrightarrow X_2^{(n,L^2)}$
($m\mapsto l$)
be the morphism defined by
\[
l_i=m_i,\q \ovl l_i=\ovl m_i \q
(i=1,\ld, n-1),\q
l_n=\frac{m_n}{1+M(m)},\q
\ovl l_n=\frac{\ovl m_n}{1+\ovl M(m)},\q
\]
where $l_i,\ovl l_i$ ($i=n+1,\ld,2n$) 
are uniquely determined by (\ref{aon-eq}) and 
then, $\eta$ is an isomorphism of 
$\TY(A,2,2n-1)$-
geometric crystal.
\begin{pro}
We have the following isomorphisms of 
$\TY(A,2,2n-1)$-geometric crystals:
\begin{equation}
\cV(\TY(A,2,2n-1))_L\mapleft{\sim}
\cB_{L}(\TY(A,2,2n-1))\mapright{\sim}
X_2^{(n,L^2)}.
\label{3iso-aon}
\end{equation}
\end{pro}
{\sl Proof.}
The second isomorphism in (\ref{3iso-aon})
is given by $\eta$. 
Then let us see the first one.
Define $\Xi\cl\cB_{L}(\TY(A,2,2n-1))
\to\cV(\TY(A,2,2n-1))_L$ ($m\mapsto x$) to be 
\[
x_i\seteq \frac{1}{\ovl m_1
\ovl m_2\cd\ovl m_i},\q
\ovl x_i=\frac{m_1\cd m_i}{L^2}\,\,
(i=1,\ld, n-1),\q
x_n=\frac{1}
{\ovl m_n(\ovl m_1\cd \ovl m_{n-1})^2},
\]
and then the inverse is 
\[
m_1=L^2\ovl x_1,\q
\ovl m_1=\frac{1}{x_1},\q
m_i=\frac{\ovl x_i}{\ovl x_{i-1}},\q
\ovl m_i=\frac{x_{i-1}}{x_i}\,\,
(2\leq i\leq n-1),\q
m_n=\frac{x_n}{L^2\ovl x_{n-1}^2},\q
\ovl m_n=\frac{x_{n-1}^2}{x_n}.
\]
Then, by direct calculations, we can check 
that $\Xi$ is an isomorphism of 
geometric crystals.\qed
%%%%%%%%%%%%%%%%%%%%%%%%%%%%%%%%%%%%%%%%%%
\subsection{Fixed-point variety-$\TY(A,2,2n)$}
\label{subsec-a2e}
As for the involution 
$\Sigma_4^{(n)}\seteq
\Sigma_2^{(n+1)}\Sigma_1^{(2n+1)}
\Sigma_0^{(2n+1)}$ 
we consider a fixed-point variety $X_4^{(n,L)}
\subset \cB_L(\TY(D,1,2n+2))$:
\[
 X_4^{(n,L)}\seteq
\{l=(l_1,\ld,l_{2n+2},\ovl l_{2n+1},
\cd,\ovl l_1)|
\Sigma_4^{(n)}(l)=l
\},
\]
where the condition $\Sigma_4^{(n)}(l)=l$ is 
equivalent to 
\begin{equation}
\begin{array}{ll}
& l_1=\ovl l_1
=\frac{l_{2n+1}\ovl l_{2n+1}}
{l_{2n+1}+\ovl l_{2n+1}},\q
l_{2n+1}=\frac{l_1\ovl l_2}{l_2}
\left(\frac{l_2}{\ovl l_2}+1\right),\q
\ovl l_{2n+1}=l_1
\left(\frac{l_2}{\ovl l_2}+1\right),\q
l_{2n+2}=1,\\
& l_{2n+2-i}=\frac{l_i\ovl l_i}{l_{i+1}}
\frac{l_{i+1}+\ovl l_{i+1}}{l_i+\ovl l_i},\q
l_{2n+2-i}=\frac{l_i\ovl l_i}{\ovl l_{i+1}}
\frac{l_{i+1}+\ovl l_{i+1}}{l_i+\ovl l_i}\q
(i=2,\ld, 2n).
\end{array}
\label{aen-eq}
\end{equation}
Note that all the equations in (\ref{aen-eq})
are not necessarily independent.
Indeed, we only need $2n$-equations:
$l_1=\cd,l_n=\cd$ and
$\ovl l_1=\cd,\ovl l_n=\cd$ in (\ref{aen-eq}).
\begin{pro}
Let 
$\cB_L(\TY(D,1,2n+2))=(\cB_L,\{e_i\},\{\gamma_i\},
\{\vep_i\})$ be the 
$\TY(D,1,2n+2)$-geometric crystal as above.
Then we have a 
$\TY(A,2,2n)$-geometric crystal
structure on $X_4^{(n,L)}$ 
as follows:
\[
\begin{array}{l}
%\hspace{-30pt}
(e^{(\TY(A,2,2n))}_i)^c\seteq
\begin{cases}
e_0^c\circ e_1^c\circ e_{2n+1}^c
\circ e_{2n+2}^c
&\text{ for }i=0,\\
e_{i+1}^c\circ e_{2n-i+1}^c
&\text{ for }1\leq i\leq n-1,\\
e_{n+1}^c&\text{ for }i=n,
\end{cases}\\
\gamma^{(\TY(A,2,2n))}_i\seteq
\begin{cases}
\gamma_0(=\gamma_1=\gamma_{2n+1}=\gamma_{2n+2})
&\text{ for }i=0,\\
\gamma_{i+1}(=\gamma_{2n-i+1})
&\text{ for }i\ne 0,
\end{cases}\\
\vep^{(\TY(A,2,2n))}_i\seteq
\begin{cases}
\vep_0(=\vep_1=\vep_{2n+1}=\vep_{2n+2})
&\text{ for } i=0,\\
\vep_{i+1}(=\vep_{2n-i+1})
&\text{ for }i\ne0.
\end{cases}
\end{array}
\]
\end{pro}
We can prove it by a similar argument to
the ones for the previous cases.

The geometric crystal $X_4^{(n,L)}$ induces 
another $\TY(A,2,2n)$-
geometric crystal $\cB_L(\TY(A,2,2n))$:
\[
 \cB_L(\TY(A,2,2n))\seteq
\{m=(m_0,m_1,\ld,m_n,\ovl m_n,\ld,\ovl m_1)
\in(\bbC^\times)^{2n}\,|\,
m_n\ovl m_n(m_0^2m_1\cd 
m_{n-1}\ovl m_{n-1}\cd\ovl m_1)^2=L^2\}.
\]
For $m\in \cB_L(\TY(A,2,2n))$ set 
\[
M(m)\seteq \frac{L}
{m_0^2m_1\cd m_{n-1}
\ovl m_n\cd \ovl m_1},\qq
\ovl M(m)\seteq \frac{L}
{m_0^2m_1\cd m_{n}
\ovl m_{n-1}\cd \ovl m_1}(=\frac{1}{M(m)}).
\]
\begin{eqnarray*}
&&\hspace{-30pt}\vep_0(m)
=m_0\left(\frac{m_1}{\ovl m_1}+1\right),\,
\vep_{n-1}(m)=\ovl m_{n-1}(M(m)+1),\,
\vep_{n}(m)=\ovl m_{n},\,
\vep_i(m)=\ovl m_i\left(\frac{m_{i+1}}
{\ovl m_{i+1}}+1\right)
,\\
&&\hspace{-30pt}
\gamma_0(m)=\frac{\ovl m_1}{m_1},
\,\,
\gamma_{n-1}(m)=\frac{m_{n-1}}{\ovl m_{n-1}}
\ovl M(m),\,\,
\gamma_{n}(m)=\frac{m_n}{\ovl m_{n}},\,\,
\gamma_i(m)
=\frac{m_i\ovl m_{i+1}}{\ovl m_im_{i+1}}
\,\,\,(1\leq i\leq n-2),\\
&&\hspace{-30pt}
e_0^c(m)=(\frac{m_0}{\xi_0},
\frac{\xi_0m_1}{c},m_2,\ld,\ovl m_2,
\frac{c\ovl m_1}{\xi_0}),\,
e_{n}^c(m)=(\cdot\cdot,cm_{n},
\frac{\ovl m_{n}}{c},\cdot\cdot),\\
&&\hspace{-30pt}
e_{n-1}^c(m)=(\cd,m_{n-1}\frac{c}{\xi_{n-1}},
m_n\frac{\xi_{n-1}^2}{c^2},\ovl m_n\xi_{n-1}^2,
\frac{\ovl m_{n-1}}{\xi_{n-1}},\ld),\\
&&\hspace{-30pt}
e_i^c(m)=(\cd,\frac{cm_i}{\xi_{i}},
\frac{\xi_{i}m_{i+1}}{c},
\cd,\xi_{i}\ovl m_{i+1},\frac{\ovl m_i}
{\xi_{i}},\ld)\q(1\leq i\leq n-2),
\end{eqnarray*}
where
\[
\xi_i\seteq
\begin{cases}
\frac{c^2\ovl m_1+m_1}{\ovl m_1+m_1}&
\text{ for }i=0,\\
\frac{c\ovl m_{i+1}+m_{i+1}}
{\ovl m_{i+1}+m_{i+1}}
&\text{ for }i\ne 0,n-1\\
\frac{c+M(m)}{1+M(m)}&\text{ for }i=n-1.
\end{cases}
\]
Let $\eta\cl\cB_L(\TY(A,2,2n))
\longrightarrow X_4^{(n,L^2)}$
($m\mapsto l$)
be the morphism defined by
\[
l_1=\ovl l_1=m_0,\q
l_i=m_{i-1},\q \ovl l_i=\ovl m_{i-1} 
(i=2,\ld, n),\q
l_{n+1}=\frac{m_n}{1+M(m)},\q
\ovl l_{n+1}=\frac{\ovl m_n}{1+\ovl M(m)},\q
\]
where $l_i,\ovl l_i$ ($i=n+2,\ld,2n+2$) 
are uniquely determined by (\ref{aen-eq})
and 
then, $\eta$ is an isomorphism of 
$\TY(A,2,2n)$-
geometric crystal.
\begin{pro}
We have the following isomorphisms of 
$\TY(A,2,2n)$-geometric crystals:
\begin{equation}
\cV(\TY(A,2,2n))_L\mapleft{\sim}
\cB_{L}(\TY(A,2,2n))\mapright{\sim}
X_4^{(n,L^2)}.
\label{3iso-aen}
\end{equation}
\end{pro}
{\sl Proof.}
The second isomorphism in (\ref{3iso-aen})
is given by $\eta$. 
Then let us see the first one.
Define $\Xi\cl\cB_{L}(\TY(A,2,2n))
\to\cV(\TY(A,2,2n-1))_L$ ($m\mapsto x$) to be 
\[
x_0=\frac{1}{m_0},\q
x_i\seteq \frac{1}{m_0^2\ovl m_1
\ovl m_2\cd\ovl m_i},\q
\ovl x_i=\frac{m_1\cd m_i}{L^2}\,\,
(i=1,\ld, n-1),\q
x_n=\frac{1}
{\ovl m_n(m_0^2\ovl m_1\cd \ovl m_{n-1})^2},
\]
and the inverse $\Xi^{-1}$ is 
\[
 m_0=\frac{1}{x_0},\,\,m_1=L^2\ovl x_1,\,\,
\ovl m_1=\frac{x_0^2}{x_1},\,\,
m_i=\frac{\ovl x_i}{\ovl x_{i-1}},\,\,
\ovl m_i=\frac{x_{i-1}}{x_i}\,\,
(2\leq i\leq n-1),\,\,
m_n=\frac{x_n}{\ovl x_{n-1}^2},\,\,
\ovl m_n=\frac{x_{n-1}^2}{x_n}.
\]
Then, by direct calculations, we can check 
that $\Xi$ is an isomorphism of 
geometric crystals.\qed

%%%%%%% Section 7 %%%%%%%%%%%%%%
%%%%%%%%%%%%%%%%%%%%%%%
\section{\bf
Product Structure on Affine Geometric Crystals}
%%% see note X pp84-86

In general, if $\bbX_1$ and $\bbX_2$ are 
geometric crystals induced from unipotent 
crystals, the product $\bbX_1\times \bbX_2$ 
possesses a geometric crystal structure
(\cite{BK},\cite{N}). 
More precisely, let 
$\bbX_1=(X,\{e_i\},\{\gamma_i\},\{\vep_i\})$ and
$\bbX_2=(Y,\{e_i\},\{\gamma_i\},\{\vep_i\})$
be geometric crystals. For $x\in X$ and $y\in Y$
set 
\begin{eqnarray}
&&\gamma_i(x,y)\seteq \gamma_i(x)\gamma_i(y),
\label{pro1}\\
&&\vep_i(x,y)\seteq \vep_i(x)+
\frac{\vep_i(x)\vep_i(y)}{\vp_i(x)}\q
(\vp_i(x)=\gamma(x)\vep_i(x)),
\label{pro2}\\
&&e_i^c(x,y)\seteq(e_i^{c_1}(x),e_i^{c_2}(y))
\qq\text{where }
c_1\seteq\frac{c\vp_i(x)+\vep_i(y)}
{\vp_i(x)+\vep_i(y)},\q
c_2\seteq \frac{c}{c_1}.
\label{pro3}
\end{eqnarray}
\begin{thm}[\cite{BK,N}]
\label{unip}
Suppose that the geometric crystals 
$\bbX_1$ and $\bbX_2$ are induced from unipotent 
crystals. Then,
\eqref{pro1}--\eqref{pro3}
endow the product $X\times Y$ 
with a structure of a
geometric crystal.
Moreover, 
if $\bbX_1$ and $\bbX_2$ are positive, 
then $\bbX_1\times \bbX_2$ is positive and 
we have the 
isomorphism of crystals:
\begin{equation}
{\mathcal UD}(\bbX_1\times\bbX_2)
\cong {\mathcal UD}(\bbX_1)\ot{\mathcal UD}(\bbX_2).
\label{times-ud-iso}
\end{equation}

\end{thm}

As for the affine geometric crystal 
$\cV(\ge)_l$ in Sect. \ref{sec5}, 
its data $e_i,\gamma_i,\vep_i$ $(i\in I\setminus\{0\})$
are obtained from the ones of the
geometric crystal $B_{\bf i}^-\cdot l^H$ 
as in \ref{schubert} which is
induced from the unipotent crystal on 
some $X_w\times l^H$ where 
${\bf i}$ is 
a reduced word for $w$ and $X_w$ is the
Schubert cell associated with $w\in W$.
%and note 
%that $\{l^H\}$ has a trivial unipotent crystal 
%structure. 
Thus, by Theorem \ref{unip}
we have a $\ge_0$-
geometric crystal structure on 
$\cV(\ge)_{L_1}\times\cd
\times\cV(\ge)_{L_k}$.
\begin{thm}
\label{prod-thm}
For any $k\in\bbZ_{\geq0}$ and 
$L_1,\cd,L_k\in \bbC^\times$, 
the product 
$\cV(\ge)_{L_1}\times\cd
\times\cV(\ge)_{L_k}$
possesses an affine 
geometric crystal structure. 
\end{thm}
{\sl Proof.}
By the argument above and Theorem \ref{unip}, 
it is enough to 
check the conditions in Definition \ref{def-gc}
related to $i=0$. 
First, let us check 
$\gamma_j(e_i^c(x_1,\cd,x_k))=
c^{a_{ij}}\gamma_j(x_1,\cd,x_k)
\q \text{for }(i,j)=(i,0),\,\,(0,i).$
%Let us assume $\ge\ne\TY(A,2,2n)$. 
For $x_j\in\cV(\ge)_{L_j}$ and $i\in I$,
we have
\begin{eqnarray*}
&&\gamma_0(e_i^c(x_1,\cd,x_k))
=\gamma_0(e_i^{c_1}(x_1))\cd\gamma_0(e_i^{c_k}(x_k))
=c_1^{a_{i0}}\gamma_0(x_1)\cd 
c_k^{a_{i0}}\gamma_0(x_k)
=c^{a_{i0}}\gamma_0(x_1,\cd,x_k),\\
&&\gamma_i(e_0^c(x_1,\cd,x_k))=
\gamma_i(e_0^{c_1}(x_1))\cd\gamma_i(e_0^{c_k}(x_k))
=c_1^{a_{0i}}\gamma_i(x_1)\cd
c_k^{a_{0i}}\gamma_i(x_k)
=c^{a_{0i}}\gamma_i(x_1,\cd,x_k),
\end{eqnarray*}
where $c_1,\cd,c_k$ are obtained by using 
(\ref{pro3}) repeatedly.

Next, let us check the relation 
$\vep_0(e_0^c(x_1,\cd,x_k))
=c^{-1}\vep_0(x_1,\cd,x_k)$ by the induction on $k$.
Assume $\vep_0(e_0^c(x'))
=c^{-1}\vep_0((x'))$ and 
$\vp_0(e_0^c(x'))
=c\vp_0((x'))$ 
where $x'=(x_1,\cd,x_{k-1})$. 
\begin{eqnarray*}
\vep_0(e_0^c(x_1,\cd,x_k))&=&
\vep_0(e_0^{c_1}(x'),e_0^{c_2}(x_k))
=\vep_0(e_0^{c_1}(x'))
+\frac{\vep_0(e_0^{c_1}(x'))
\vep_0(e_0^{c_2}(x_k))}
{\vp_0(e_0^{c_1}(x'))}\\
&=&c_1^{-1}\vep_0(x')
+\frac{\vep_0(x')\vep_0(x_k)}
{cc_1\vp_0(x')}
=c_1^{-1}\left(
\vep_0(x')
+\frac{\vep_0(x')\vep_0(x_k)}
{c\vp_0(x')}\right)\\
&=&\frac{\vp_0(x')
+\vep_0(x_k)}{c\vp_0(x')+\vep_0(x_k)}
\cdot
\frac{c\vp_0(x')+\vep_0(x_k)}
{c\vp_0(x')}\cdot
\vep_0(x')
=c^{-1}\left(
\vep_0(x')
+\frac{\vep_0(x')\vep_0(x_k)}
{\vp_0(x')}\right)\\
&=&c^{-1}\vep_0(x_1,\cd,x_k).
\end{eqnarray*}
Finally, let us check the Verma relations:
we need to see the following cases
\begin{enumerate}
\item
The case 
$\ge=\TY(A,1,n)$
$\begin{cases}
e_0^ae_i^{ab}e_0^b=e_i^be_0^{ab}e_i^a&
\text{ if }i=1,n\\
e_0^ce_i^d=e_i^de_0^c&\text{ if }i\ne1,n.
\end{cases}
$
\item
The case 
$\ge=\TY(B,1,n),\TY(D,1,n),\TY(A,2,2n-1)$.
$\begin{cases}
e_0^ae_i^{ab}e_0^b=e_i^be_0^{ab}e_i^a&
\text{ if }i=2\\
e_0^ce_i^d=e_i^de_0^c&\text{ if }i\ne2.
\end{cases}
$
%\item
%The case 
%$\ge={\TY(A,2,2n)}^\dagger$.
%$\begin{cases}
%e_i^ae_0^{a^2b}e_i^{ab}e_0^b
%=e_0^be_i^{ab}e_0^{a^2b}e_i^a&
%\text{ if }i=1\\
%e_0^ce_i^d=e_i^de_0^c&\text{ if }i\ne1.
%\end{cases}$
\item
The case 
$\ge=\TY(D,2,n+1)$, $\TY(A,2,2n)$.
$\begin{cases}
e_0^ae_i^{a^2b}e_0^{ab}e_i^b
=e_i^be_0^{ab}e_i^{a^2b}e_0^a&
\text{ if }i=1\\
e_0^ce_i^d=e_i^de_0^c&\text{ if }i\ne1.
\end{cases}$
\end{enumerate}
By the result in \cite{KOTY}, we have 
the product structure on $\{\cB_l(\TY(A,1,n))\}$ 
and $\{\cB_l(\TY(D,1,N))\}$.
Since $\cV(\TY(A,1,n))_l\cong\cB_l(\TY(A,1,n))$
and $\cV(\TY(D,1,n))_l\cong\cB_l(\TY(D,1,n))$,
we also have the product structure
on $\{\cV(\TY(A,1,n))_l\}$ 
and $\{\cV(\TY(D,1,n))_l\}$.
Hence,  
we have the Verma relations on their product.
Thus, we obtain the case (i) and (ii) $\ge=\TY(D,1,n)$.

For the case $\sigma(i)\ne0$, it is easy to check
the relation. Indeed, {\it e.g.}, 
for the case (iii) $i=1$, we have:
\begin{eqnarray*}
%&&{\rm (i)}\,\, i\ne1,n\,\,:
%e_0^ce_i^d=\osigma^{-1}e_1^c\osigma e_i^d
%=\osigma^{-1}(e_1^c\osigma e_{i}^d
%\osigma^{-1})\osigma
%=\osigma^{-1} (e_1^c e_{i+1}^d)\osigma
%=\osigma^{-1} (e_{i+1}^d e_1^c )\osigma
%=e_i^d e_0^c,\\
&&
e_0^ae_1^{a^2b}e_0^{ab}e_1^b
=(\osigma^{-1}e_{n}^a\osigma)
(\osigma^{-1} e_{n-1}^{a^2b}\osigma)
(\osigma^{-1}e_n^{ab}\osigma)
(\osigma^{-1} e_{n-1}^b\osigma)
=\osigma^{-1}(e_n^ae_{n-1}^{a^2b}e_n^{ab}e_{n-1}^b)
\osigma\\
&&\qq\qq
=\osigma^{-1}(e_{n-1}^be_{n}^{ab}e_{n-1}^{a^2b}
e_{n}^a)\osigma
=e_1^be_0^{ab}e_1^{a^2b}e_0^a.
\end{eqnarray*}
Similarly, the other cases with 
$\sigma(i)\ne0$ can be shown.

To complete (ii) and (iii), 
it suffices to check the cases $\sigma(i)=0$
: {\it i.e.,}
(ii) $i=1$, (iii) $i=n$. 
\\
In the previous section 
we see that for $\ge\ne\TY(A,1,n)$ the 
geometric crystal $\cV(\ge)_l$ is obtained from 
the geometric crystal $\cB(\TY(D,1,N))_{l'}$
($l'=l$ or $l^2$)
by the method of foldings.
Thus, in the case (ii) we have 
$(e_0^{\ge})^c=(e_0^{\TY(D,1,N)})^c$ and 
$(e_1^{\ge})^c=(e_1^{\TY(D,1,N)})^c$ for 
$\ge=\TY(B,1,n), \TY(A,2,2n-1)$. 
Since $(e_0^{\TY(D,1,N)})^c
(e_1^{\TY(D,1,N)})^d=(e_1^{\TY(D,1,N)})^d
(e_0^{\TY(D,1,N)})^c$, we have 
$(e_0^{\ge})^c(e_1^{\ge})^d=
(e_1^{\ge})^d(e_0^{\ge})^c$ 
and then we completed (ii).
In the case 
(iii), we have 
$(e_0^{\ge})^c=(e_0^{\TY(D,1,N)})^c$ and 
\[
(e_n^{\ge})^c=\begin{cases}
(e_{n+1}^{\TY(D,1,n+2)})^c\circ 
(e_{n+2}^{\TY(D,1,n+2)})^c &\ge=\TY(D,2,n+1),\\
(e_{n+1}^{\TY(D,1,2n+2)})^c&\ge=\TY(A,2,2n).
\end{cases}
\]
Then $(e_0^{\ge})^c(e_n^{\ge})^d
=(e_n^{\ge})^d(e_0^{\ge})^c$, which 
completes (iii).\qed

By Theorem \ref{ud-geo}, Theorem \ref{unip} 
and Theorem \ref{prod-thm}, 
we obtain
\begin{cor}
\label{multi-ud}
For positive real numbers $L_j$ 
($j=1,\cd,k$), 
let $\theta_{L_j}:T'\to\cV(\ge)_{L_j}$ be 
the positive structure as in the previous section. Then 
\[
\Theta\seteq(\theta_{L_1},\cd,\theta_{L_k})
:{T'}^{\times k}\to\cV(\ge)_{L_1}\times\cd
\times\cV(\ge)_{L_k}
\] 
defines a positive structure 
on $\cV(\ge)_{L_1}\times\cd
\times\cV(\ge)_{L_k}$ and 
we have the 
isomorphism of geometric crystals:
\begin{equation}
{\mathcal UD}_\Theta(\cV(\ge)_{L_1}\times\cd
\times\cV(\ge)_{L_k})\cong
B_\ify(\ge^L)^{\ot k}.
\end{equation}
\end{cor}

%%%%%%%%%% Section 8 %%%%%%%%%%%
%%%%%%%%%%%%%%%%%%%%%%%%%%%%%%%%%%%%%%%%%%%%%%%%%%
\section{\bf M-matrices and Automorphisms}\label{Mmat-sec}
\setcounter{equation}{0}
\renewcommand{\theequation}{\thesection.\arabic{equation}}

\subsection{Definition of M-matrices}
An M-matrix is an important object to realize 
\tropical R map
from geometric crystals, though 
its exact definition is not yet fixed.
In this paper, we take the following as a 
temporary definition of M-matrices.
\begin{df}
Let $\ge$ be an affine Lie algebra and 
$\cB_L\subset (\bbC^\times)^k$ 
be a certain $\ge$-geometric crystal depending
on $L\in\bbC^\times$.
For $x\in\cB_L$ and an indeterminate $z$, 
a square matrix $M_L(x,z)$ with entries in
$\bbC(x,z)$ is an {\it M-matrix}
if 
\begin{enumerate}
\item
For $x\in \cB_L$, $i\in I$ and $c\in\bbC^\times$, 
there exist non-singular
matrices $X_i(x,c,z)$ and 
$Y_i(x,c,z)$ whose entries are in
$\bbC(x,c,z)$ and satisfying
\[
M_L({e_i^c(x)},z)=
%x_i(\frac{c-1}{z^{\del_{i0}}\vep_i(x)})
X_i(x,c,z)M_L(x,z)Y_i(x,c,z).
%x_i(\frac{c^{-1}-1}{z^{\del_{i0}}\vp_i(x)})
\]
%(Then, $\cM_L$ holds GC-structure.)
\item
For a given $(x,y)\in \cB_L\times \cB_K$ and 
$L,K\in \bbC^\times$, 
the solution $(x',y')\in \cB_K\times \cB_L$ 
of the equation
$$
M_L(x,z)M_K(y,z)=M_K(x',z)M_L(y',z),\,
$$
uniquely exists for any $z$ and the correspondence 
$(x,y)\mapsto (x',y')$ defines a 
birational map \\
$R:\cB_L\times \cB_K\to
\cB_K\times \cB_M$.
\item
Let $x_j,y_j$ be in 
$\cB_{L_j}$ $(j=1,2,\ld,m)$. Suppose 
$L_i\ne L_j$ ($i\ne j$) and 
\[
 M_{L_1}(x_1,z)\cd M_{L_m}(x_m,z)
= M_{L_1}(y_1,z)\cd M_{L_m}(y_m,z).
\]
Then $x_j=y_j$ $(j=1,\ld,m)$.
\end{enumerate}
\end{df}

%%%%%%%%%%%%%%%%%%%%%%%
%\newpage
\begin{ex}[\cite{KOTY}]
$\TY(A,1,n)$-case: \\
Let $\cB_L(\TY(A,1,n))$ be the geometric crystal
as in Sect.\ref{geo-a}. 
For 
$l=(l_1,\ld,l_{n+1})\in\cB_L(\TY(A,1,n))$, set
$$
M_L(l,z)\seteq
\begin{pmatrix}
l_1^{-1}&&&&&-z\\
-1&l_2^{-1}&&&&\\
&&\cd &&&\\
&&&-1&l_n^{-1}&\\
&&&&-1&l_{n+1}^{-1}
\end{pmatrix}^{-1}.
$$
\end{ex}
Then the matrix $M_L(l,z)$ is an M-matrix.
%\newpage
%%%%%%%%%%%%%%%%%%%%%%%%%%%%%%%%%%%%%%
\subsection{\bf Explicit form of the 
M-matrix of type $\TY(D,1,n)$}

An explicit form of the M-matrix $M_L(l,z)$
for $\cB_L(\TY(D,1,n))$
is described \cite{KOTY} 
where the geometric crystal
$\cB_L(\TY(D,1,n))$ is as in Sect.\ref{fold}.

For $l=(l_1,\ld,
l_n,\ovl l_{n-1},\ld,\ovl l_1)\in
\cB_L(\TY(D,1,n))$, the M-matrix $M_L(l,z)$ 
is given
as follows \cite{KOTY}:
$M(l,z)$ is a $2n\times 2n$ matrix in the form
$M_L(l,z)=A(l)+zB(l)+z^2 C(l)$ where each 
matrix $A(l), B(l)$ and $C(l)$ are given by 
$C(l)=E_{1,2n}$, 
\begin{eqnarray}
&&A(l)_{i,i}=\begin{cases}
\frac{l_i}{\ovl l_i}&1\leq i\leq n-1,\\
l_n&i=n,\\
\frac{1}{l_n}&i=n+1,\\
\frac{\ovl l_{2n+1-i}}{l_{2n+1-i}}
&n+2\leq i\leq 2n,
\end{cases}\\
&&A(l)_{i,j}=l_j\cd l_{i-1}
\left(1+\frac{l_i}{\ovl l_i}\right)
\q 1\leq j<i\leq n-1,
\end{eqnarray}
\begin{eqnarray}
&&A(l)_{2n+1-j,2n+1-i}=\ovl l_j\cd\ovl l_{i-1}
\left(1+\frac{\ovl l_i}{l_i}\right)\q
1\leq j<i\leq n-1,\\
&&A(l)_{n,j}=l_j\cd l_n\q 1\leq j\leq n-1,\\
&&A(l)_{2n+1-j,n}=\ovl l_j\cd \ovl l_{n-1}l_n
\q 1\leq j\leq n-1,\\
&&A(l)_{2n+1-j,n+1}=\ovl l_j\cd\ovl l_{n-1}
\q 1\leq j\leq n-1,\\
&&A(l)_{2n+1-i,j}=(l_j\cd l_n)
(\ovl l_i\cd \ovl l_{n-1})\qq 1\leq i,j\leq n-1,\\
&&A(l)_{i,j}=0\qq \text{otherwise}.
\end{eqnarray}
The matrix element $B(l)_{i,j}$ is given by
$B(l)_{i,j}=A(l)_{i,1}A(l)_{2n,j}-LA(l)_{i,j}
-\del_{i,1}\del_{j,2n}$ (\cite{KOTY}).

%%%%%%%%%%%%%%%%%%%%%%%%%%%%%%%%%%%%%%%
%%\newpage
\subsection{Matrix Realization}

In terms of {M-matrices}, 
the involutions $\Sigma_0^{(n)}, 
\Sigma_1^{(n)}, \Sigma_2^{(n)}$
in Sect.\ref{fold} are realized by 
the adjoint action of certain matrices.
\begin{pro}
For each involution 
$\Sigma_i^{(n)}$ $(i=0,1,2)$,
there exists a non-singular 
matrix $J^{(n)}_i=J^{(n)}_i(z)$  such that
$$
M_L(\Sigma_i^{(n)}(l), z)
=J^{(n)}_i(z)M_L(l,z){J^{(n)}_i(z)}^{-1},
$$
where the matrix $J^{(n)}_i$ 
does not depend on $l$.
Each $J_i^{(n)}$ is as follows:
\begin{eqnarray*}
J^{(n)}_0=
\left(
\begin{tabular}{c|c|c}
0&0&z\\
\hline0&$E_{2n}$
%{$\begin{array}{ccccc}
%1&&&&\\&1&&&\\&&\ddots&&\\
%&&&1&\\&&&&1
%\end{array}$}
&0\\
\hline $z^{-1}$&0&0
\end{tabular}
\right),
\q
J^{(n)}_1= 
 \left(
\begin{tabular}{c|c|c}
%{$\begin{array}{ccc}
%1&&\\&\ddots&\\
%&&1
%\end{array}$}
$E_n$&0&0\\
\hline0&
{$\begin{array}{cc}
0&1\\1&0
\end{array}$}
&0\\
\hline 0&0&$E_n$
%{$\begin{array}{ccc}
%1&&\\&\ddots&\\
%&&1
%\end{array}$}
\end{tabular}
\right),
%\end{eqnarray*}
\q
J^{(n)}_2=
 \left(
\begin{tabular}{c|c}
{0}&$z\cdot E_{2n}$
%{$\begin{array}{ccccc}
%z&&&&\\&z&&&\\&&\ddots&&\\
%&&&z&\\&&&&z
%\end{array}$}
\\
\hline
%{$\begin{array}{ccccc}
%1&&&&\\&1&&&\\&&\ddots&&\\
%&&&1&\\&&&&1
%\end{array}$}
$E_{2n}$&{0}
\end{tabular}\right),\q
\begin{array}{c}
%J_2=J_0J_1,\q J_4=J_0J3.\\ \\
%\text{Note: 
%$\forall J_i$ does not depend on $L$.}
\end{array}
\end{eqnarray*}
where $E_m$ is the identity matrix of 
size $m$. 
The size of each matrix 
$J^{(n)}_i$ and $M_L(l,z)$ is
$2n+2$ for $i=0$ and $i=1$, 
and $4n$ for $i=2.$
Setting $J^{(n)}_3=J^{(n+1)}_0J^{(n+1)}_1
(=J^{(n+1)}_1J^{(n+1)}_0)$ 
and $J^{(n)}_4
=(J^{(2n+1)}_0J^{(2n+1)}_1)J^{(n+1)}_2
(=J^{(n+1)}_2(J^{(2n+1)}_0J^{(2n+1)}_1))$, 
we also have 
\begin{equation}
M_L(\Sigma_k^{(n)}(l), z)=J^{(n)}_k
M_L(l,z){J^{(n)}_k}^{-1}
\q(k=3,4).
\label{34}
\end{equation}
\end{pro}
{\sl Proof.} 
Since the involution $\Sigma^{(n)}_0$ 
(resp. $\Sigma^{(n)}_1$) coincides with 
the involution $\sigma_1$ (resp. $\sigma_n$)
in \cite{KOTY} and 
the matrix $J^{(n)}_0$ (resp. $J^{(n)}_1$) 
is identified with the matrix $J_1(z)$ 
(resp. $J_n(z)$) as in \cite{KOTY}. 
Thus, Lemma 3.10 in \cite{KOTY} shows that 
our assertions for $\Sigma^{(n)}_k$ $(k=0,1)$
are right. Then let us show the case $k=2$.
Since we have the explicit form of $M_L(l,z)$
as in the last subsection,
it is carried out by case-by-case calculations.
For example, let us see the diagonal entries:
\begin{eqnarray*}
M_L(\Sigma^{(n)}_2(l),z)_{i,i}
=\begin{cases}
1/l_{2n}&i=1,\\
\frac{\ovl l_{2n+1-i}}{l_{2n+1-i}}+zL
&2\leq i\leq 2n-1,\\
\ovl l_1/l_1&i=2n,\\
l_1/\ovl l_1&i=2n+1,\\
\frac{l_{i-2n}}{\ovl l_{i-2n}}+zL
&2n+2\leq i\leq 4n-1,\\
l_{2n}&i=4n.
\end{cases}
\end{eqnarray*}
Here note that the matrix $M_L(l,z)$ is a 
$4n\times 4n$-matrix.

On the other-hand, the diagonal entries of 
the matrix 
\[
M'\seteq
J^{(2)}_nM_L(l,z){J^{(2)}_n}^{-1}
=J^{(2)}_n\begin{pmatrix}M_1&M_2\\M_3&M_4
\end{pmatrix}{J^{(2)}_n}^{-1}=
\begin{pmatrix}M_4&zM_3\\ z^{-1}M_2&M_1
\end{pmatrix}
\]
 are given by:
\begin{eqnarray*}
&&\text{For }i=1,\ld,2n,\qq\q
M'_{i,i}=M_L(l,z)_{i+2n,i+2n}=\begin{cases}
1/l_{2n}&i=1,\\
\frac{\ovl l_{2n+1-i}}{l_{2n+1-i}}+zL
&2\leq i\leq 2n-1,\\
\ovl l_1/l_1&i=2n,
\end{cases}\\
&&\text{For }i=2n+1,\ld,4n,\qq\q
M'_{i,i}=M_L(l,z)_{i-2n,i-2n}=\begin{cases}
l_1/\ovl l_1&i=2n+1,\\
\frac{l_{i-2n}}{\ovl l_{i-2n}}+zL
&2n+2\leq i\leq 4n-1,\\
l_{2n}&i=4n.
\end{cases}
\end{eqnarray*}
Then, we have $M_L(\Sigma^{(2)}_n(l),z)_{i,i}
=M'_{i,i}
(=(J^{(2)}_nM_L(l,z){J^{(2)}_n}^{-1})_{i,i})$.
The other cases are also obtained similarly.
\qed
%%%%%%%%%%%%%%%%%%%%%%%%%%%%%%%%%%%%%%%
\subsection{Birational maps on 
fixed point varieties }\label{birat-fix}

Let $\cR_{LK}$ be the birational map
defined by 
\[
\begin{array}{ccccc}
\cR_{LK}:&\cB_L(\TY(D,1,N))\times 
\cB_K(\TY(D,1,N))&\to &
\cB_K(\TY(D,1,N))\times \cB_L(\TY(D,1,N))
&(L,K\in\bbC^\times)\\
&(l,m)&\mapsto&(l',m'),&
\end{array}
\]
where $(l',m')$ be  the unique solution of the
equation $M_L(l,z)M_K(m,z)=M_K(l',z)M_L(m',z)$.
Let $X_i^{(n,L)}$ be one of the
fixed point subvarieties
in $\cB_L(\TY(D,1,N))$. 
\begin{thm}
\label{rxx}
Let us denote $\cR^{(i)}_{LK}$ the restriction 
of $\cR_{LK}$ on 
$X_i^{(n,L)}\times X_i^{(n,K)}$.
Then 
$\cR^{(i)}_{LK}$ is a well-defined 
birational map 
$X_i^{(n,L)}\times X_i^{(n,K)}\to 
X_i^{(n,K)}\times X_i^{(n,L)}$ 
$(L,K\in \CC^\times)$.
\end{thm}
{\sl Proof.}
For $(l,m)\in X_i^{(n,L)}\times X_i^{(n,K)}$, 
set 
$(l',m')\seteq \cR(l,m)$, {\it i.e.,}
\begin{equation}
M_L(l,z)M_K(m,z)=M_K(l',z)M_L(m',z).
\label{a}
\end{equation}
 We have
\begin{eqnarray*}
M_L(l,z)M_K(m,z)&=&
M_L(\Sigma_i^{(n)}(l),z)M_K(\Sigma_i^{(n)}(m),z)
\\
&=&J_i^{(n)}M_L(l,z)M_K(m,z){J_i^{(n)}}^{-1}\\
&=&J_i^{(n)}M_K(l',z)M_L(m',z){J_i^{(n)}}^{-1}\\
&=&M_K(\Sigma_i^{(n)}(l'),z)
M_L(\Sigma_i^{(n)}(m'),z).
\end{eqnarray*}
Then, we have 
\[
 M_K(l',z)M_L(m',z)
=M_K(\Sigma_i^{(n)}(l'),z)
M_L(\Sigma_i^{(n)}(m'),z).
\]
It follows from
the uniqueness of the solution 
for (\ref{a}) that 
$(l',m')=(\Sigma_i^{(n)}(l'),
\Sigma_i^{(n)}(m'))$ and then we have
$(l',m')\in X_i^{(n,K)}\times X_i^{(n,L)}$.
\qed

%\newpage
%%%%%%%%%%%%%%%%%%%%%%%%%%%%%%%%%%%%%%%%%
%%%%%%%%%% section 9 %%%%%%%%%%%%%%%%%%%%%%%%%%%
\section{\bf \Tropical R 
Maps}\label{tro-r}
\setcounter{equation}{0}
\renewcommand{\theequation}{\thesection.\arabic{equation}}

In this section, we define the notion of 
\tropical R map and give explicit forms of the 
affine \tropical R maps 
on the geometric crystals 
constructed above.
%%%%%%%%%%%%%%%%%%%%%%%%%%%%%%%%
\subsection{Definition of \tropical R map}
\begin{df}
\label{trop-r}
Let $\{(X_{\lm},\{e_i^\lm\},
\{\gamma_i^\lm\},\{\vep_i^\lm\}
)\}_{\lm\in\Lm}$ be a family 
of geometric crystals equipped with the product 
structures, where $\Lm$ is a certain
index set and its element is called a
{\it spectral parameter}. 
A birational map 
${\mathcal R_{\lm\mu}}:X_\lm\times X_\mu\longrightarrow 
X_\mu\times X_\lm$ $(\lm,\mu\in\Lm)$
is said to be a {\it \tropical R map}
(or shortly, tropical R)
if it satisfies the following conditions:
\begin{eqnarray}
&&(e^{X_\mu\times X_\lm}_i)^c
\circ{\mathcal R_{\lm\mu}}
={\mathcal R_{\lm\mu}}
\circ (e^{X_\lm\times X_\mu}_i)^c,
\label{R-e}\\
&&\vep^{X_\lm\times X_\mu}_i
=\vep^{X_\mu\times X_\lm}_i
\circ{\mathcal R}_{\lm\mu},\label{R-ep}\\
&&\gamma_i^{X_\lm\times X_\mu}
=\gamma_i^{X_\mu\times X_\lm}
\circ{\mathcal R}_{\lm\mu},\label{R-g}\\
&&{\mathcal R}^{(12)}
{\mathcal R}^{(23)}{\mathcal R}^{(12)}
={\mathcal R}^{(23)}
{\mathcal R}^{(12)}{\mathcal R}^{(23)}
\q\text{on }
X_{\lm}\times X_{\mu}\times X_{\nu},
\label{YB},\\
&&{\mathcal R}_{\mu\lm}{\mathcal R}_{\lm\mu}
=\operatorname{id}_{\lm\mu}.\label{inver}
\end{eqnarray}
for any $i\in I$ 
and any $\lm,\mu,\nu\in\Lm$. 
Here ${\mathcal R}^{(ij)}$ 
means that it acts on 
$i$-th and $j$-th components of the product.
\end{df}
In the rest of this section, we give the 
explicit forms of the
\tropical $\mathcal R$ for the affine 
geometric crystals.
%%%%%%%%%%%%%%%%%%%%%%%%%%%%%%%%%%%%%%%%%%
\subsection{$D^{(1)}_n$ case $(n\geq4)$}%8.2

Let $\cR_{LK}:\cB_L(\TY(D,1,N))\times 
\cB_K(\TY(D,1,N))\to
\cB_K(\TY(D,1,N))\times \cB_L(\TY(D,1,N))$
$(L,K\in\bbC^\times)$ be the birational 
map as in \ref{birat-fix}.
In \cite{KOTY}, it is shown that 
the morphism $\cR_{LM}$ is 
a tropical R map for 
$\{\cB_L(\TY(D,1,n))\}_{L\in\bbC^\times}$.
We shall describe the explicit form of 
$\cR_{LM}$.
Let $\sharp:\cB_L(\TY(D,1,n))\rightarrow 
\cB_L(\TY(D,1,n))$
be an involution defined by 
\begin{equation}
\sharp(l_1,l_2,\ld,l_n,\ovl l_{n-1},
\ld,\ovl l_2,\ovl l_1)
=(\ovl l_1,l_2,\ld,l_n,
\ovl l_{n-1},\ld,\ovl l_2,l_1),
\end{equation}
that is, $\sharp:l_1\leftrightarrow \ovl l_1$
and $*:\cB_L(\TY(D,1,n))\times \cB_M(\TY(D,1,n))
\rightarrow \cB_M(\TY(D,1,n))\times 
\cB_L(\TY(D,1,n))$ 
an involution defined by 
\[
((l_1,l_2,\ld,\ovl l_2,\ovl l_1),
(m_1,m_2,\ld,\ovl m_2,\ovl m_1))^*
=((\ovl m_1,\ovl m_2,\ld,m_2,m_1),
(\ovl l_1,\ovl l_2,\ld, l_2,l_1))
\]
that is, 
$*:l_i\leftrightarrow \ovl m_i, \q
  \ovl l_i\leftrightarrow m_i\,\,(1\leq i\leq n-1),\q
l_n\leftrightarrow m_n.$

Following \cite{KOTY}, 
we define the rational functions 
$V_i$ $(i=0,1,\ld,n-1)$ and $W_i$ 
($i=1,\ld,n-1$) on 
$\cB_L(\TY(D,1,n))\times
\cB_M(\TY(D,1,n))$ $(L,M\in\bbC^\times)$ by 
\[
 W_i:=V_i{V_i}^*
+(M-L){V_i}^*+(L-M)V_i
\q(1\leq i\leq n-2),\q
W_{n-1}:=V_{n-1}{V_{n-1}}^*.
\]
\begin{equation}
V_i=
\sum_{j=1}^{n-2}(\theta_{i,j}(l,m)
+\theta'_{i,j}(l,m))
+\sum_{j=1}^{n}(\eta_{i,j}(l,m)
+\eta'_{i,j}(l,m)),
\end{equation}
where we set $L=l_1l_2\cd \ovl l_2\ovl l_1$,
$M=m_1m_2\cd \ovl m_2\ovl m_1$,
\begin{eqnarray}
&&\theta_{i,j}(l,m)
=\begin{cases}
\displaystyle L\prod_{k=j+1}^i
\frac{\ovl m_k}{\ovl l_k}&\text{for }
1\leq j\leq i,\\
\displaystyle M\prod_{k=i+1}^j
\frac{\ovl l_k}{\ovl m_k}&\text{for }
i+1\leq j\leq n-2,
\end{cases}\\
&&\theta'_{i,j}(l,m)=L\left(\prod_{k=1}^i
\frac{\ovl m_k}{\ovl l_k}\right)
\left(\prod_{k=1}^j\frac{m_k}{l_k}\right)
\q\text{for}\q j=1,\ld,n-2,
\end{eqnarray}
\begin{eqnarray}
&&\eta_{i,j}(l,m)=
\begin{cases}
\displaystyle L\left(\prod_{k=j+1}^i
\frac{\ovl m_k}{\ovl l_k}\right)
\left(\frac{\ovl m_j}{l_j}\right)&\text{for}\q 1\leq j\leq i,\\
\displaystyle M\left(\prod_{k=i+1}^j\frac{\ovl l_k}{\ovl m_k}\right)
\left(\frac{\ovl m_j}{l_j}\right)&\text{for}\q i+1\leq j\leq n-1,\\
\displaystyle M
\left(\prod_{k=i+1}^{n-1}\frac{\ovl l_k}{\ovl m_k}\right)l_n
&\text{for}\q j=n,
\end{cases}\\
&&\eta'_{i,j}(l,m)=
\begin{cases}
\displaystyle L\left(\prod_{k=1}^i\frac{\ovl m_k}{\ovl l_k}\right)
\left(\prod_{k=1}^j\frac{m_k}{l_k}\right)
\left(\frac{l_j}{\ovl m_j}\right)&\text{for}\q 1\leq j\leq n-1,\\
\displaystyle L\left(\frac{L}{M}\right)^{\del_{i,n-1}}
\left(\prod_{k=1}^i\frac{\ovl m_k}{\ovl l_k}\right)
\left(\prod_{k=1}^{n-1}\frac{m_k}{l_k}\right)
\left(\frac{1}{l_n}\right)&\text{for}\q j=n.
\end{cases}
\end{eqnarray}
{\sl Remark.}
Though $\TTY(W,n,i)$ seems
not to be a positive rational function, 
it is, indeed, positive by the formula
( \cite{KOTY} Lemma 4.14):
\[
 \left(\frac{1}{l_i}+\frac{1}{\ovl m_i}\right)
\TTY(W,n,i)=\frac{1}{m_i}
\TTY(V,n,i){\STY(V,*,i-1)}+\frac{1}{\ovl l_i}
\TTY(V,n,i-1){\STY(V,*,i)}
\]
Now, we introduce the 
following rational transformation 
$\mathcal R$ on 
$(\bbC^\times)^{2n-1}\times(\bbC^\times)^{2n-1}$:
\begin{eqnarray}
&&{\mathcal R}(l,m)=(l',m'),\nn\\
{\rm where}
 &&l'_1=
m_1\frac{{V_0^{\sharp}}}{V_1}, \q 
\ovl l'_1=\ovl m_1\frac{V_0}{V_1},\q
l'_i=m_i\frac{V_{i-1}W_i}
{V_iW_{i-1}},
\q \ovl l'_i=\ovl m_i\frac{V_{i-1}}{V_i}
\q(2\leq i\leq n-1),\nn\\
&&l'_n=m_n\frac{V_{n-1}}{{V^*_{n-1}}},\q
m'_1=l_1\frac{V_0}{{V^*_1}}, \q 
\ovl m'_1=\ovl l_1
\frac{V_0^\sharp}{{V^*_1}},\label{explicit-R}
\\
&&m'_i=l_i\frac{{V^*_{i-1}}}
{{V^*_i}},
\q \ovl m'_i=\ovl l_i
\frac{{V^*_{i-1}}W_i}
{{V^*_i}W_{i-1}}
\q(2\leq i\leq n-1),\q
m'_n=l_n\frac{{V^*_{n-1}}}
{V_{n-1}}.\nn
\end{eqnarray}
Here note that 
for $(l',m')={\mathcal R}(l,m)$ we have 
$l_1'l_2'\cd \ovl l_2'\ovl l_1'
=M$ and $m_1'm_2'\cd\ovl m_2'\ovl m_1'
=L$. 
Then 
$\mathcal R$ defines a rational map 
${\mathcal B}_L(D^{(1)}_n)\times
{\mathcal B}_M(D^{(1)}_n)\longrightarrow 
{\mathcal B}_M(D^{(1)}_n)
\times{\mathcal B}_L(D^{(1)}_n)$.
The following is one of the main results in \cite{KOTY}:
\begin{thm}[\cite{KOTY}]
The rational map ${\mathcal R}=\cR_{LM}:
{\mathcal B}_L(D^{(1)}_n)\times
{\mathcal B}_M(D^{(1)}_n)
\longrightarrow 
{\mathcal B}_M(D^{(1)}_n)\times
{\mathcal B}_L(D^{(1)}_n)$ gives a \tropical R map
of type $\TY(D,1,n)$
on the family of affine geometric crystals,
which will be denoted by $\cR(\TY(D,1,n))$
in the sequel.
\end{thm}
Here we describe the \tropical R on 
$\cV(\TY(D,1,n))_L\times\cV(\TY(D,1,n))_M$:
$\ovl\cR(x,y)\seteq
(\Xi,\Xi)\circ\cR\circ(\Xi^{-1},\Xi^{-1})(x,y)$
($x\in\cV(\TY(D,1,n))_L,\,y\in\cV(\TY(D,1,n))_M$).
We define the rational functions
$\ovl V_i$ $(i=0,1,\ld,n-1)$ and 
$\ovl W_i$ ($i=1,\ld,n-1$) on 
$\cV(\TY(D,1,n))_L\times
\cV(\TY(D,1,n))_M$ $(L,M\in\bbC^\times)$ by 
\[
\ovl W_i:=\ovl V_i{\ovl V^*_i}
+(M-L){\ovl V^*_i}+(L-M)\ovl V_i
\q(1\leq i\leq n-2),\q
\ovl W_{n-1}:=\ovl V_{n-1}
{\ovl V^*_{n-1}}.
\] 
\begin{equation}
\ovl V_i=
\sum_{j=1}^{n-2}(\ovl\theta_{i,j}(x,y)
+\ovl\theta'_{i,j}(x,y))
+\sum_{j=1}^{n}(\ovl\eta_{i,j}(x,y)
+\ovl\eta'_{i,j}(x,y)),
\end{equation}
where $\sharp$ is the involution on 
$\cV(\TY(D,1,n))_L$ and $*$ is the involution 
$\cV(\TY(D,1,n))_L\times \cV(\TY(D,1,n))_M
\to \cV(\TY(D,1,n))_M\times \cV(\TY(D,1,n))_L$
defined by
\begin{eqnarray*}
\sharp&:&x_i\mapsto \frac{x_i}{Lx_1\ovl x_1}\q
(i=1\cd, n)\q,
\ovl x_i\mapsto \frac{\ovl x_i}{Lx_1\ovl x_1}\q
(i=1,\cd,n-2),\\
*&:&x_i\mapsto \frac{1}{M\ovl y_i},\q
\ovl x_i\mapsto \frac{1}{M y_i}\q(i=1,\cd,n-2),\q
x_{n-1}\mapsto \frac{1}{M y_n},\q
x_n\mapsto \frac{1}{M y_{n-1}},\\
&&y_i\mapsto \frac{1}{L\ovl x_i},\q
\ovl y_i\mapsto \frac{1}{L x_i}\q(i=1,\cd,n-2),\q
y_{n-1}\mapsto \frac{1}{L x_n},\q
y_n\mapsto \frac{1}{L x_{n-1}},
\end{eqnarray*}
and 
\begin{eqnarray*}
&&\hspace{-30pt}\ovl\theta_{i,j}(x,y)\seteq
L\frac{x_iy_j}{x_jy_i}\q(1\leq j\leq i<n-1),\q
\ovl\theta_{n-1,j}(x,y)\seteq
L\frac{x_ny_j}{x_jy_n}\q(1\leq j<n-1),\\
&&\hspace{-30pt}\ovl\theta_{i,j}(x,y)\seteq
M\frac{x_iy_j}{x_jy_i}\q(i+1\leq j\leq n-2),\,\,
\ovl\theta'_{i,j}(x,y)\seteq
M\frac{x_i\ovl y_j}{\ovl x_jy_i}\,\,
(0\leq i\leq n-2),\,\,
\ovl\theta'_{n-1,j}(x,y)\seteq
M\frac{x_n\ovl y_j}{\ovl x_jy_n},\\
&&\hspace{-30pt}\ovl\eta_{i,j}(x,y)\seteq
L^{1-\del_{j,1}}
\frac{x_i\ovl x_{j-1}y_{j-1}}{x_j\ovl x_jy_i}
\q(1\leq j\leq i<n-1),\q
\ovl\eta_{n-1,j}(x,y)\seteq
L^{1-\del_{j,1}}
\frac{x_n\ovl x_{j-1}y_{j-1}}
{x_j\ovl x_jy_n}\q(1\leq j\leq n-1),\\
&&\hspace{-30pt}\ovl\eta_{i,j}(x,y)\seteq
\frac{M}{L^{\del_{j,1}}}
\frac{x_i\ovl x_{j-1}y_{j-1}}{x_j\ovl x_jy_i}\q
(0\leq i<j\leq n-2),\q
\ovl\eta_{i,n-1}(x,y)=
M\frac{x_i\ovl x_{n-2}y_{n-2}}
{x_{n-1}x_ny_i}\q(i+1\leq n-1),\\
&&\hspace{-30pt}\ovl\eta_{i,n}(x,y)\seteq
M\frac{x_iy_{n}}
{x_{n-1}y_i}\q(i\leq n-2),\q
\ovl\eta_{n-1,n}(x,y)\seteq 
M\frac{x_n}{x_{n-1}},\\
&&\hspace{-30pt}\ovl\eta'_{i,j}(x,y)\seteq
L^{\del_{j,1}}M
\frac{x_iy_j\ovl y_j}{\ovl x_{j-1}y_{j-1}y_i}\q
(1\leq j\leq n-2,\,i<n-1),\q
\ovl\eta'_{n-1,j}(x,y)\seteq
L^{\del_{j,1}}M
\frac{x_ny_j\ovl y_j}{\ovl x_{j-1}y_{j-1}y_n}\q
(1\leq j\leq n-2),\\
&&\hspace{-30pt}\ovl\eta'_{i,n-1}(x,y)\seteq
M\frac{x_iy_{n-1}y_n}{\ovl x_{n-2}y_{n-2}y_i}\q
(0\leq i<n-1),\q
\ovl\eta'_{n-1,n-1}(x,y)\seteq
M\frac{x_ny_{n-1}}{\ovl x_{n-2}y_{n-2}},\\
&&\hspace{-30pt}\ovl\eta'_{i,n}(x,y)\seteq
M\frac{x_iy_{n-1}}{y_ix_{n}}\q
(0\leq i<n-1),\q
\ovl\eta'_{n-1,n}(x,y)\seteq 
L\frac{y_{n-1}}{y_n},
\end{eqnarray*}
where we understand $x_0=y_0=\ovl x_0=\ovl y_0
=1$.
Note that in the rest of this section 
if we write 
$\ovl\theta_{ij}=\ovl\theta^{LM}_{i,j}$, 
$\ovl\theta'_{ij}={\ovl\theta'}^{LM}_{i,j}$,
$\ovl\eta_{ij}=\ovl\eta^{LM}_{i,j}$ and 
$\ovl\eta'_{ij}={\ovl\eta'}^{LM}_{i,j}$, 
we understand that 
$\ovl\theta_{i,j}^{LM}(x,y)^*
= \ovl\theta_{i,j}^{ML}((x,y)^*),$
${\ovl\theta'}_{i,j}^{LM}(x,y)^*
= {\ovl\theta'}_{i,j}^{ML}((x,y)^*),$
$\ovl\eta_{i,j}^{LM}(x,y)^*
= \ovl\eta_{i,j}^{ML}((x,y)^*)$ and 
${\ovl\eta'}_{i,j}^{LM}(x,y)^*
= {\ovl\eta'}_{i,j}^{ML}((x,y)^*).$

Set  $\ovl\cR=\ovl\cR(\TY(D,1,n))$:
$\ovl\cR(x,y)=(x',y')$
where
\begin{eqnarray*}
&&\hspace{-20pt}
x'_i\seteq y_i\frac{\ovl V_i}
{\ovl V_0 },\q
\ovl x'_i\seteq \ovl y_i\frac{\ovl V_0^\sharp 
\ovl W_i}{\ovl V_i\ovl W_1},\,
(1\leq i\leq n-2)\q
x'_{n-1}\seteq y_{n-1}\frac{{\ovl V^*_{n-1}}}
{\ovl V_{0}},\q
x'_n\seteq y_n\frac{\ovl V_{n-1}}
{\ovl V_0},\\
&&\hspace{-20pt}
%y'_1\seteq x_1\frac{{\ovl V^*_1}}
%{{\ovl V_0}^\sharp},\,\,
y'_i\seteq x_i\frac{{\ovl V^*_i}\ovl W_1}
{\ovl V_0^{\sharp}\ovl W_i},\q
\ovl y'_i\seteq \ovl x_i\frac{\ovl V_0}
{\ovl V_i^*}\,\,\,(1\leq i\leq n-2),\,\,
y'_{n-1}\seteq x_{n-1}\frac{\ovl V_{n-1}
\ovl W_1}{\ovl V_0^\sharp\ovl W_{n-1}},
\,\,
y'_n
\seteq x_n\frac{{\ovl V^*_{n-1}}
\ovl W_1}
{\ovl V_0^\sharp \ovl W_{n-1}}.
\end{eqnarray*}

%%%%%%%%%%%%%%%%%%%%%%%%%%%%%%%%%%%%%%%
%%%%%%%%%%%%%%%%%%%%%%%%%%%%%%%%%%%%%%
\subsection{\Tropical R for $\TY(B,1,n)$
$(n\geq2)$}%8.3
\label{trop-b1n}

By applying the method of folding, we shall
obtain an explicit form of the \tropical R for 
$\{\cB_L(\TY(B,1,n))\}_{L\in \bbC^\times}$
and $\{\cV(\TY(B,1,n))_L\}_{L\in\bbC^\times}$.
Indeed, it follows from Theorem \ref{rxx} that 
we have the birational map 
$R(\TY(B,1,n))=R_{LK}(\TY(B,1,n))
:\TY(X,{n,L},1)\times \TY(X,{n,K},1)
\to\TY(X,{n,K},1)\times \TY(X,{n,L},1)$.
Let us see that $R_{LK}(\TY(B,1,n))$ is 
a \tropical R.
\begin{lem}
\label{rx1}
The birational map $R(\TY(B,1,n))$ is 
a \tropical R on the 
$\TY(B,1,n)$-geometric crystals
$\{\TY(X,{n,L},1)\}_{L}$.
\end{lem}
{\sl Proof.}
It suffices to show (\ref{R-e})-(\ref{inver}) 
in Definition \ref{trop-r}. 
The relations  (\ref{YB}) and (\ref{inver}) 
for $\cR(\TY(B,1,n))$ are obtained 
from the one 
for $\cR(\TY(D,1,n+1))$. 
The others are easily shown
since $\cR(\TY(D,1,n+1))$
commutes with the actions of $e_i^c$ and 
preserves $\gamma_i$ and $\vep_i$ $(i\in I)$, 
and the data 
$(e_i^{\TY(B,1,n)})^c$, $\gamma_i^{\TY(B,1,n)}$
and $\vep_i^{\TY(B,1,n)}$ on $\TY(X,{(n,L)},1)$
are defined as in Proposition \ref{fix-b}.
\qed

Now let us describe the explicit form of tropical 
R on $\cB_L(\TY(B,1,n))$.
Set $\cR(\TY(B,1,n))\seteq (\eta^{-1},\eta^{-1})
\circ R(\TY(B,1,n))\circ(\eta,\eta)$ where
$\eta$ is as in \ref{subsec-bn}.
Let $\sharp:\cB_L(\TY(B,1,n))\rightarrow 
\cB_L(\TY(B,1,n))$
be an involution defined by 
\begin{equation}
\sharp(l_1,l_2,\ld,l_n,\ovl l_{n},\ld,\ovl l_2,\ovl l_1)
=(\ovl l_1,l_2,\ld,l_n,\ovl l_{n},\ld,\ovl l_2,l_1),
\end{equation}
that is, $\sharp:l_1\leftrightarrow \ovl l_1$.
Let $*:\cB_L(\TY(B,1,n))\times 
\cB_M(\TY(B,1,n))
\rightarrow \cB_M(\TY(B,1,n))\times 
\cB_L(\TY(B,1,n))$ 
be an involution defined by 
\[
((l_1,l_2,\ld,\ovl l_2,\ovl l_1),
(m_1,m_2,\ld,\ovl m_2,\ovl m_1))^*
=((\ovl m_1,\ovl m_2,\ld,m_2,m_1),
(\ovl l_1,\ovl l_2,\ld, l_2,l_1))
\]
that is, 
$*:l_i\leftrightarrow \ovl m_i, \q
\ovl l_i\leftrightarrow m_i\,\,
(1\leq i\leq n).$

Restricting the functions $\TTY(V,n+1,i)$
and $\TTY(W,n+1,i)$ for $\TY(D,1,n+1)$ to 
$l_{n+1}=m_{n+1}=1$, 
we define the rational functions 
$V_i$ $(i=0,1,\ld,n)$ and $W_i$ 
($i=1,\ld,n$) on 
$\cB_L(\TY(B,1,n))\times
\cB_M(\TY(B,1,n))$ $(L,M\in\bbC^\times)$ as 
\[
 W_i:=V_i{V^*_i}
+(M-L){V^*_i}+(L-M)V_i
\q(1\leq i\leq n-1),\q
W_{n}:=V_{n}{V^*_{n}}.
\]
\begin{equation}
V_i=
\sum_{j=1}^{n-1}(\theta_{i,j}(l,m)
+\theta'_{i,j}(l,m))
+\sum_{j=1}^{n+1}(\eta_{i,j}(l,m)
+\eta'_{i,j}(l,m)),
\end{equation}
where $L=l_1l_2\cd \ovl l_2\ovl l_1$,
$M=m_1m_2\cd \ovl m_2\ovl m_1$,
\begin{eqnarray}
&&\theta_{i,j}(l,m)
=\begin{cases}
\displaystyle L\prod_{k=j+1}^i
\frac{\ovl m_k}{\ovl l_k}&\text{for }
1\leq j\leq i,\\
\displaystyle M\prod_{k=i+1}^j
\frac{\ovl l_k}{\ovl m_k}&\text{for }
i+1\leq j\leq n-1,
\end{cases}\\
&&\theta'_{i,j}(l,m)=L\left(\prod_{k=1}^i
\frac{\ovl m_k}{\ovl l_k}\right)
\left(\prod_{k=1}^j\frac{m_k}{l_k}\right)
\q\text{for}\q j=1,\ld,n-1,
\end{eqnarray}
\begin{eqnarray}
&&\eta_{i,j}(l,m)=
\begin{cases}
\displaystyle L\left(\prod_{k=j+1}^i
\frac{\ovl m_k}{\ovl l_k}\right)
\left(\frac{\ovl m_j}{l_j}\right)&\text{for}\q 1\leq j\leq i,\\
\displaystyle M\left(\prod_{k=i+1}^j\frac{\ovl l_k}{\ovl m_k}\right)
\left(\frac{\ovl m_j}{l_j}\right)&\text{for}\q i+1\leq j\leq n,\\
\displaystyle M
\left(\prod_{k=i+1}^{n}
\frac{\ovl l_k}{\ovl m_k}\right)
&\text{for}\q j=n+1,
\end{cases}\\
&&\eta'_{i,j}(l,m)=
\begin{cases}
\displaystyle L\left(\prod_{k=1}^i\frac{\ovl m_k}{\ovl l_k}\right)
\left(\prod_{k=1}^j\frac{m_k}{l_k}\right)
\left(\frac{l_j}{\ovl m_j}\right)
&\text{for}\q 1\leq j\leq n,\\
\displaystyle L
\left(\frac{L}{M}\right)^{\del_{i,n}}
\left(\prod_{k=1}^i\frac{\ovl m_k}{\ovl l_k}\right)
\left(\prod_{k=1}^{n}\frac{m_k}{l_k}\right)
&\text{for}\q j=n+1.
\end{cases}
\end{eqnarray}
Now, we define the \tropical R map 
$\cR(\TY(B,1,n))$ on $\cB_L(\TY(B,1,n))
\times \cB_M(\TY(B,1,n))$ by
\begin{eqnarray}
&&{\cR(\TY(B,1,n))}(l,m)=(l',m') 
\q {\rm where}
\nn\\
 &&l'_1=
m_1\frac{V_0^{\sharp}}{V_1}, \q 
\ovl l'_1=\ovl m_1\frac{V_0}{V_1},\q
l'_i=m_i\frac{V_{i-1}W_i}
{V_iW_{i-1}},
\q \ovl l'_i=\ovl m_i\frac{V_{i-1}}{V_i}
\q(2\leq i\leq n),\label{explicit-R-b}\\
&&
m'_1=l_1\frac{V_0}{{V^*_1}}, \q 
\ovl m'_1=\ovl l_1
\frac{V_0^\sharp}{{V^*_1}},\q
m'_i=l_i\frac{{V^*_{i-1}}}
{{V^*_i}},
\q \ovl m'_i=\ovl l_i
\frac{{V^*_{i-1}}W_i}
{{V^*_i}W_{i-1}}
\q(2\leq i\leq n).
\nn
\end{eqnarray}
Here note that for $(l',m')={\mathcal R}(l,m)$ we have 
$l_1'l_2'\cd \ovl l_2'\ovl l_1'
=M$ and $m_1'm_2'\cd\ovl m_2'\ovl m_1'
=L$.

We shall describe the \tropical R map on 
$\cV(\TY(B,1,n))_L\times \cV(\TY(B,1,n))_M$ defined
by $\ovl\cR(x,y)\seteq
(\Xi,\Xi)\circ\cR(\TY(B,1,n))
\circ(\Xi^{-1},\Xi^{-1})(x,y)$
($x\in\cV(\TY(B,1,n))_L,\,
y\in\cV(\TY(B,1,n))_M$).
Restricting the rational functions
$\TTY(\ovl V,n+1,i)$ and $\TTY(\ovl W,n+1,i)$ 
for $\TY(D,1,n+1)$
to $x_n=x_{n+1}$ and $y_n=y_{n+1}$, 
we define the rational functions
$\ovl V_i$ $(i=0,1,\ld,n)$ and 
$\ovl W_i$ ($i=1,\ld,n$) on 
$\cV(\TY(B,1,n))_L\times
\cV(\TY(B,1,n))_M$ $(L,M\in\bbC^\times)$ by 
\[
\ovl W_i:=\ovl V_i{\ovl V^*_i}
+(M-L){\ovl V^*_i}+(L-M)\ovl V_i
\q(1\leq i\leq n-1),\q
\ovl W_{n}:=\ovl V_{n}
{\ovl V^*_{n}}.
\] 
\begin{equation}
\ovl V_i=
\sum_{j=1}^{n-1}(\ovl\theta_{i,j}(x,y)
+\ovl\theta'_{i,j}(x,y))
+\sum_{j=1}^{n+1}(\ovl\eta_{i,j}(x,y)
+\ovl\eta'_{i,j}(x,y)),
\end{equation}
where $\sharp$ is the involution on 
$\cV(\TY(B,1,n))_L$ and $*$ is the involution 
$\cV(\TY(B,1,n))_L\times \cV(\TY(B,1,n))_M
\to \cV(\TY(B,1,n))_M\times \cV(\TY(B,1,n))_L$
defined by
\begin{eqnarray*}
\sharp&:&x_i\mapsto \frac{x_i}{Lx_1\ovl x_1}\q
(i=1,\cd, n),\q
\ovl x_i\mapsto \frac{\ovl x_i}{Lx_1\ovl x_1}\q
(i=1,\cd,n-1), \\
*&:&x_i\mapsto \frac{1}{M\ovl y_i},\q
\ovl x_i\mapsto \frac{1}{M y_i}\q(i=1,\cd,n-1),\q
x_n\mapsto \frac{1}{M y_n},\\
&&y_i\mapsto \frac{1}{L\ovl x_i},\q
\ovl y_i\mapsto \frac{1}{L x_i}\q(i=1,\cd,n-1),\q
y_n\mapsto \frac{1}{L x_n},
\end{eqnarray*}
and
\begin{eqnarray*}
&&\ovl\theta_{i,j}(x,y)\seteq
L\frac{x_iy_j}{x_jy_i}\q(j\leq i\leq n,j<n),\q
\ovl\theta_{i,j}(x,y)\seteq
M\frac{x_iy_j}{x_jy_i}\q(i< j<n),\\
&&\ovl\theta'_{i,j}(x,y)\seteq
M\frac{x_i\ovl y_j}{\ovl x_jy_i}\q
(1\leq j<n),\q
\ovl\eta_{i,j}(x,y)\seteq
L^{1-\del_{j,1}}
\frac{x_i\ovl x_{j-1}y_{j-1}}{x_j\ovl x_jy_i}
\q(1\leq j\leq i\leq n,j<n),\\
&&\ovl\eta_{i,j}(x,y)\seteq
\frac{M}{L^{\del_{j,1}}}
\frac{x_i\ovl x_{j-1}y_{j-1}}{x_j\ovl x_jy_i}\q
(0\leq i<j<n),\q
\ovl\eta_{i,n}(x,y)=
M\frac{x_i\ovl x_{n-1}y_{n-1}}
{x_n^2y_i}\q(0\leq i\leq n),\\
&&
\ovl\eta_{i,n+1}(x,y)\seteq
M\frac{x_iy_{n}}
{x_{n}y_i}\q(i\leq n),\\
&&\ovl\eta'_{i,j}(x,y)\seteq
L^{\del_{j,1}}M
\frac{x_iy_j\ovl y_j}{\ovl x_{j-1}y_{j-1}y_i}\q
(1\leq j<n,\,i\leq n),\q
\ovl\eta'_{n,n+1}(x,y)\seteq L,\\
&&\ovl\eta'_{i,n}(x,y)\seteq
M\frac{x_iy_n^2}{\ovl x_{n-1}y_{n-1}y_i}\q
(0\leq i\leq n),\q
\ovl\eta'_{i,n+1}(x,y)\seteq
M\frac{x_iy_{n}}
{y_ix_{n}}\q
(0\leq i<n),
\end{eqnarray*}
where we understand $x_0=\ovl x_0=y_0=\ovl y_0=1$.
Note that as above {\it e.g.,}
\[
\ovl\theta_{i,j}(x,y)^*:=
\ovl\theta_{i,j}^{LM}(x,y)^*
=\ovl\theta_{i,j}^{ML}((x,y)^*).
\]
Here we define $\ovl\cR(x,y)=
\ovl\cR(\TY(B,1,n))(x,y)=(x',y')$
where
\begin{eqnarray*}
&&\hspace{-20pt}
x'_i\seteq y_i\frac{\ovl V_i}{\ovl V_0 }\,\,
(1\leq i\leq n-1)\q
x'_{n}\seteq y_{n}\frac{{\ovl V_{n}}}
{\ovl V_{0}},\q
\ovl x'_1\seteq \ovl y_1
\frac{\ovl V_0^\sharp}
{\ovl V_1},\q
\ovl x'_i\seteq \ovl y_i\frac{\ovl V_0^\sharp 
\ovl W_i}{\ovl V_i\ovl W_1}\,\,\,
(i\geq2),\\
&&\hspace{-20pt}
y'_1\seteq x_1\frac{{\ovl V^*_1}}
{\ovl V_0^\sharp},\,\,
y'_i\seteq x_i\frac{{\ovl V^*_i}\ovl W_1}
{\ovl V_0^{\sharp}\ovl W_i}\,\,\,(i>1),
\,\,\,
\ovl y'_i\seteq \ovl x_i\frac{\ovl V_0}
{{\ovl V^*_i}}\,\,\,
(1\leq i\leq n-1),\q
y'_{n}\seteq x_{n}\frac{\ovl V_{n}
\ovl W_1}{\ovl V_0^\sharp \ovl W_{n}}.
\end{eqnarray*}

%%%%%%%%%%%%%%%%%%%%%%%%%%%%%%%%%%%%%%%%
\subsection{\Tropical R for $\TY(D,2,n+1)$
$(n\geq2)$}

As in the previous subsection, 
we shall describe 
\tropical R maps of type $\TY(D,2,n+1)$.
We see the following lemma for 
$\{\TY(X,{n,L},3)\}_{L\in\bbC^\times}$.
\begin{lem}
\label{rx3}
The birational map $R(\TY(D,2,n+1))$ is 
a \tropical R map on 
%the $\TY(D,2,n+1)$-geometric crystals
$\{\TY(X,{n,L},3)\}_{L}$.
\end{lem}
The proof is the 
same as the one for Lemma \ref{rx1}.

Let us describe the explicit form of tropical 
R on $\cB_L(\TY(D,2,n+1))$ as in the previous
subsection.
Set $\cR(\TY(D,2,n+1))\seteq (\eta^{-1},\eta^{-1})
\circ R(\TY(D,2,n+1))\circ(\eta,\eta)$ where
$\eta$ is as in \ref{subsec-d2n}.
Let $*:\cB_L(\TY(D,2,n+1))\times 
\cB_M(\TY(D,2,n+1))
\rightarrow \cB_M(\TY(D,2,n+1))\times 
\cB_L(\TY(D,2,n+1))$ 
be an involution defined by 
\[
((l_0,l_1,\ld,\ovl l_2,\ovl l_1),
(m_0,m_1,\ld,\ovl m_2,\ovl m_1))^*
=((m_0,\ovl m_1,\ovl m_2,\ld,m_2,m_1),
(l_0,\ovl l_1,\ovl l_2,\ld, l_2,l_1))
\]
that is, $*:l_0\leftrightarrow m_0$ and 
$l_i\leftrightarrow \ovl m_i, \q
\ovl l_i\leftrightarrow m_i\,\,
(1\leq i\leq n)$.

Restricting the functions $\TTY(V,n+1,i)$
and $\TTY(W,n+1,i)$ for $\TY(B,1,n+1)$ to 
$l_1=\ovl l_1$ and $m_1=\ovl m_1$, and replacing
$X_i$ with $X_{i-1}$ 
where $X=l,\ovl l,m,$ and $\ovl m$, 
we define the rational functions 
$V_i$ $(i=0,1,\ld,n+1)$ and $W_i$ 
($i=1,\ld,n+1$) on 
$\cB_L(\TY(D,2,n+1))\times
\cB_M(\TY(D,2,n+1))$ $(L,M\in\bbC^\times)$ as 
\[
 W_i:=V_i{V^*_i}
+(M-L){V^*_i}+(L-M)V_i
\q(1\leq i\leq n),\q
W_{n+1}:=V_{n+1}{V^*_{n+1}}.
\]
\begin{equation}
V_i=
\sum_{j=1}^{n}(\theta_{i,j}(l,m)
+\theta'_{i,j}(l,m))
+\sum_{j=1}^{n+2}(\eta_{i,j}(l,m)
+\eta'_{i,j}(l,m)),
\end{equation}
where $L=l_0^2l_1l_2\cd \ovl l_2\ovl l_1$,
$M=m_0^2m_1m_2\cd \ovl m_2\ovl m_1$,
\begin{eqnarray}
&&\theta_{i,j}(l,m)
=\begin{cases}
\displaystyle L\prod_{k=j}^{i-1}
\frac{\ovl m_k}{\ovl l_k}&\text{for }
1\leq j\leq i\leq n+1,\\
\displaystyle M\prod_{k=i}^{j-1}
\frac{\ovl l_k}{\ovl m_k}&\text{for }
0\leq i<j\leq n,
\end{cases}\\
&&\theta'_{i,j}(l,m)=L\left(\prod_{k=0}^{i-1}
\frac{\ovl m_k}{\ovl l_k}\right)
\left(\prod_{k=0}^{j-1}\frac{m_k}{l_k}\right)
\q\text{for}\q j=1,\ld,n,
\end{eqnarray}
\begin{eqnarray}
&&\eta_{i,j}(l,m)=
\begin{cases}
\displaystyle L\left(\prod_{k=j}^{i-1}
\frac{\ovl m_k}{\ovl l_k}\right)
\left(\frac{\ovl m_{j-1}}{l_{j-1}}\right)
&\text{for}\q 1\leq j\leq i\leq n+1,\\
\displaystyle M\left(\prod_{k=i}^{j-1}
\frac{\ovl l_k}{\ovl m_k}\right)
\left(\frac{\ovl m_{j-1}}{l_{j-1}}\right)
&\text{for}\q i+1\leq j\leq n+1,\\
\displaystyle M
\left(\prod_{k=i}^{n}
\frac{\ovl l_k}{\ovl m_k}\right)
&\text{for}\q j=n+2,
\end{cases}\\
&&\eta'_{i,j}(l,m)=
\begin{cases}
\displaystyle L\left(\prod_{k=0}^{i-1}
\frac{\ovl m_k}{\ovl l_k}\right)
\left(\prod_{k=0}^{j-1}\frac{m_k}{l_k}\right)
\left(\frac{l_{j-1}}{\ovl m_{j-1}}\right)
&\text{for}\q 1\leq j\leq n+1,\\
\displaystyle L
\left(\frac{L}{M}\right)^{\del_{i,n+1}}
\left(\prod_{k=0}^{i-1}
\frac{\ovl m_k}{\ovl l_k}\right)
\left(\prod_{k=0}^{n}\frac{m_k}{l_k}\right)
&\text{for}\q j=n+2,
\end{cases}
\end{eqnarray}
where we understand $\ovl l_0=l_0,\ovl m_0=m_0$.
Now, we define the \tropical R map 
$\cR(\TY(D,2,n+1))$ on $\cB_L(\TY(D,2,n+1))
\times \cB_M(\TY(D,2,n+1))$ by
\begin{eqnarray}
&&{\cR(\TY(D,2,n+1))}(l,m)=(l',m') 
\q {\rm where}
\nn\\
 &&l'_0=
m_0\frac{{V_0}}{V_1}, \q 
l'_i=m_i\frac{V_{i}W_{i+1}}
{V_{i+1}W_i},
\q \ovl l'_i
=\ovl m_i\frac{V_i}{V_{i+1}}
\q(1\leq i\leq n),\label{explicit-R-d2}\\
&&
m'_0=l_0\frac{V_0}{{V^*_1}}, \q 
m'_i=l_i\frac{{V^*_{i}}}
{{V^*_{i+1}}},
\q \ovl m'_i=\ovl l_i
\frac{{V^*_{i}}W_{i+1}}
{{V^*_{i+1}}W_{i}}
\q(1\leq i\leq n).
\nn
\end{eqnarray}
Here note that for 
$(l',m')={\mathcal R}(\TY(D,2,n+1))(l,m)$ we have 
${l'_0}^2l_1'l_2'\cd \ovl l_2'\ovl l_1'
=M$ and ${m'_0}^2m_1'm_2'\cd\ovl m_2'\ovl m_1'
=L$.

Next, we shall describe the \tropical R on 
$\cV(\TY(D,2,n+1))_L\times \cV(\TY(D,2,n+1))_M$
defined by
$\ovl\cR(x,y)\seteq
(\Xi,\Xi)\circ\cR(\TY(D,2,n+1))
\circ(\Xi^{-1},\Xi^{-1})(x,y)$
($x\in\cV(\TY(D,2,n+1))_L,\,
y\in\cV(\TY(D,2,n+1))_M$).

In this case, we replace $L$ (resp. $M$)
for $\cB_L(\TY(D,2,n+1))$
(resp. $\cB_M(\TY(D,2,n+1))$)
with $L^2$ (resp. $M^2$) since 
$\cV(\TY(D,2,n+1))_L\cong \cB_{L^2}(\TY(D,2,n+1))$
as in \ref{subsec-d2n}.

%Restricting the rational functions
%$\TTY(\ovl V,n+1,i)$ and $\TTY(\ovl W,n+1,i)$ 
%for $\TY(B,1,n+1)$ 
%to $x_1=L^2\ovl x_1$ and $y_1=M^2y_1$ and 
%replacing $X_i$ with $X_{i-1}$ 
%($X=x,y,\ovl x,\ovl y$),
We define the rational functions
$\ovl V_i(x,y)$ $(i=0,1,\ld,n+1)$ and 
$\ovl W_i(x,y)$ ($i=1,\ld,n+1$) on 
$\cV(\TY(D,2,n+1))_L\times
\cV(\TY(D,2,n+1))_M$ $(L,M\in\bbC^\times)$ by 
\[
\ovl W_i:=\ovl V_i{\ovl V^*_i}
+(M^2-L^2){\ovl V^*_i}
+(L^2-M^2)\ovl V_i
\q(1\leq i\leq n),\q
\ovl W_{n+1}:=\ovl V_{n+1}
{\ovl V^*_{n+1}}.
\] 
\begin{equation}
\ovl V_i=
\sum_{j=1}^{n}(\ovl\theta_{i,j}(x,y)
+\ovl\theta'_{i,j}(x,y))
+\sum_{j=1}^{n+2}(\ovl\eta_{i,j}(x,y)
+\ovl\eta'_{i,j}(x,y)),
\end{equation}
where $*$ is the involution 
$\cV(\TY(D,2,n+1))_L\times \cV(\TY(D,2,n+1))_M
\to \cV(\TY(D,2,n+1))_M\times \cV(\TY(D,2,n+1))_L$
defined by
\begin{eqnarray*}
*&:&x_0\mapsto y_0,\,\,
x_i\mapsto \frac{y_0^2}{M^2 \ovl y_i},\,\,
\ovl x_i\mapsto \frac{y_0^2}{M^2 y_i}\,\,
(i=1,\cd,n-1),\,\,
x_n\mapsto \frac{y_0^2}{M^2 y_n},\\
&&y_0\mapsto x_0,\,\,
y_i\mapsto \frac{x_0^2}{L^2 \ovl x_i},\,\,
\ovl y_i\mapsto \frac{x_0^2}{L^2 x_i}\,\,
(i=1,\cd,n-1),\,\,
y_n\mapsto \frac{x_0^2}{L^2 x_n}.
\end{eqnarray*}
and
\begin{eqnarray*}
&&\ovl\theta_{i,j}(x,y)\seteq
\begin{cases}
L^2\left(\frac{x_{i-1}}{y_{i-1}}
\right)^{1+\del_{i,1}}
\left(\frac{y_{j-1}}{x_{j-1}}
\right)^{1+\del_{j,1}}
&(1\leq j\leq i\leq n+1,\,j\leq n),\\
M^2\left(\frac{x_0}{y_0}
\right)^{\del_{i,0}+\del_{i,1}-\del_{j,1}}
\frac{x_{i-1}y_{j-1}}{x_{j-1}y_{i-1}}
&(0\leq i<j\leq n),
\end{cases}\\
&&\ovl\theta'_{i,j}(x,y)\seteq
M^2\left(\frac{L^2}{M^2}\right)^{\del_{j,1}}
\left(\frac{x_0}{y_0}
\right)^{\del_{i,0}+\del_{i,1}}
\frac{x_{i-1}\ovl y_{j-1}}{\ovl x_{j-1}y_{i-1}}\q
(1\leq j\leq n),\\
&&\ovl\eta_{i,j}(x,y)\seteq
\begin{cases}
L^{2(1-\del_{j,2})}y_0^{\del_{j,1}+\del_{j,2}}
\left(\frac{x_0}{y_0}\right)^{\del_{i,1}\del_{j,1}}
\frac{x_{i-1}\ovl x_{j-2}y_{j-2}}
{x_{j-1}\ovl x_{j-1}y_{i-1}}&
(1\leq j\leq i\leq n+1),\\
\frac{M^2}{L^{2\del_{j,2}}}
\left(\frac{x_0}{y_0^{1-\del_{j,1}-\del_{j,2}}}
\right)^{\del_{i,0}+\del_{i,1}}
\frac{x_{i-1}\ovl x_{j-2}y_{j-2}}
{x_{j-1}\ovl x_{j-1}y_{i-1}}&
(i+1\leq j\leq n+1),\\
M^2\left(\frac{x_0}{y_0}
\right)^{\del_{i,0}+\del_{i,1}}
\frac{x_{i-1}y_{n}}
{x_ny_{i-1}}&(j=n+2),
\end{cases}\\
&&
\ovl\eta'_{i,j}(x,y)\seteq
\begin{cases}
\ovl\theta_{i,j}'(x,y)\times
\left(\frac{L^2}{y_0}\right)^{\del_{j,2}}
\frac{\ovl x_{j-1}y_{j-1}}
{\ovl x_{j-2}y_{j-2}}&
(1\leq j\leq n),\\
M^2\left(\frac{x_0}{y_0}
\right)^{\del_{i,0}+\del_{i,1}}
\frac{x_{i-1}y_n^2}{y_{i-1}\ovl x_{n-1}y_{n-1}}&
(j=n+1),\\
M^2\left(\frac{L^2}{M^2}\right)^{\del_{i,n+1}}
\left(\frac{x_0}{y_0}\right)^{\del_{i,0}+\del_{i,1}}
\frac{x_{i-1}y_n}{y_{i-1}x_n}&(j=n+2).
\end{cases}
\end{eqnarray*}
where we understand 
$x_{-1}=\ovl x_0=y_{-1}=\ovl y_0=1$ 
and $\ovl x_n=x_n$, $\ovl y_n=y_n$.

Here we define $\ovl\cR(x,y)=
\ovl\cR(\TY(D,2,n+1))(x,y)=(x',y')$
by
\begin{eqnarray*}
&&\hspace{-20pt}
x'_0\seteq y_0\frac{\ovl V_1}
{\ovl V_0 }\,\,
x'_i\seteq y_i\frac{\ovl V_1
\TTY(\ovl V,n,i+1)}{{\TTY(\ovl V,n,0)}^2 }\,\,
(1\leq i\leq n)\q
\ovl x'_i\seteq \ovl y_i
\frac{\ovl V_1\TTY(\ovl W,n,i+1)}
{\ovl V_{i+1}\ovl W_1}\,\,\,
(1\leq i\leq n-1),\\
&&\hspace{-20pt}
y'_0\seteq x_0\frac{{\ovl V^*_1}}
{\ovl V_0},\,\,
y'_i\seteq x_i\frac{{\ovl V^*_1}
{\ovl V^*_{i+1}}\ovl W_1}
{{\ovl V_0}^{2}\ovl W_{i+1}}\,\,\,
(1\leq i\leq n),
\,\,\,
\ovl y'_i\seteq \ovl x_i\frac{{\ovl V^*_1}}
{{\ovl V^*_{i+1}}}\,\,\,
(1\leq i\leq n-1).
\end{eqnarray*}
%%%%%%%%%%%%%%%%%%%%%%%%%%%%%%
\subsection{\Tropical R for $\TY(A,2,2n-1)$
$(n\geq3)$}

We shall describe 
\tropical R's of type $\TY(A,2,2n-1)$.
We see the following lemma for 
$\{\TY(X,{n,L},2)\}_{L\in\bbC^\times}$.
\begin{lem}
\label{rx2}
The birational map $R(\TY(A,2,2n-1))$ is 
a \tropical R map on 
%the $\TY(A,2,2n-1)$-geometric crystals
$\{\TY(X,{n,L},2)\}_{L}$.
\end{lem}
The proof is the 
same as the one for Lemma \ref{rx1}.

Let us describe the explicit form of tropical 
R on $\cB_L(\TY(A,2,2n-1))$.
Set $\cR(\TY(A,2,2n-1))\seteq (\eta^{-1},\eta^{-1})
\circ R(\TY(A,2,2n-1))\circ(\eta,\eta)$ where
$\eta$ is as in \ref{subsec-a2o}.

Let $\sharp$ be the involution on 
$\cB_L(\TY(A,2,2n-1))$ defined by
$\sharp:l_1\leftrightarrow \ovl l_1$ for 
$l=(l_1,\cd,\ovl l_1)\in\cB_L(\TY(A,2,2n-1))$ and 
$*:\cB_L(\TY(A,2,2n-1))\times 
\cB_M(\TY(A,2,2n-1))
\rightarrow \cB_M(\TY(A,2,2n-1))\times 
\cB_L(\TY(A,2,2n-1))$ 
an involution defined by 
\[
((l_1,l_2,\ld,\ovl l_2,\ovl l_1),
(m_1,m_2,\ld,\ovl m_2,\ovl m_1))^*
=((\ovl m_1,\ovl m_2,\ld,m_2,m_1),
(\ovl l_1,\ovl l_2,\ld, l_2,l_1))
\]
that is, 
$*:l_i\leftrightarrow \ovl m_i, \q
\ovl l_i\leftrightarrow m_i\,\,
(1\leq i\leq n)$.

Restricting the functions $\TTY(V,2n,i)$
and $\TTY(W,2n,i)$ for $\TY(D,1,2n)$ to 
$\TY(X,{n,L^2},2)\times \TY(X,{n,M^2},2)$, 
we define the rational functions 
$V_i$ $(i=0,1,\ld,n)$ and $W_i$ 
($i=1,\ld,n$) on 
$\cB_L(\TY(A,2,2n-1))\times
\cB_M(\TY(A,2,2n-1))$ $(L,M\in\bbC^\times)$ as 
\[
 W_i:=V_i{V^*_i}
+(M^2-L^2){V^*_i}+(L^2-M^2)V_i
\q(1\leq i\leq n).
\]
\begin{equation}
V_i=
\sum_{j=1}^{2n-2}(\theta_{i,j}(l,m)
+\theta'_{i,j}(l,m))
+\sum_{j=1}^{2n}(\eta_{i,j}(l,m)
+\eta'_{i,j}(l,m)),
\end{equation}
where 
$L^2=l_n\ovl l_n(\prod_{i=1}^{n-1}l_i\ovl l_i)^2$,
$M^2=m_n\ovl m_n(\prod_{i=1}^{n-1}m_i\ovl m_i)^2$,
$\Del\seteq\frac{1+\ovl \mu(l)}{1+\ovl \mu(m)}$, and 
\begin{eqnarray*}
&&\hspace{-30pt}\theta_{i,j}(l,m)\\
&&\hspace{-30pt}
=\begin{cases}
\displaystyle L^2\prod_{k=j+1}^{i}
\frac{\ovl m_k}{\ovl l_k}
\Del^{\del_{i,n}-\del_{j,n}}&
1\leq j\leq i\leq n,\\
\displaystyle M^2\prod_{k=i+1}^{j}
\frac{\ovl l_k}{\ovl m_k}
\Del^{-\del_{j,n}}&
0\leq i<j\leq n,\nn \\
\displaystyle LM\left(\prod_{k=1}^i
\frac{\ovl m_k}{\ovl l_k}\right)
\left(\prod_{k=1}^{2n-j-1}\frac{m_k}{l_k}\right)
\Del^{\del_{i,n}}
\frac{\ovl l_{2n-j}(m_{2n-j}+\ovl m_{2n-j})}
{(l_{2n-j}+\ovl l_{2n-j})\ovl m_{2n-j}}
&0\leq i\leq n<j\leq 2n-2.
\end{cases}
\end{eqnarray*}
\[
%\begin{equation}
\hspace{-30pt}\theta'_{i,j}(l,m)
=\begin{cases}
\displaystyle L^2\left(\prod_{k=1}^{i}
\frac{\ovl m_k}{\ovl l_k}\right)
\left(\prod_{k=1}^j\frac{m_k}{l_k}\right)
\Del^{\del_{i,n}}
\left(\frac{1+\mu(l)}{1+ \mu(m)}\right)
^{\del_{j,n}}
&j\leq n,\\
\displaystyle LM
\left(\prod_{k=1}^i
\frac{\ovl m_k}{\ovl l_k}\right)
\left(\prod_{k=1}^{2n-j-1}\frac{\ovl l_k}
{\ovl m_k}\right)\Del^{\del_{i,n}}
\frac{m_{2n-j}(l_{2n-j}+\ovl l_{2n-j})}
{l_{2n-j}(m_{2n-j}+\ovl m_{2n-j})}&n<j.
\end{cases}
\]
%\end{equation}
\begin{eqnarray*}
&&\hspace{-20pt}\eta_{i,j}(l,m)=
\begin{cases}
\displaystyle L^2\left(\prod_{k=j+1}^{i}
\frac{\ovl m_k}{\ovl l_k}\right)
\left(\frac{\ovl m_{j}}{l_{j}}\right)
\Del^{\del_{i,n}}\mu(l)^{\del_{j,n}}
&1\leq j\leq i\leq n,\\
\displaystyle M^2\left(\prod_{k=i+1}^{j}
\frac{\ovl l_k}{\ovl m_k}\right)
\left(\frac{\ovl m_{j}}{l_{j}}\right)
\mu(l)^{\del_{j,n}}
&i<j\leq n,\\
\displaystyle LM
\left(\prod_{k=1}^{i}
\frac{\ovl m_k}{\ovl l_k}\right)
\left(\prod_{k=1}^{n-1}
\frac{m_k}{l_k}\right)\Del^{\del_{i,n}}
\frac{1+\mu(m)}{1+\ovl \mu(l)}
& i\leq n,\, j=n+1\\
\displaystyle LM
\left(\prod_{k=1}^{i}
\frac{\ovl m_k}{\ovl l_k}\right)
\left(\prod_{k=1}^{2n-j}\frac{m_{k}}{l_{k}}\right)
\Del^{\del_{i,n}}
\frac{(m_{2n-j+1}+\ovl m_{2n-j+1})l_{2n-j+1}}
{(l_{2n-j+1}+\ovl l_{2n-j+1})\ovl m_{2n-j+1}}
&i<n+1<j\leq  2n-1,\\
%\displaystyle LM
%\left(\prod_{k=1}^i
%\frac{\ovl m_k}{\ovl l_k} \right)
%\frac{m_1(\frac{m_2}{\ovl m_2}+1)}
%{\frac{l_1\ovl l_2}{l_2}(\frac{l_2}{\ovl l_2}+1)}
%\Del^{\del_{i,n}}
%&i<j=2n-1,\\
\displaystyle LM
\left(\prod_{k=1}^i
\frac{\ovl m_k}{\ovl l_k} \right)
\frac{\ovl l_1}{l_1}\Del^{\del_{i,n}}&i<j=2n.
\end{cases}
\end{eqnarray*}
\[
\eta'_{i,j}(l,m)=
\begin{cases}
\displaystyle L^2\left(\prod_{k=1}^{i}
\frac{\ovl m_k}{\ovl l_k}\right)
\left(\prod_{k=1}^{j}\frac{m_k}{l_k}\right)
\left(\frac{l_{j}}{\ovl m_{j}}\right)
\Del^{\del_{i,n}}\ovl \mu(m)^{\del{j,n}}
&1\leq j\leq n,\\
\displaystyle L^2\left(\prod_{k=1}^{i}
\frac{\ovl m_k}{\ovl l_k}\right)
\left(\prod_{k=1}^{n}\frac{m_k}{l_k}\right)
\Del^{\del_{i,n}}
\frac{1+\mu(l)}{(1+\mu(m))\mu(m)}
&j=n+1,\\
\displaystyle LM
\left(\prod_{k=1}^{i}
\frac{\ovl m_k}{\ovl l_k}\right)
\left(\prod_{k=1}^{2n-j}\frac{\ovl l_{k}}
{\ovl m_{k}}\right)
\Del^{\del_{i,n}}
\frac{(l_{2n-j+1}+\ovl l_{2n-j+1})\ovl m_{2n-j+1}}
{(m_{2n-j+1}+\ovl m_{2n-j+1})l_{2n-j+1}}
&n+1<j\leq 2n-1,\\
%\displaystyle LM
%\left(\prod_{k=1}^i
%\frac{\ovl m_k}{\ovl l_k}\right)
%\Del^{\del_{i,n}}
%\frac{\ovl l_1\ovl m_2(l_2+\ovl l_2)}
%{\ovl m_1l_2(m_2+\ovl m_2)}
%&j=2n-1,\\
\displaystyle LM
\left(\prod_{k=1}^i
\frac{\ovl m_k}{\ovl l_k}\right)
\Del^{\del_{i,n}}
\frac{m_1}{\ovl m_1}
&j=2n.
\end{cases}
\]
Now, we define the \tropical R 
$\cR(\TY(A,2,2n-1))\cl\cB_L(\TY(A,2,2n-1))\times 
\cB_M(\TY(A,2,2n-1))\to \cB_M(\TY(A,2,2n-1))
\times \cB_L(\TY(A,2,2n-1))$ by 
\begin{eqnarray}
&&{\cR(\TY(A,2,2n-1))}(l,m)=(l',m') 
\q {\rm where}
\nn\\
 &&l'_1=
m_1\frac{V_0^\sharp}{V_1}, \q 
%\ovl l'_1=
%\ovl m_1\frac{{V_0}}{V_1}, \q 
l'_i=m_i\frac{V_{i-1}W_{i}}
{V_{i}W_{i-1}},\q
\ovl l'_i
=\ovl m_i\frac{V_{i-1}}{V_i}
\q(1\leq i\leq n-1),\nn\\
&&
l'_n=\frac{m_n\TTY(V,n,n-1)\TTY(W,n,n)}
{(1+\mu(m))\TTY(V,n,n)\TTY(W,n,n-1)}
\left(1+\frac{m_n\TTY(W,n,n)}
{\ovl m_n\TTY(W,n,n-1)}\ovl \mu(m)\right),\q
\ovl l'_n=\frac{\ovl m_n\TTY(V,n,n-1)}
{(1+\ovl \mu(m))\TTY(V,n,n)}
\left(1+\frac{\ovl m_n\TTY(W,n,n-1)}
{m_n\TTY(W,n,n)}\mu(m)\right),
%\label{explicit-R-aon}
\nn\\
&&
m'_1=l_1\frac{V_0}{{V^*_1}}, \q 
\ovl m'_1=\ovl l_1
\frac{V_0^\sharp}{{V^*_1}}, \q 
m'_i=l_i\frac{{V^*_{i-1}}}
{{V^*_{i}}},
\q \ovl m'_i=\ovl l_i
\frac{{V^*_{i-1}}W_{i}^{(n)}}
{{V^*_{i}}W_{i-1}^{(n)}}
\q(2\leq i\leq n-1),\nn\\
&&
m'_n=\frac{l_n{\STY(V,*,n-1)}}
{(1+\mu(l)){\STY(V,*,n)}}
\left(1+\frac{l_n\TTY(W,n,n-1)}
{\ovl l_n\TTY(W,n,n)}\ovl \mu(l)\right),\q
\ovl m'_n=\frac{\ovl l_n{\STY(V,*,n-1)}\TTY(W,n,n)}
{(1+\ovl \mu(l)){\STY(V,*,n)}\TTY(W,n,n-1)}
\left(1+\frac{\ovl l_n\TTY(W,n,n)}
{l_n\TTY(W,n,n-1)}\mu(l)\right),
\nn
\end{eqnarray}
Here note that for 
$(l',m')={\mathcal R}(\TY(A,2,2n-1))(l,m)$ we have 
$(l_1'l_2'\cd l'_{n-1}\ovl l'_{n-1}\cd
 \ovl l_2'\ovl l_1')^2 l'_n\ovl l'_n
=M^2$ and \\ 
$(m_1'm_2'\cd m'_{n-1}\ovl m'_{n-1}\cd
\ovl m_2'\ovl m_1')^2m'_n\ovl m'_n
=L^2$.

Next, we shall describe \tropical R on 
$\cV(\TY(A,2,2n-1))_L\times \cV(\TY(A,2,2n-1))_M$.
Let $*$ be the involution 
$\cV(\TY(A,2,2n-1))_L\times \cV(\TY(A,2,2n-1))_M
\to \cV(\TY(A,2,2n-1))_M\times \cV(\TY(A,2,2n-1))_L$
defined by
\begin{eqnarray*}
*&:&
x_i\mapsto \frac{1}{M^2 \ovl y_i},\,\,
\ovl x_i\mapsto \frac{1}{M^2 y_i},\,\,
y_i\mapsto \frac{1}{L^2 \ovl x_i},\,\,
\ovl y_i\mapsto \frac{1}{L^2 x_i}\,\,
(i=1,\cd,n)
\end{eqnarray*}
and $\sharp$ the involution on 
$\cV(\TY(A,2,2n-1))_L$ defined by
\[
 \sharp:x_i\mapsto \frac{x_i}{L^2 x_1\ovl x_1},\q
\ovl x_i\mapsto \frac{\ovl x_i}{L^2x_1\ovl x_1}
\q(i=1,\cd,n-1),\q
x_n\mapsto \frac{x_n}{(L^2x_1\ovl x_1)^2}.
\]
We define the rational functions
$\ovl V_i(x,y)$ $(i=0,1,\ld,n)$ and 
$\ovl W_i(x,y)$ ($i=1,\ld,n$) on 
$\cV(\TY(A,2,2n-1))_L\times
\cV(\TY(A,2,2n-1))_M$ $(L,M\in\bbC^\times)$ by 
\begin{eqnarray*}
&&\ovl V_i=
\sum_{j=1}^{2n-2}(\ovl\theta_{i,j}(x,y)
+\ovl\theta'_{i,j}(x,y))
+\sum_{j=1}^{2n}(\ovl\eta_{i,j}(x,y)
+\ovl\eta'_{i,j}(x,y)),\\
&&\ovl W_i:=\ovl V_i{\ovl V^*_i}
+(M^2-L^2){\ovl V^*_i}
+(L^2-M^2)\ovl V_i
\q(1\leq i\leq n),\\
&&
\mu(x)\seteq\frac{x_n}{Lx_{n-1}\ovl x_{n-1}}=
\ovl\mu(x)^{-1},\,\,
\mu(y)\seteq\frac{y_n}{My_{n-1}\ovl y_{n-1}}=
\ovl\mu(y)^{-1},\\
&&
\Del\seteq\frac{1+\ovl\mu(x)}{1+\ovl\mu(y)},\q
\nabla\seteq\frac{1+\mu(x)}{1+\mu(y)},\q
\square\seteq\frac{1+\mu(y)}{1+\ovl\mu(x)}.
%\Del\seteq\frac{1+\frac{Lx_{n-1}\ovl x_{n-1}}
%{x_n}}{1+\frac{M y_{n-1}\ovl y_{n-1}}{y_n}},\q
%\nabla\seteq\frac{1+\frac{x_n}
%{Lx_{n-1}\ovl x_{n-1}}}
%{1+\frac{y_n}{M y_{n-1}\ovl y_{n-1}}},\q
%\square\seteq
%\frac{1+\frac{y_n}{M y_{n-1}\ovl y_{n-1}}}
%{1+\frac{Lx_{n-1}\ovl x_{n-1}}{x_n}}.
\end{eqnarray*}
\begin{eqnarray*}
&&\ovl\theta_{i,j}(x,y)\seteq
\begin{cases}
L^2\frac{x_iy_j}{x_jy_i}\left(\frac{y_{n-1}}
{x_{n-1}}\right)^{\del_{i,n}(1-\del_{i,j})}
\left(\Del\right)^{\del_{i,n}-\del_{j,n}}
&(1\leq j\leq i\leq n),\\
M^2\frac{x_iy_j}{x_jy_i}
\left(\frac{x_{n-1}}{y_{n-1}}\Del^{-1}\right)
^{\del_{j,n}}
&(0\leq i<j\leq n),\\
\frac{M^3}{L}\frac{x_i x_{2n-j-1}(y_{2n-j-1}
\ovl y_{2n-j-1}+y_{2n-j}\ovl y_{2n-j})}
{y_i y_{2n-j-1}(x_{2n-j-1}
\ovl x_{2n-j-1}+x_{2n-j}\ovl x_{2n-j})}
\left(\frac{y_{n-1}}{x_{n-1}}
\Del\right)^{\del_{i,n}}
&(i\leq n<j\leq 2n-2),
\end{cases}
\end{eqnarray*}
\begin{eqnarray*}
&&\ovl\theta'_{i,j}(x,y)\seteq
\begin{cases}
M^2\frac{x_i\ovl y_j}{y_i\ovl x_j}
\left(\frac{y_{n-1}}{x_{n-1}}\Del
\right)^{\del_{i,n}}
\left(\frac{L^2\ovl x_{n-1}}{M^2 \ovl y_{n-1}}
\nabla\right)^{\del_{j,n}}
&(1\leq j\leq n),\\
LM \frac{x_i y_{2n-j-1}
(1+\frac{x_{2n-j-1}\ovl x_{2n-j-1}}
{x_{2n-j}\ovl x_{2n-j}})}
{y_i x_{2n-j-1}
(1+\frac{y_{2n-j-1}\ovl y_{2n-j-1}}
{y_{2n-j}\ovl y_{2n-j}})}
\left(\frac{y_{n-1}}{x_{n-1}}\Del\right)^{\del_{i,n}}
&(n<j\leq 2n-2),
\end{cases}
\end{eqnarray*}
\begin{eqnarray*}
&&\ovl\eta_{i,j}(x,y)\seteq
\begin{cases}
L^{2-2\del_{j,1}}\frac{x_iy_{j-1}\ovl x_{j-1}}
{y_ix_j\ovl x_j}\left(\frac{y_{n-1}}{x_{n-1}}
\right)^{\del_{i,n}(1-\del_{i,j})}
\Del^{\del_{i,n}}\left(\frac{L x_ny_{n-1}}
{x_{n-1}}\right)^{\del_{j,n}}
&(1\leq j\leq i\leq n),\\
L^{-2\del_{j,1}}M^2\frac{x_iy_{j-1}\ovl x_{j-1}}
{y_ix_j\ovl x_j}(Lx_n)^{\del_{j,n}}
&(i< j\leq n),\\
\frac{M^3}{L}\frac{x_i\ovl y_{n-1}}
{y_i\ovl x_{n-1}}\left(
\frac{y_{n-1}}{x_{n-1}}\Del\right)^{\del_{i,n}}
\square
&(0\leq i\leq n,j=n+1),\\
\frac{M^3}{L}\frac{x_i\ovl y_{2n-j}}
{y_i\ovl x_{2n-j}}
\frac{1+\frac{y_{2n-j+1}
\ovl y_{2n-j+1}}{y_{2n-j}\ovl y_{2n-j}}}
{1+\frac{x_{2n-j}
\ovl x_{2n-j}}{x_{2n-j+1}\ovl x_{2n-j+1}}}
\left(\frac{y_{n-1}}{x_{n-1}}
\Del\right)^{\del_{i,n}}
&(0\leq i\leq n<j\leq 2n-1),\\
\frac{M}{L}
\frac{x_i}{y_ix_1\ovl x_1}\left(\frac{y_{n-1}}
{x_{n-1}}\Del\right)^{\del_{i,n}}
&(0\leq i\leq n<j=2n),
\end{cases}\\
&&
\ovl\eta'_{i,j}(x,y)\seteq
\begin{cases}
L^{2\del_{j,1}}M^2\frac{x_iy_{j}\ovl y_{j}}
{y_iy_{j-1}\ovl x_{j-1}}
\left(\frac{y_{n-1}}{x_{n-1}}\Del
\right)^{\del_{i,n}}\left(\frac{1}{{M y_n}}
\right)^{\del_{j,n}}
&(1\leq j\leq n),\\
L^2M\frac{x_iy_{n-1}\ovl x_{n-1}}
{y_ix_n}
\left(\frac{y_{n-1}}{x_{n-1}}\Del
\right)^{\del_{i,n}}\nabla
&(j=n+1),\\
LM\frac{x_i y_{2n-j}}{y_i x_{2n-j}}
\frac{\left(1+\frac{x_{2n-j}\ovl x_{2n-j}}
{x_{2n-j+1}\ovl x_{2n-j+1}}\right)}
{\left(1+\frac{y_{2n-j+1}\ovl y_{2n-j+1}}
{y_{2n-j}\ovl y_{2n-j}}\right)}
\left(\frac{y_{n-1}}{x_{n-1}}\Del
\right)^{\del_{i,n}}
&(n+1<j\leq 2n-1),\\
LM^3\frac{x_iy_1\ovl y_1}{y_i}
\left(\frac{y_{n-1}}{x_{n-1}}\Del
\right)^{\del_{i,n}}
&(j=2n),
\end{cases}
\end{eqnarray*}
where we understand 
$x_0=\ovl x_0=y_0=\ovl y_0=1$,
$\ovl x_n=x_n$ and $\ovl y_n=y_n$.

Here we define $\ovl\cR(x,y)=
\ovl\cR(\TY(A,2,2n-1))(x,y)=(x',y')$
by
\begin{eqnarray*}
&&\hspace{-20pt}
x'_i\seteq y_i\frac{\ovl V_i}
{\TTY(\ovl V,n,0)}\,\,
,\q
\ovl x'_i\seteq \ovl y_i
\frac{\ovl V_0^\sharp\TTY(\ovl W,n,i)}
{\ovl V_{i}\ovl W_1}\,\,\,
(1\leq i\leq n-1),\q
x'_n=y_n
\frac{{{\STY(\ovl V,\sharp,0)}}^2\ovl W_{n-1}
\ovl W_n}
{\TTY(\ovl V,n,n-1)\TTY(\ovl V,n,n)\ovl W_1^2}
\frac{1+\frac{\ovl W_{n}}{\ovl W_{n-1}}\mu(y)}
{1+\mu(y)},
\\
&&\hspace{-20pt}
y'_i\seteq x_i\frac{{\ovl V^*_i}
\TTY(\ovl W,n,1)}
{\ovl V_{0}^\sharp W_i}
\,\,\,
\ovl y'_i\seteq \ovl x_i\frac{{\ovl V_0}}
{{\ovl V^*_{i}}}\,\,\,
(1\leq i\leq n-1),\q
y'_n=x_n\frac{{\TTY(\ovl V,n,0)}^2}
{{\STY(\ovl V,*,n-1)}{\STY(\ovl V,*,n)}}
\frac{1+\frac{\ovl W_{n-1}}{\ovl W_n}\mu(x)}
{1+\mu(x)}.
\end{eqnarray*}

%%%%%%%%%%%%%%%%%%%%%%%%%%%%%%
\subsection{\Tropical R for $\TY(A,2,2n)$
$(n\geq2)$}

We shall describe 
\tropical R's of type $\TY(A,2,2n)$.
We see the following lemma for 
$\{\TY(X,{n,L},2)\}_{L\in\bbC^\times}$.
\begin{lem}
\label{rx22n}
The birational map $R(\TY(A,2,2n))$ is 
a \tropical R map on the 
$\TY(A,2,2n)$-geometric crystal
$\{\TY(X,{n,L},4)\}_{L}$.
\end{lem}
The proof is the 
same as the one for Lemma \ref{rx1}.

Let us describe the explicit form of tropical 
R on $\cB_L(\TY(A,2,2n))$.
Set $\cR(\TY(A,2,2n))\seteq (\eta^{-1},\eta^{-1})
\circ R(\TY(A,2,2n))\circ(\eta,\eta)$ where
$\eta$ is as in \ref{subsec-a2e}.

Let  
$*:\cB_L(\TY(A,2,2n))\times 
\cB_M(\TY(A,2,2n))
\rightarrow \cB_M(\TY(A,2,2n))\times 
\cB_L(\TY(A,2,2n))$ 
be the involution defined by 
\[
((l_0,l_1,l_2,\ld,\ovl l_2,\ovl l_1),
(m_0,m_1,m_2,\ld,\ovl m_2,\ovl m_1))^*
=((m_0,\ovl m_1,\ovl m_2,\ld,m_2,m_1),
(l_0,\ovl l_1,\ovl l_2,\ld, l_2,l_1))
\]
that is, 
$*:l_0\leftrightarrow m_0,\q
l_i\leftrightarrow \ovl m_i, \q
\ovl l_i\leftrightarrow m_i\,\,
(1\leq i\leq n)$.

Restricting the functions $\TTY(V,n+1,i)$
and $\TTY(W,n+1,i)$ for $\TY(A,2,2n+1)$ to 
$\TY(X,{n,L},4)\times \TY(X,{n,M},4)$, 
we define the rational functions 
$V_i$ $(i=0,1,\ld,n+1)$ and $W_i$ 
($i=1,\ld,n+1$) on 
$\cB_L(\TY(A,2,2n))\times
\cB_M(\TY(A,2,2n))$ $(L,M\in\bbC^\times)$ as 
\begin{eqnarray*}
&& W_i:=V_i{V^*_i}
+(M^2-L^2){V^*_i}+(L^2-M^2)V_i
\q(1\leq i\leq n+1),\\
&&
V_i=
\sum_{j=1}^{2n}(\theta_{i,j}(l,m)
+\theta'_{i,j}(l,m))
+\sum_{j=1}^{2n+2}(\eta_{i,j}(l,m)
+\eta'_{i,j}(l,m)),
\end{eqnarray*}
where 
$L^2=l_n\ovl l_n(l_0^2
\prod_{i=1}^{n-1}l_i\ovl l_i)^2$,
$M^2=m_n\ovl m_n(m_0^2\prod_{i=1}^{n-1}m_i\ovl m_i)^2$,
$\Del\seteq\frac{1+\ovl \mu(l)}{1+\ovl \mu(m)}$, and 
\begin{eqnarray*}
&&\theta_{i,j}(l,m)\\
&&
=\begin{cases}
\displaystyle L^2\prod_{k=j}^{i-1}
\frac{\ovl m_k}{\ovl l_k}
\Del^{\del_{i,n+1}-\del_{j,n+1}}&
1\leq j\leq i\leq n+1,\\
\displaystyle M^2\prod_{k=i}^{j-1}
\frac{\ovl l_k}{\ovl m_k}
\Del^{-\del_{j,n+1}}&
0\leq i<j\leq n+1,\nn \\
\displaystyle LM
\left(\prod_{k=0}^{i-1}
\frac{\ovl m_k}{\ovl l_k}\right)
\left(\prod_{k=0}^{2n-j}\frac{m_k}{l_k}\right)
\Del^{\del_{i,n+1}}
\frac{\ovl l_{2n-j+1}
(m_{2n-j+1}+\ovl m_{2n-j+1})}
{(l_{2n-j+1}+\ovl l_{2n-j+1})\ovl m_{2n-j+1}}
&0\leq i\leq n+1<j.
\end{cases}
\end{eqnarray*}
\[
%\begin{equation}
\hspace{-15pt}\theta'_{i,j}(l,m)
=\begin{cases}
\displaystyle L^2\left(\prod_{k=0}^{i-1}
\frac{\ovl m_k}{\ovl l_k}\right)
\left(\prod_{k=0}^{j-1}\frac{m_k}{l_k}\right)
\Del^{\del_{i,n+1}}
\left(\frac{1+\mu(l)}{1+ \mu(m)}\right)
^{\del_{j,n+1}}
&j\leq n+1,\\
\displaystyle LM
\left(\prod_{k=0}^{i-1}
\frac{\ovl m_k}{\ovl l_k}\right)
\left(\prod_{k=0}^{2n-j}\frac{\ovl l_k}
{\ovl m_k}\right)
\Del^{\del_{i,n+1}}
\frac{m_{2n-j+1}(l_{2n-j+1}+\ovl l_{2n-j+1})}
{l_{2n-j+1}(m_{2n-j+1}+\ovl m_{2n-j+1})}&n+1<j.
\end{cases}
\]
%\end{equation}
\begin{eqnarray*}
&&\hspace{-30pt}\eta_{i,j}(l,m)=
\begin{cases}
\displaystyle L^2\left(\prod_{k=j}^{i-1}
\frac{\ovl m_k}{\ovl l_k}\right)
\left(\frac{\ovl m_{j-1}}{l_{j-1}}\right)
\Del^{\del_{i,n+1}}\mu(l)^{\del_{j,n+1}}
&1\leq j\leq i,\\
\displaystyle M^2\left(\prod_{k=i}^{j-1}
\frac{\ovl l_k}{\ovl m_k}\right)
\left(\frac{\ovl m_{j-1}}{l_{j-1}}\right)
\mu(l)^{\del_{j,n+1}}
&0\leq i<j\leq n+1,\\
\displaystyle LM
\left(\prod_{k=0}^{i-1}
\frac{\ovl m_k}{\ovl l_k}\right)
\left(\prod_{k=0}^{n-1}\frac{m_k}{l_k}\right)
\Del^{\del_{i,n+1}}\frac{1+\mu(m)}{1+\ovl \mu(l)}
& 0\leq i\leq n+1,j=n+2\\
\displaystyle LM
\left(\prod_{k=0}^{i-1}
\frac{\ovl m_k}{\ovl l_k}\right)
\left(\prod_{k=0}^{2n-j+1}
\frac{m_{k}}{l_{k}}\right)
\Del^{\del_{i,n+1}}
\frac{\frac{m_{2n-j+2}}{\ovl m_{2n-j+2}}+1}
{\frac{\ovl l_{2n-j+2}}{l_{2n-j+2}}+1}
&0\leq i<n+2<j\leq  2n+1,\\
\displaystyle LM
\left(\prod_{k=0}^{i-1}
\frac{\ovl m_k}{\ovl l_k} \right)
\Del^{\del_{i,n+1}}&i\leq n+1,j=2n+2.
\end{cases}
\end{eqnarray*}
\[
\eta'_{i,j}(l,m)=
\begin{cases}
\displaystyle L^2\left(\prod_{k=0}^{i-1}
\frac{\ovl m_k}{\ovl l_k}\right)
\left(\prod_{k=0}^{j-1}
\frac{m_k}{l_k}\right)
\left(\frac{l_{j-1}}{\ovl m_{j-1}}\right)
\Del^{\del_{i,n+1}}\ovl \mu(m)^{\del_{j,n+1}}
&1\leq j\leq n+1,\\
\displaystyle L^2\left(\prod_{k=0}^{i-1}
\frac{\ovl m_k}{\ovl l_k}\right)
\left(\prod_{k=0}^n
\frac{m_k}{l_k}\right)
\Del^{\del_{i,n+1}}
\frac{1+\mu(l)}{(1+\mu(m))\mu(m)}
&j=n+2,\\
\displaystyle LM
\left(\prod_{k=0}^{i-1}
\frac{\ovl m_k}{\ovl l_k}\right)
\left(\prod_{k=0}^{2n-j+1}\frac{\ovl l_{k}}
{\ovl m_{k}}\right)
\Del^{\del_{i,n+1}}
\frac{(l_{2n-j+2}+\ovl l_{2n-j+2})
\ovl m_{2n-j+2}}
{(m_{2n-j+2}+\ovl m_{2n-j+2})l_{2n-j+2}}
&n+2<j\leq 2n+1,\\
\displaystyle LM
\left(\prod_{k=0}^{i-1}
\frac{\ovl m_k}{\ovl l_k}\right)
\Del^{\del_{i,n+1}}
&j=2n+2,
\end{cases}
\]
where we understand $\ovl l_0=l_0$ and 
$\ovl m_0=m_0$.

Now, we define the \tropical R 
$\cR(\TY(A,2,2n))\cl\cB_L(\TY(A,2,2n))\times 
\cB_M(\TY(A,2,2n))\to \cB_M(\TY(A,2,2n))
\times \cB_L(\TY(A,2,2n))$ by 
\begin{eqnarray}
&&{\cR(\TY(A,2,2n))}(l,m)=(l',m') 
\q {\rm where}
\nn\\
&&l'_0=
m_0\frac{{V_0}}{V_1}, \q 
l'_i=m_i\frac{V_{i}W_{i+1}}
{V_{i+1}W_{i}},\q
\ovl l'_i
=\ovl m_i\frac{V_{i}}{V_{i+1}}
\q(1\leq i\leq n-1),\nn\\
&&
l'_n=\frac{m_n\TTY(V,n,n)\TTY(W,n,n+1)}
{(1+\mu(m))\TTY(V,n,n+1)\TTY(W,n,n)}
\left(1+\frac{m_n\TTY(W,n,n+1)}
{\ovl m_n\TTY(W,n,n)}\ovl \mu(m)\right),\q
\ovl l'_n=\frac{\ovl m_n\TTY(V,n,n)}
{(1+\ovl \mu(m))\TTY(V,n,n+1)}
\left(1+\frac{\ovl m_n\TTY(W,n,n)}
{m_n\TTY(W,n,n+1)}\mu(m)\right),
%\label{explicit-R-aen}
\nn\\
&&
m'_0=l_0\frac{V_0}{{V^*_1}}, 
\q 
m'_i=l_i\frac{{V^*_{i}}}
{{V^*_{i+1}}},
\q \ovl m'_i=\ovl l_i
\frac{{V^*_{i}}W_{i+1}^{(n)}}
{{V^*_{i+1}}W_{i}^{(n)}}
\q(1\leq i\leq n-1),\nn\\
&&
m'_n=\frac{l_n{\STY(V,*,n)}}
{(1+\mu(l)){\STY(V,*,n+1)}}
\left(1+\frac{l_n\TTY(W,n,n)}
{\ovl l_n\TTY(W,n,n+1)}\ovl \mu(l)\right),\q
\ovl m'_n=
\frac{\ovl l_n{\STY(V,*,n)}\TTY(W,n,n+1)}
{(1+\ovl \mu(l)){\STY(V,*,n+1)}\TTY(W,n,n)}
\left(1+\frac{\ovl l_n\TTY(W,n,n+1)}
{l_n\TTY(W,n,n)}\mu(l)\right),\q
\nn
\end{eqnarray}
Here note that for 
$(l',m')={\mathcal R}(\TY(A,2,2n))(l,m)$ we have 
$(l_1'l_2'\cd l'_{n-1}\ovl l'_{n-1}\cd
 \ovl l_2'\ovl l_1')^2 l'_n\ovl l'_n
=M^2$ and \\ 
$(m_1'm_2'\cd m'_{n-1}\ovl m'_{n-1}\cd
\ovl m_2'\ovl m_1')^2m'_n\ovl m'_n
=L^2$.

Next, we shall describe \tropical R on 
$\cV(\TY(A,2,2n))_L\times \cV(\TY(A,2,2n))_M$. 
Let $*$ be the involution 
$\cV(\TY(A,2,2n))_L\times \cV(\TY(A,2,2n))_M
\to \cV(\TY(A,2,2n))_M\times \cV(\TY(A,2,2n))_L$
defined by $x_0\leftrightarrow y_0$ and 
\begin{eqnarray*}
x_i\mapsto \frac{y_0^2}{M^2 \ovl y_i},\,\,
\ovl x_i\mapsto \frac{y_0^2}{M^2 y_i},\,\,
x_n\mapsto\frac{y_0^4}{M^2 y_n},\,\,
y_i\mapsto \frac{x_0^2}{L^2 \ovl x_i},\,\,
\ovl y_i\mapsto \frac{x_0^2}{L^2 x_i}\,\,
y_n\mapsto\frac{x_0^4}{L^2 x_n},\,\,
(i=1,\cd,n-1).
\end{eqnarray*}
We define the rational functions
$\ovl V_i(x,y)$ $(i=0,1,\ld,n+1)$ 
and 
$\ovl W_i(x,y)$ ($i=1,\ld,n+1$) on 
$\cV(\TY(A,2,2n))_L\times
\cV(\TY(A,2,2n))_M$ $(L,M\in\bbC^\times)$ by 
\begin{eqnarray*}
&&\ovl V_i=
\sum_{j=1}^{2n}(\ovl\theta_{i,j}(x,y)
+\ovl\theta'_{i,j}(x,y))
+\sum_{j=1}^{2n+2}(\ovl\eta_{i,j}(x,y)
+\ovl\eta'_{i,j}(x,y)),\\
&&
\ovl W_i:=\ovl V_i{\ovl V*_i}
+(M^2-L^2){\ovl V^*_i}
+(L^2-M^2)\ovl V_i
\q(1\leq i\leq n+1),\\
&&
\mu(x)\seteq\frac{x_n}{Lx_{n-1}\ovl x_{n-1}}=
\ovl\mu(x)^{-1},\,\,
\mu(y)\seteq\frac{y_n}{My_{n-1}\ovl y_{n-1}}=
\ovl\mu(y)^{-1},\\
&&
\Del\seteq\frac{1+\ovl\mu(x)}{1+\ovl\mu(y)},\q
\nabla\seteq\frac{1+\mu(x)}{1+\mu(y)},\q
\square\seteq\frac{1+\mu(y)}{1+\ovl\mu(x)}.
%\Del\seteq\frac{1+\frac{Lx_{n-1}\ovl x_{n-1}}
%{x_n}}{1+\frac{M y_{n-1}\ovl y_{n-1}}{y_n}},\q
%\nabla\seteq\frac{1+\frac{x_n}
%{Lx_{n-1}\ovl x_{n-1}}}
%{1+\frac{y_n}{M y_{n-1}\ovl y_{n-1}}},\q
%\square\seteq
%\frac{1+\frac{y_n}{M y_{n-1}\ovl y_{n-1}}}
%{1+\frac{Lx_{n-1}\ovl x_{n-1}}{x_n}}.
\end{eqnarray*}
\begin{eqnarray*}
&&\hspace{-30pt}\ovl\theta_{i,j}(x,y)\seteq\\
&&\hspace{-30pt}\begin{cases}\displaystyle
L^2\left(\frac{x_{i-1}}{y_{i-1}}
\right)^{1+\del_{i,1}}
\left(\frac{y_{j-1}}
{x_{j-1}}\right)^{1+\del_{j,1}}
\left(\frac{y_{n-1}}
{x_{n-1}}\Del\right)^{\del_{i,n+1}-\del_{j,n+1}}
&(1\leq j\leq i\leq n+1),\\ \displaystyle
M^2\frac{x_{i-1}y_{j-1}}{x_{j-1}y_{i-1}}
\left(\frac{x_{n-1}}{y_{n-1}}
\Del^{-1}\right)^{\del_{j,n+1}}
\left(\frac{x_0}{y_0}\right)^{
\del_{i,0}+\del_{i,1}-\del_{j,1}}
&(0\leq i<j\leq n+1),\\ \displaystyle
LM\frac{x_{i-1}}{y_{i-1}}\left(
\frac{M^2 \ovl y_{2n-j}}{L^2 \ovl x_{2n-j}}
\right)^{1-\del_{j,2n}}
\frac{1+M^{2\del_{j,2n}}
\frac{y_{2n+1-j}\ovl y_{2n+1-j}}
{y_{2n-j}\ovl y_{2n-j}}}
{1+L^{2\del_{j,2n}}
\frac{x_{2n+1-j}\ovl x_{2n+1-j}}
{x_{2n-j}\ovl x_{2n-j}}}
\left(\frac{x_0}{y_0}\right)^{
\del_{i,0}+\del_{i,1}}
\left(\frac{y_{n-1}}{x_{n-1}}
\Del\right)^{\del_{i,n+1}}
&(i\leq n+1<j<2n),
\end{cases}
\end{eqnarray*}
\begin{eqnarray*}
&&\hspace{-30pt}
\ovl\theta'_{i,j}(x,y)\seteq\\
&&\hspace{-30pt}\begin{cases}\displaystyle
L^2
\frac{x_{i-1}}{y_{i-1}}
\left(\frac{y_{n-1}}{x_{n-1}}\Del
\right)^{\del_{i,n+1}}
\left(\frac{M^2 \ovl y_{j-1}}
{L^2 \ovl x_{j-1}}\right)^{1-\del_{j,1}}
\left(\frac{L^2 \ovl x_{n-1}}{M^2 \ovl y_{n-1}}
\nabla\right)^{\del_{j,n+1}}
\left(\frac{x_0}{y_0}\right)^{
\del_{i,0}+\del_{i,1}}
&(1\leq j\leq n+1),\\ \displaystyle
LM
\frac{x_{i-1}y_{2n-j}}
{y_{i-1}x_{2n-j}}
\frac{1+\frac{x_{2n-j}\ovl x_{2n-j}}
{L^{2\del_{j,2n}}x_{2n+1-j}\ovl x_{2n+1-j}}}
{1+\frac{y_{2n-j}\ovl y_{2n-j}}
{M^{2\del_{j,2n}}y_{2n+1-j}\ovl y_{2n+1-j}}}
\left(\frac{y_n}{x_n}\Del\right)^{\del_{i,n+1}}
\left(\frac{y_0}{x_0}\right)^{
\del_{i,0}+\del_{i,1}-\del_{j,2n}}
&(n+1<j\leq 2n),\\ 
\end{cases}
\end{eqnarray*}
\begin{eqnarray*}
&&\hspace{-30pt}\ovl\eta_{i,j}(x,y)\seteq\\
&&\hspace{-30pt}\begin{cases}\displaystyle
L^{2+\del_{j,n+1}-2\del_{j,2}}
\frac{x_{i-1}y_{j-2}\ovl x_{j-2}}
{y_{i-1}x_{j-1}\ovl x_{j-1}}
\left(\frac{y_{n-1}}{x_{n-1}}
\Del\right)^{\del_{i,n+1}}
\frac{x_0^{\del_{i,1}+\del_{j,1}-\del_{j,2}}}
{y_0^{\del_{i,1}-\del_{j,1}-\del_{j,2}}}
x_n^{\del_{j,n+1}}
&(1\leq j\leq i\leq n+1),\\ \displaystyle
M^2L^{\del_{j,n+1}-2\del_{j,2}}
\frac{x_{i-1}y_{j-2}\ovl x_{j-2}}
{y_{i-1}x_{j-1}\ovl x_{j-1}}
\frac{x_0^{\del_{i,0}
+\del_{i,1}+\del_{j,1}-\del_{j,2}}}
{y_0^{\del_{i,0}+
\del_{i,1}-\del_{j,1}-\del_{j,2}}}
x_n^{\del_{j,n+1}}
&(0\leq i< j\leq n+1),\\ \displaystyle
\frac{M^3}{L}
\frac{x_{i-1}\ovl y_{n-1}}
{y_{i-1}\ovl x_{n-1}}
\left(\frac{y_{n-1}}{x_{n-1}}
\Del\right)^{\del_{i,n+1}}
\left(\frac{x_0}{y_0}\right)^{
\del_{i,0}+\del_{i,1}}\square
&(0\leq i<j=n+2),\\ \displaystyle
\frac{M^3}{L}
\frac{x_{i-1}\ovl y_{2n-j+1}}
{y_{i-1}\ovl x_{2n-j+1}}
\frac{1+\frac{y_{2n-j+2}
\ovl y_{2n-j+2}}{y_{2n-j+1}
\ovl y_{2n-j+1}}}
{1+\frac{x_{2n-j+1}
\ovl x_{2n-j+1}}
{x_{2n-j+2}\ovl x_{2n-j+2}}}
\left(\frac{y_{n-1}}{x_{n-1}}
\Del\right)^{\del_{i,n+1}}
\left(\frac{x_0}{y_0}\right)^{
\del_{i,0}+\del_{i,1}}
&(i\leq n+1,n+2<j\leq 2n)\\ \displaystyle
LM\frac{x_{i-1}}{y_{i-1}}
\frac{1+\frac{M^2y_1\ovl y_1}{y_0^2}}
{1+\frac{x_0^2}{L^2x_1\ovl x_1}}
\left(\frac{y_{n-1}}{x_{n-1}}\Del
\right)^{\del_{i,n+1}}
\left(\frac{x_0}{y_0}\right)^{
\del_{i,0}+\del_{i,1}} 
&(i\leq n+1,j=2n+1),\\ \displaystyle
LM\frac{x_{i-1}}{y_{i-1}}
\left(\frac{y_{n-1}}
{x_{n-1}}\Del\right)^{\del_{i,n+1}}
\left(\frac{x_0}{y_0}\right)^{
-1+\del_{i,0}+\del_{i,1}}
&(i\leq n+1,j=2n+2),
\end{cases}
\end{eqnarray*}
\begin{eqnarray*}
&&\hspace{-30pt}
\ovl\eta'_{i,j}(x,y)\seteq\\
&&\hspace{-30pt}\begin{cases}\displaystyle
\ovl\theta'_{i,j}(x,y)\frac{\ovl x_{j-1}y_{j-1}}
{\ovl x_{j-2}y_{j-2}}
\frac{L^{2\del_{j,2}}x_0^{\del_{j,2}}}
{L^{2\del_{j,n+1}}x_0^{2\del_{j,1}}
y_0^{\del_{j,2}}}
\left(\frac{M \ovl y_{n-1}}{x_{n-1}y_n\nabla}
\right)^{\del_{j,n+1}}
&(1\leq j\leq n+1),\\ \displaystyle
L^2M\frac{x_{i-1}y_{n-1}\ovl x_{n-1}}
{y_{i-1}x_n}
\left(\frac{y_{n-1}}
{x_{n-1}}\Del\right)^{\del_{i,n+1}}
\left(\frac{x_0}{y_0}\right)^{
\del_{i,0}+\del_{i,1}}\nabla
&(j=n+2)\\ \displaystyle
LM\frac{x_{i-1} y_{2n-j+1}}
{y_{i-1}x_{2n-j+1}}
\frac{\left(1+\frac{x_{2n-j+1}
\ovl x_{2n-j+1}}
{x_{2n-j+2}\ovl x_{2n-j+2}}\right)}
{\left(1+\frac{
y_{2n-j+2}\ovl y_{2n-j+2}}
{y_{2n-j+1}\ovl y_{2n-j+1}}\right)}
\left(\frac{y_{n-1}}{x_{n-1}}\Del
\right)^{\del_{i,n+1}}
\left(\frac{x_0}{y_0}\right)^{
\del_{i,0}+\del_{i,1}}
&(n+1<j\leq 2n),\\ \displaystyle
LM\frac{x_{i-1}}{y_{i-1}}
\frac{1+\frac{x_0^2}
{L^2x_1\ovl x_1}}
{1+\frac{M^2 y_1\ovl y_1}{y_0^2}}
\left(\frac{y_{n-1}}{x_{n-1}}\Del
\right)^{\del_{i,n+1}}
\left(\frac{x_0}{y_0}\right)^{
-1+\del_{i,0}+\del_{i,1}}
&(j=2n+1),\\ \displaystyle
LM\frac{x_{i-1}}{y_{i-1}}
\left(\frac{y_{n-1}}{x_{n-1}}\Del
\right)^{\del_{i,n+1}}
\left(\frac{x_0}{y_0}\right)^{
-1+\del_{i,0}+\del_{i,1}}
&(j=2n+2),
\end{cases}
\end{eqnarray*}
where we understand that $\ovl x_0=x_0$, 
$\ovl y_0=y_0$ and $x_{-1}=\ovl x_{-1}=
y_{-1}=\ovl y_{-1}=1$.
Here we define $\ovl\cR(x,y)=
\ovl\cR(\TY(A,2,2n))(x,y)=(x',y')$
by
\begin{eqnarray*}
&&
x'_0\seteq y_0\frac{\TTY(\ovl V,n,1)}
{\TTY(\ovl V,n,0)},\q
x'_i\seteq y_i\frac{\TTY(\ovl V,n,1)
\ovl V_{i+1}}
{{\TTY(\ovl V,n,0)}^2},\q
\ovl x'_i\seteq \ovl y_i
\frac{{\ovl V_1}\TTY(\ovl W,n,i+1)}
{\ovl V_{i+1}\ovl W_1}\,\,\,
(1\leq i\leq n-1),\\
&&x'_n=y_n\frac{{\TTY(\ovl V,n,1)}^2
\TTY(\ovl W,n,n)\TTY(\ovl W,n,n+1)}
{{\TTY(\ovl V,n,n)}\TTY(\ovl V,n,n+1)
{\TTY(\ovl W,n,1)}^2}
\left(1+\frac{\mu(y)\TTY(\ovl W,n,n+1)}
{\TTY(\ovl W,n,n)}
\right)(1+\mu(y))^{-1},
\\
&&
y'_0\seteq x_0\frac{{\STY(\ovl V,*,1)}}
{\TTY(\ovl V,n,0)},\q
y'_i\seteq x_i\frac{{\STY(\ovl V,*,1)}
{\ovl V^*_{i+1}}\TTY(\ovl W,n,1)}
{{\ovl V_{0}}^2 W_{i+1}},\q
\ovl y'_i\seteq \ovl x_i
\frac{{\ovl V^*_1}}
{{\ovl V^*_{i+1}}}\,\,\,
(1\leq i\leq n-1),\\
&&y'_n=x_n\frac{({\STY(\ovl V,*,1)})^2}
{{\STY(\ovl V,*,n)}{\STY(\ovl V,*,n+1)}}
\left(1+\frac{\mu(x)
\TTY(\ovl W,n,n)}{
\TTY(\ovl W,n,n+1)}\right)(1+\mu(x))^{-1}
\end{eqnarray*}
%%%%%%%%%%%%%%%%%%%%%%%%%%%%%%%%%%%%%%%%%%%%%
\subsection{Uniqueness of the \tropical R-maps}

%We shall show the following theorem in Sect.5.
\begin{thm}
\label{uniqueness}
Let $\ovl\cR$ 
be the \tropical R as introduced in
this section.
Set $z_0\seteq\cR(\bf 1,1)$. 
Let $\cR'$ be a tropical
R such that $\cR'(\bf 1,1)=z_0$. Then we have 
$\ovl\cR=\cR'$ as birational maps.
\end{thm}
{\sl Proof.}
Let $\cV_l,\cV_m$ $(l,m\in\bbC^\times)$ 
be the 
affine geometric crystals constructed in 
Sect.\ref{sec5}.
By Theorem \ref{ud-geo} and 
Theorem \ref{unip}, we have that 
${\mathcal UD}_{(\theta_l,\theta_m)}
(\cV_l(\ge)\times \cV_m(\ge))$
is isomorphic to the crystal $B_\ify(\ge^L)
\ot B_\ify(\ge^L)$, where 
$\theta_l$ is the positive structure as in
 \ref{subsec-ud}.
Since $B_\ify(\ge^L)\ot B_\ify(\ge^L)$ is connected, 
$\cV_l(\ge)\times \cV_m(\ge)$ is prehomogeneous
by Theorem \ref{conn-preh}.
Therefore, by Lemma \ref{uniq} we obtain 
$\ovl\cR=\cR'$, which completes the proof of 
Theorem \ref{uniqueness}. \qed

%%%%%%%%%%%%%%%%%%%%%%%%%%%%%%%%%%%%%%%%%%%

\end{document}